\documentclass[11pt, makeidx]{amsart}

\usepackage{enumerate}
\usepackage{amssymb,amsfonts,amsthm}
\usepackage{srcltx}
\usepackage{float}
\usepackage[OT2,OT1]{fontenc}
\newcommand\cyr{%
\renewcommand\rmdefault{wncyr}%
\renewcommand\sfdefault{wncyss}%
\renewcommand\encodingdefault{OT2}%
\normalfont
\selectfont}
\DeclareTextFontCommand{\textcyr}{\cyr}
\usepackage{color}
\makeindex

\makeatletter
\def\sssub{\@startsection{paragraph}{4}}

\renewcommand\paragraph{\@startsection{paragraph}{4}{\z@}{1.25ex}{0.0001pt}{\normalfont\normalsize\em}}

\makeatother

\numberwithin{paragraph}{subsubsection}

\newcommand{\N}{\mathcal{N}}
\newcommand{\Z}{{\mathbb Z}}

\newcommand{\e }{\varepsilon }

\newcommand{\la}{\langle}
\newcommand{\ra}{\rangle}

\newcommand{\bbb}{{\mathcal B}}

\newtheorem{theorem}{Theorem}[section]

\newtheorem{lemma}[theorem]{Lemma}

\newtheorem{prop}[theorem]{Proposition}
\theoremstyle{definition}
\newtheorem{df}[theorem]{Definition}

\newtheorem{rk}[theorem]{Remark}
\newtheorem{prob}[theorem]{Problem}
\newtheorem{notation}[theorem]{Notation}
\newtheorem{assumption}[theorem]{Assumption}

\newcommand{\Lab}{{\mathrm{Lab}}}

\newcommand{\gcal}{{\mathcal G}}

\newcommand{\tool}{\stackrel{\ell}{\too} }

\newcommand{\ttt}{{\mathcal T}}
\newcommand{\kkk}{{\mathcal K}}

\newcommand{\aaa}{{\mathcal A}}

\newcommand{\gi}{\Gamma}

\newcommand{\vk}{van Kampen }

\newcommand{\iv}{^{-1}}

\newcommand{\too}{\to }

\newcommand{\xxx}{{\mathcal X} }

\newcommand{\pp}{{\mathcal P} }

\newcommand{\sss}{{\mathcal S} }

\newcommand{\ncn}{{\ll\!\!\N\!\!\gg}}

\newcommand{\da}{\gcal}
\newcommand{\p}{\mathfrak{p}}
\newcommand{\q}{\mathfrak{q}}
\newcommand{\rf}{\mathfrak{r}}
\newcommand{\s}{\mathfrak{s}}

\numberwithin{equation}{section}

\begin{document}

\def\thesubsubsection       {\thesubsection.\Alph{subsubsection}}

\def\theparagraph       {\thesubsubsection\arabic{paragraph}}

\title[Aspherical groups and manifolds]{Aspherical groups and manifolds with extreme properties}
\author{Mark Sapir}\thanks{The research was supported in part by NSF grant DMS-0700811.}
\address{Department of Mathematics, Vanderbilt University, Nashville, TN 37240, U.S.A.}
\email{m.sapir@vanderbilt.edu}

\begin{abstract}
We prove that every finitely generated group with recursive
aspherical presentation embeds into a group with finite aspherical presentation. This and several known facts  about
groups and manifolds imply that there exists a 4-dimensional closed
aspherical manifold $M$ such that the fundamental group $\pi_1(M)$ coarsely contains an expander, and so it has infinite asymptotic dimension, is not coarsely embeddable into a Hilbert
space, does not satisfy G. Yu's property A, and does not satisfy the Baum-Connes conjecture with
coefficients. Closed aspherical manifolds with any of these properties were previously unknown.
\end{abstract}
\maketitle

\tableofcontents

\section{Introduction}

The main result of this paper is the following theorem (see Theorem \ref{th:1}). For definitions see Section \ref{ss:aafc}.

\begin{theorem}\label{main} Every finitely generated group with combinatorially aspherical recursive presentation complex embeds into a group with finite combinatorially aspherical presentation complex.
\end{theorem}

Using Davis' construction \cite{Davis} (see also \cite[Chapter 11]{Davis1}) this allows one to create closed aspherical manifolds of dimension 4 and higher with some previously unknown ``extreme" properties.\footnote{As explained in \cite[Chapter 11]{Davis1}, in dimension 5 and higher, one can assume that the aspherical manifold we construct is smooth. This was communicated to us by I. Belegradek.} For example, by Gromov \cite{GrRandom} (for more details see Arzhantseva and Delzant \cite{AD}), there exists a finitely generated group with recursive combinatorially aspherical presentation whose Cayley graph coarsely contains an expander. Hence Theorem \ref{main} imply that there exist closed aspherical manifolds of dimension 4 and higher whose fundamental groups coarsely contain expanders. These groups and manifolds are not coarsely embeddable into a Hilbert space \cite{GrRandom}, do not satisfy G.Yu's property A \cite{GYu}, and are counterexamples to the Baum-Connes conjecture with coefficients (see Higson, Lafforgue, Skandalis \cite{HLS}). They also have infinite asymptotic dimension.  That solves a problem, formulated first by G. Yu in \cite{GYu1} asking whether the fundamental group of a closed aspherical manifold can have infinite asymptotic dimension. A weaker problem of whether the asymptotic dimension of a closed aspherical manifold can exceed its (ordinary) dimension was mentioned by Gromov in \cite[Page 33]{Grai} and was open till now also. Note that Dranishnikov's problem whether the asymptotic dimension of an aspherical $n$-manifold is always $n$ or infinity is still open (see \cite[Problem 3.4]{Dra}).  It seems that  Gromov's random groups and our Theorem \ref{main} cannot give an example where the dimensions are different while both are finite in view of the recent paper by Willett \cite{Wil}.\footnote{For a recent survey about aspherical manifolds see L\"uck \cite{Luck}, for constructions of aspherical manifolds with other ``exotic" properties see Davis \cite{Davis1}.}

As another corollary one can deduce that a torsion-free Tarski monster (that is a finitely generated group all of whose proper subgroups are infinite cyclic, see Olshanskii \cite{Olbook}) embeds into the fundamental group  of a closed aspherical manifold. Indeed, by \cite{Olbook}, the torsion-free Tarski monsters constructed by Olshanskii have recursive combinatorially aspherical presentations, it remains to apply Theorem \ref{main} and \cite{Davis}. More generally, every {\em lacunary hyperbolic} group  given by a recursive {\em graded small cancelation} presentation (see \cite{OOS}) embeds into the fundamental group of a closed aspherical manifold.

Note that the fact that  any recursively presented group embeds into a finitely presented group is the celebrated Higman embedding theorem \cite{Hi} which is one of the main results in the algorithmic group theory. There exist at least ten versions of Higman's embedding construction preserving various properties of groups, from relatively easy constructions to very complicated (see \cite{Rot, BDM, BORS, OScol, SaICM, Sasur} and references there). But all the previous versions could not produce finitely presented groups with finite $K(.,1)$. The reason is also common to all these constructions and can be roughly described as follows. Let $G$ be a finitely presented group containing a copy of a recursively presented group $\Gamma$ and constructed using one of the existing proofs of the Higman embedding theorem. Then there exists $1\ne g\in G$ which centralizes $\Gamma$ (it is true for each of the previously known constructions). Consider two copies $\Delta_1, \Delta_2$ of a \vk diagram  over the presentation of $G$ with boundaries labeled by a relation $r$ of $\Gamma$ (such diagrams must exist by the \vk lemma since $r=1$ in $G$). Consider the \vk diagram $\Delta_3$ for the commutativity relation $gr=rg$ obtained by gluing together diagrams for $ga=ag$ for every letter $a$ in $r$ (see Figure \ref{com}). It has the form of a rectangle with sides labeled by $r, g, r\iv, g\iv$. Glue $\Delta_1, \Delta_2, \Delta_3$ together to form a spherical \vk diagram over the presentation of $G$. That spherical diagram is not combinatorially homotopic to a trivial diagram. Moreover it is basically clear that if the defining relations of $\Gamma$ are independent enough (say, satisfy some form of  small cancelation as in many interesting cases), then these spherical diagrams cannot be ``generated" by finitely many spherical diagrams over the presentation of $G$. Hence the relation module \cite{LS} of the group $G$ is infinitely generated. A similar problem arises even if we assume only $g\Gamma g\iv\subset \Gamma$. The main difficulty that we had to overcome in this paper was to avoid this type of ``trivial" spherical diagrams. In our construction, we are using $S$-machines (which can be viewed as multiple HNN-extensions of free groups) first introduced in \cite{SBR} and used for some versions of Higman embedding in \cite{BORS, OScol} and other papers. The finitely presented group in this paper is built from two (different) $S$-machines, and several hyperbolic and close to hyperbolic groups that ``glue" these $S$-machines together. One of the main tools of the proof is the congruence extension property of certain subgroups of hyperbolic groups first established by  Olshanskii in \cite{Ol95}.

\begin{figure}
\unitlength .6 mm 
\linethickness{0.4pt}
\ifx\plotpoint\undefined\newsavebox{\plotpoint}\fi 
\begin{picture}(119,26.25)(0,0)
\put(10.25,14.625){\vector(0,1){.07}}\put(10.25,6.25){\line(0,1){16.75}}
\put(19.25,14.625){\vector(0,1){.07}}\put(19.25,6.25){\line(0,1){16.75}}
\put(28.25,14.625){\vector(0,1){.07}}\put(28.25,6.25){\line(0,1){16.75}}
\put(37.25,14.625){\vector(0,1){.07}}\put(37.25,6.25){\line(0,1){16.75}}
\put(46.25,14.625){\vector(0,1){.07}}\put(46.25,6.25){\line(0,1){16.75}}
\put(55,14.625){\vector(0,1){.07}}\put(55,6.25){\line(0,1){16.75}}
\put(63.125,23.125){\vector(1,0){.07}}\multiput(10.25,23)(13.21875,.03125){8}{\line(1,0){13.21875}}
\put(62.875,6.5){\vector(1,0){.07}}\multiput(10.25,6.25)(7.0166667,.0333333){15}{\line(1,0){7.0166667}}
\put(115.625,15){\vector(0,1){.07}}\multiput(115.5,6.75)(.03125,2.0625){8}{\line(0,1){2.0625}}
\put(7.75,14.25){\makebox(0,0)[cc]{$g$}}
\put(17,14.25){\makebox(0,0)[cc]{$g$}}
\put(26.75,14.25){\makebox(0,0)[cc]{$g$}}
\put(35,14.25){\makebox(0,0)[cc]{$g$}}
\put(45,14.25){\makebox(0,0)[cc]{$g$}}
\put(53,14.25){\makebox(0,0)[cc]{$g$}}
\put(55.25,26.25){\makebox(0,0)[cc]{$r$}}
\put(54,3){\makebox(0,0)[cc]{$r$}}
\put(119,14.25){\makebox(0,0)[cc]{$g$}}
\put(84,14.25){\makebox(0,0)[cc]{$\ldots$}}
\end{picture}

\caption{}\label{com}
\end{figure}

One can view this paper as giving examples of  finitely presented groups with extreme properties. There are now many examples of finitely generated extreme groups, sometimes called monsters, but finitely presented monsters are much more rare. There are some objective reasons for it. For example a hyperbolic group version of the Cartan-Hadamard theorem proved by Gromov  (see the Appendix of \cite{OOS} by Kapovich and Kleiner) shows that a lacunary hyperbolic finitely presented group is hyperbolic and so it cannot be torsion, have few subgroups, or be of intermediate growth, etc. Since many methods of creating extreme groups are based on lacunar hyperbolicity \cite{Olbook,OOS}, these methods cannot produce finitely presented groups. Among many outstanding problems of that kind, we mention the following four here.

\begin{prob}Is there a finitely presented infinite torsion group?
\end{prob}
There is an idea due to Rips for constructing such an example: Embed the infinite free Burnside group $B(m,n)$ \cite{Ad} into a finitely presented group $H=\la X\mid R\ra$, then impose relations $x=w_x$ where $w_x$ are ``random" words in $B(m,n)$. Clearly the factor-group is torsion (being also a factor-group of $B(m,n)$). The question is whether it is infinite. So far there was no good enough embedding of $B(m,n)$ into a finitely presented group. It may be possible that the embedding from this paper helps. Note that if we replace the equalities $x=w_x$ by conjugacy $x^t=w_x$ ($t$ is a new generator), then the method can work and produce a finitely presented group containing a copy of $B(2,n)$ (provided $n\gg1$ and odd) and with bounded torsion derived subgroup. This way a finitely presented counterexample to the von Neumann conjecture was constructed in \cite{OSamen}.
\begin{prob}
Is there a finitely presented Tarski monster?
\end{prob}
Currently there are no ideas how to construct such a group and
no ideas how to prove that such a group does not exist.

\begin{prob}
Is there an infinite finitely presented group all of whose non-trivial elements are conjugate?
\end{prob}
Osin \cite{Osin} has recently constructed an infinitely presented finitely generated example (that was a major breakthrough since the problem was open for more than 60 years). As with Tarski monsters, there are no approaches to proving or disproving that a finitely presented example exists.

\begin{prob}
Is there a finitely presented group of intermediate growth?
\end{prob}

All known examples of groups of intermediate growth are close to automata groups which are in a sense relatives of the first Grigorchuk group \cite{Grig}. There are no general ideas how to construct a finitely presented example except by checking more and more groups generated by automata in hope to find an example. There are also no ideas of how to prove that such groups do not exist.  There are several embeddings of Grigorchuk's group into finitely presented groups \cite{GrigL}, \cite{Rov}. All these finitely presented groups have exponential growth. A version of Higman embedding construction preserving or almost preserving growth currently is completely out of reach because all known constructions use HNN extensions and/or amalgamated products.

{\bf Acknowledgement.} This paper completes a project started about ten years ago after M. Gromov first announced the construction of his random finitely generated groups containing expanders, and a question whether these groups embed into groups with finite $K(.,1)$ arose. During these ten years, I was helped by several people including D. Osin and E. Rips. I am especially grateful to A.Yu. Olshanskii for many useful discussions.

\section{Preliminaries}

\subsection{Bands and annuli}
\vskip 0.1 in
The next definition of a band in a diagram is crucial for our paper.
\vskip 0.1 in

Let $\la\xxx\mid R\ra$ be a group presentation.

Let $S$ be a subset of ${\xxx}$. An \label{sband}$S$-band $\bbb$
is a sequence of cells $\pi_1,...,\pi_n$ in a \vk diagram such
that

\begin{itemize}
\item Each two consecutive cells in this sequence have a common edge
labeled by a letter from $S$.
\item Each cell $\pi_i$, $i=1,...,n$ has exactly two $S$-edges
(i.e. edges labeled by a letter from $S$) having opposite orientations.
\end{itemize}

Figure \ref{f-3} illustrates this
concept. In this Figure edges $e, e_1,...,e_{n-1},f$ are
$S$-edges, the lines $l(\pi_i,e_i), l(\pi_i,e_{i-1})$ connect
fixed points in the cells with fixed points of the corresponding
edges.

\begin{figure}[ht]
\unitlength=0.90mm \special{em:linewidth 0.4pt}
\linethickness{0.4pt}
\begin{picture}(149.67,30.11)(5,0)
\put(19.33,30.11){\line(1,0){67.00}}
\put(106.33,30.11){\line(1,0){36.00}}
\put(142.33,13.11){\line(-1,0){35.67}}
\put(86.33,13.11){\line(-1,0){67.00}}
\put(33.00,21.11){\line(1,0){50.00}}
\put(110.00,20.78){\line(1,0){19.33}}
\put(30.00,8.78){\vector(1,1){10.33}}
\put(52.33,8.44){\vector(0,1){10.00}}
\put(76.33,8.11){\vector(-1,1){10.33}}
\put(105.66,7.78){\vector(1,2){5.33}}
\put(132.66,8.11){\vector(-1,1){10.00}}
\put(16.00,21.11){\makebox(0,0)[cc]{$e$}}
\put(29.66,25.44){\makebox(0,0)[cc]{$\pi_1$}}
\put(60.33,25.44){\makebox(0,0)[cc]{$\pi_2$}}
\put(133.00,25.78){\makebox(0,0)[cc]{$\pi_n$}}
\put(145.33,21.44){\makebox(0,0)[cc]{$f$}}
\put(77.00,25.11){\makebox(0,0)[cc]{$e_2$}}
\put(122.00,25.44){\makebox(0,0)[cc]{$e_{n-1}$}}
\put(100.33,25.44){\makebox(0,0)[cc]{$S$}}
\put(96.33,30.11){\makebox(0,0)[cc]{$\dots$}}
\put(96.33,13.11){\makebox(0,0)[cc]{$\dots$}}
\put(26.66,4.78){\makebox(0,0)[cc]{$l(\pi_1,e_1)$}}
\put(52.66,4.11){\makebox(0,0)[cc]{$l(\pi_2,e_1)$}}
\put(78.33,4.11){\makebox(0,0)[cc]{$l(\pi_2,e_2)$}}
\put(104.66,3.78){\makebox(0,0)[cc]{$l(\pi_{n-1},e_{n-1})$}}
\put(134.00,3.11){\makebox(0,0)[cc]{$l(\pi_n,e_{n-1})$}}
\put(88.66,30.11){\makebox(0,0)[cc]{$q_2$}}
\put(88.66,13.11){\makebox(0,0)[cc]{$q_1$}}
\put(49.33,25.11){\makebox(0,0)[cc]{$e_1$}}
\put(33.33,21.11){\circle*{0.94}}
\put(46.66,21.11){\circle*{1.33}}
\put(60.33,21.11){\circle*{0.94}}
\put(73.66,21.11){\circle*{1.33}}
\put(95.00,20.56){\makebox(0,0)[cc]{...}}
\put(19.33,21.00){\line(1,0){14.00}}
\put(129.67,20.67){\line(1,0){12.67}}
\put(19.00,13.00){\vector(0,1){17.00}}
\put(46.67,13.00){\vector(0,1){17.00}}
\put(73.67,13.00){\vector(0,1){17.00}}
\put(117.00,13.00){\vector(0,1){17.00}}
\put(117.00,20.67){\circle*{1.33}}
\put(129.67,20.67){\circle*{0.67}}
\put(142.33,13.00){\vector(0,1){17.00}}
\put(142.33,20.67){\circle*{1.33}}
\put(19.00,21.00){\circle*{1.33}}
\put(14.00,8.33){\vector(1,1){11.00}}
\put(10.33,5.00){\makebox(0,0)[cc]{$l(\pi_1,e)$}}
\put(148.00,8.00){\vector(-1,1){10.33}}
\put(154.67,3.11){\makebox(0,0)[cc]{$l(\pi_n,f)$}}
\end{picture}
\caption{}\label{f-3}
\end{figure}

The broken line formed by the lines $l(\pi_i,e_i)$,
$l(\pi_i,e_{i-1})$ connecting points inside neighboring cells is
called the \label{cl}{\em median} of the band $\bbb$. The $S$-edges
$e$ and $f$ are called the \label{start}{\em start} and {\em end}
edges of the band. The boundary of the subdiagram $\cup\pi_i$ has the form $e\p f\iv\q\iv$. The paths $\p, \q$ are called the {\em sides} of the band.

A band $\pi_1,...,\pi_t$ is called \label{redd}{\em reduced} if
$\pi_{i+1}$ is not a mirror image of $\pi_i$, $i=1,...,t-1$
(otherwise cells $\pi_i$ and $\pi_{i+1}$ {\em cancel} and there
exists a diagram with the same boundary label as $\cup_i\pi_i$ and
containing fewer cells).

We say that two bands {\em intersect} if their medians
intersect. We say that a band is an {\em annulus} if
its median is a closed  curve (see
Figure \ref{f-2} a).

\begin{center}
\begin{figure}[ht]
\unitlength=1.5mm
\linethickness{0.4pt}
\begin{picture}(101.44,22.89)(5,0)
\put(30.78,13.78){\oval(25.33,8.44)[]}
\put(30.78,13.89){\oval(34.67,18.00)[]}
\put(39.67,9.56){\line(0,-1){4.67}}
\put(25.67,9.56){\line(0,-1){4.67}}
\put(18.56,14.89){\line(-1,0){5.11}}
\put(32.56,7.11){\makebox(0,0)[cc]{$\pi_1=\pi_n$}}
\put(21.00,7.33){\makebox(0,0)[cc]{$\pi_2$}}
\put(15.89,11.78){\makebox(0,0)[cc]{$\pi_3$}}
\put(43.44,12.89){\line(1,0){4.67}}
\put(43.00,8.44){\makebox(0,0)[cc]{$\pi_{n-1}$}}
\put(19.22,10.89){\line(-1,-1){4.00}}
\put(25.89,20.44){\circle*{0.00}}
\put(30.78,20.44){\circle*{0.00}}
\put(30.33,1.56){\makebox(0,0)[cc]{a}}
\put(35.44,20.44){\circle*{0.00}}
\put(62.78,4.89){\line(0,1){4.67}}
\put(62.78,9.56){\line(1,0){38.67}}
\put(101.44,9.56){\line(0,-1){4.67}}
\put(101.44,4.89){\line(-1,0){38.67}}
\put(82.11,16.22){\oval(38.67,13.33)[t]}
\put(62.78,15.78){\line(0,-1){7.33}}
\put(101.44,16.44){\line(0,-1){8.22}}
\put(82.00,12.22){\oval(27.78,14.67)[t]}
\put(67.89,12.89){\line(0,-1){8.00}}
\put(95.89,13.11){\line(0,-1){8.22}}
\put(67.89,13.78){\line(-1,1){4.67}}
\put(95.99,14.00){\line(1,1){5.11}}
\put(73.67,9.56){\line(0,-1){4.67}}
\put(89.67,9.56){\line(0,-1){4.67}}
\put(76.56,21.33){\circle*{0.00}}
\put(81.44,21.33){\circle*{0.00}}
\put(86.11,21.33){\circle*{0.00}}
\put(65.22,12.22){\makebox(0,0)[cc]{$\pi_1$}}
\put(65.22,6.89){\makebox(0,0)[cc]{$\pi$}}
\put(71.00,6.89){\makebox(0,0)[cc]{$\gamma_1$}}
\put(92.78,7.11){\makebox(0,0)[cc]{$\gamma_m$}}
\put(98.56,7.11){\makebox(0,0)[cc]{$\pi'$}}
\put(98.33,12.67){\makebox(0,0)[cc]{$\pi_n$}}
\put(94.11,20.22){\makebox(0,0)[cc]{$S$}}
\put(85.89,6.89){\makebox(0,0)[cc]{$T$}}
\put(76.56,6.89){\circle*{0.00}} \put(79.44,6.89){\circle*{0.00}}
\put(82.33,6.89){\circle*{0.00}}
\put(83.44,1.56){\makebox(0,0)[cc]{b}}
\end{picture}
\caption{}\label{f-2}

\end{figure}
\end{center}

Let $S$ and $T$ be two disjoint sets of letters, let ($\pi$,
$\pi_1$, \ldots, $\pi_n$, $\pi'$) be an $S$-band and let ($\pi$,
$\gamma_1$, \ldots, $\gamma_m$, $\pi'$) be a $T$-band. Suppose
that:
\begin{itemize}
\item the medians of these bands form a simple closed
curve,
\item on the boundary  of $\pi$ and on the boundary of $\pi'$ the pairs of $S$-edges separate the pairs of $T$-edges,
\item the start and end edges of these bands are not contained in the
region bounded by the medians of the bands.
\end{itemize}
Then we say that these bands form an \label{stanl}{\em
$(S,T)$-annulus} and the closed curve formed by the parts of
medians of these bands is the {\em median} of this annulus (see
Figure 3b). For every annulus we define the {\em inside}
diagram of the annulus as the subdiagram bounded by the
median of the annulus. The union of the inside diagram and the annulus is called the {\em subdiagram bounded by the annulus}.

We shall call an $S$-band {\em maximal} if it is not contained in any other
$S$-band.

\subsection{Asphericity and free constructions}\label{ss:aafc}

Let $\pp$ be a group presentation.

\begin{assumption}\label{a:1} We shall always assume that all words in $R$ are cyclically reduced, no relator is a proper power in the free group and no relator is a cyclic shift of another relator or its inverse.
\end{assumption}
A group presentation $\pp$ is called (topologically) {\em aspherical} if the universal cover of the presentation complex\footnote{The presentation complex of $\pp$ has one vertex,  edges corresponding to the generators of $\pp$ and 2-cells corresponding to the relators of $\pp$ \cite{LS}.} of $\pp$ is contractible. The presentation $\pp$ is {\em combinatorially aspherical} \cite[Section 6]{CH} (sometimes also called {\em Peiffer aspherical} \cite{BV}) if every spherical \vk diagram over $\pp$ (i.e. a map from the sphere $S^2$ to the presentation complex of the group) is (combinatorially) homotopic to an empty diagram \cite{CCH, CH}.

Recall that a combinatorial homotopy can insert and delete subdiagrams consisting of two mirror image cells that share an edge and also can make {\em diamond moves} \cite{CH}  as shown on Figure \ref{f-7} below.

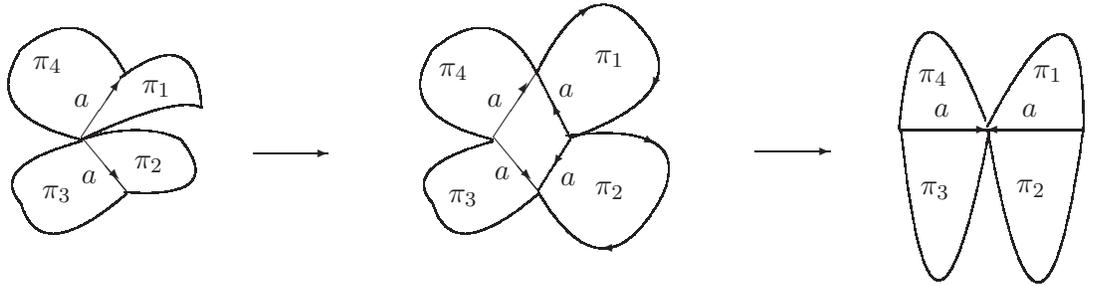
\begin{figure}[ht]
\unitlength 1mm 
\linethickness{0.4pt}
\ifx\plotpoint\undefined\newsavebox{\plotpoint}\fi 
\begin{picture}(148.875,42.375)(23,5)

\put(38.33,25.08){\line(2,3){5.67}}
\put(93.08,25.08){\line(2,3){5.67}}
\put(38.33,25.08){\line(5,-6){6.1}}
\put(93.08,25.08){\line(5,-6){6.1}}
\put(41,29.08){\vector(2,3){2.33}}
\put(95.25,28.83){\vector(2,3){2.33}}
\put(41,21.75){\vector(1,-1){2.33}}
\put(95.75,21.75){\vector(1,-1){2.33}}

\put(38.33,29.92){\makebox(0,0)[cc]{$a$}}
\put(39.33,19.83){\makebox(0,0)[cc]{$a$}}
\put(61.08,22.92){\vector(1,0){10}}
\put(127.83,23.17){\vector(1,0){10}}

\put(147,26){\vector(1,0){11.67}}
\put(171.33,26){\vector(-1,0){12.67}}

\put(152.67,28.67){\makebox(0,0)[cc]{$a$}}
\put(164.33,28.67){\makebox(0,0)[cc]{$a$}}
\qbezier(43.5,33)(53.625,40.5)(54.25,29)
\qbezier(54.25,29)(50.75,30.625)(38.25,24.75)
\qbezier(38.25,24.75)(46.125,27.25)(51.5,24.75)
\qbezier(51.5,24.75)(57.75,16.75)(44,17.75)
\put(48.25,31.5){\makebox(0,0)[cc]{$\pi_1$}}
\put(47.25,21.5){\makebox(0,0)[cc]{$\pi_2$}}
\qbezier(147,26.25)(149.25,51.25)(158.5,27.25)
\qbezier(158.75,26.5)(171.875,51.375)(171.5,25.75)
\put(151.5,33){\makebox(0,0)[cc]{$\pi_4$}}
\put(166.75,33.5){\makebox(0,0)[cc]{$\pi_1$}}
\qbezier(30,36.5)(38.5,43.875)(44,33.75)
\qbezier(84.75,36.5)(93.25,43.875)(98.75,33.75)
\qbezier(30,36.25)(24.875,28)(38.25,24.75)
\qbezier(84.75,36.25)(79.625,28)(93,24.75)
\qbezier(38.25,24.5)(25.75,21.125)(30.25,16.25)
\qbezier(93,24.5)(80.5,21.125)(85,16.25)
\qbezier(30.25,16.25)(32.75,7.625)(44.25,17.5)
\qbezier(85,16.25)(87.5,7.625)(99,17.5)
\put(35,17.5){\makebox(0,0)[cc]{$\pi_3$}}
\put(33.75,34.75){\makebox(0,0)[cc]{$\pi_4$}}
\qbezier(147.25,26)(151.375,-14)(159,26)
\qbezier(158.75,26)(165.875,-14.625)(171.5,26.25)
\put(151.75,18){\makebox(0,0)[cc]{$\pi_3$}}
\put(164.5,18.25){\makebox(0,0)[cc]{$\pi_2$}}
\put(101,29.375){\vector(-1,2){.07}}\multiput(103.25,25)(-.03358209,.065298507){134}{\line(0,1){.065298507}}
\put(101.125,21.375){\vector(-2,-3){.07}}\multiput(103.25,24.75)(-.033730159,-.053571429){126}{\line(0,-1){.053571429}}
\put(105.75,42.25){\vector(3,2){.07}}\qbezier(98.75,33.5)(105.875,47.875)(112.5,39.75)
\put(114.25,31.75){\vector(-1,-2){.07}}\qbezier(112.5,39.75)(120.5,31.125)(103.5,25)
\put(114.5,24.25){\vector(2,-1){.07}}\qbezier(102.75,25.25)(119.875,27.25)(115.5,17.25)
\put(107.5,10.25){\vector(-1,0){.07}}\qbezier(115.5,17.25)(107.75,3)(99,17.75)
\put(93.25,30){\makebox(0,0)[cc]{$a$}}
\put(94.25,20.25){\makebox(0,0)[cc]{$a$}}
\put(102.75,31.25){\makebox(0,0)[cc]{$a$}}
\put(103,19.5){\makebox(0,0)[cc]{$a$}}
\put(108.5,35.5){\makebox(0,0)[cc]{$\pi_1$}}
\put(108.5,18){\makebox(0,0)[cc]{$\pi_2$}}
\put(89,17){\makebox(0,0)[cc]{$\pi_3$}}
\put(87.75,33.75){\makebox(0,0)[cc]{$\pi_4$}}
\end{picture}

\caption{Diamond move: cut along edges labeled by $a$, then fold; $\pi_1,\pi_2,\pi_3,\pi_4$ are cells.}\label{f-7}
\end{figure}

The following statement was proved in \cite{CCH} using different terminology (pictures instead of \vk diagrams). We are not going to show here that our terminology here is equivalent to the one in \cite{CCH} (but see a translation of that terminology into the language of diagrams in \cite[Section 32]{Olbook}). In fact we only need the ``only if" implication of that lemma which is obvious.

\begin{lemma}[Proposition 1.3, \cite{CCH}]\label{l:a} Under Assumption \ref{a:1} a group presentation is combinatorially aspherical if and only if it is aspherical.
\end{lemma}

We shall use the following two ``combination" statements for asphericity from \cite{CCH} (see Theorems 3.7 and 4.3 there). Since the terminology used in \cite{CCH} differs from ours, we provide a proof of the second statement here (the idea of the proof is common for many proofs involving \vk diagrams). The proof of the first statement is similar.

\begin{df} A (connected) planar subcomplex without cut-points $\Delta'$ of a \vk diagram $\Delta$ will be called a {\em holey} $M$-subdiagram. If $\Delta$ does not have holes (i.e. is homeomorphic to a disc), it will be called a {\em disc} subdiagram or simply a {\em subdiagram}.
\end{df}

\begin{lemma}[HNN-extensions]\label{cch1}
Let $\pp_1=\la M \mid L\ra$ be a presentation of a group $T$ and let
$\pp_2=\la M, t\mid  L, t x_i t\iv y_i\iv, i\in I\ra$ be the {\em standard} presentation of an HNN-extension of $T$. Suppose
that no $tx_i t\iv y_i\iv$ is conjugate to any other relator of $\pp_2$. Then $\pp_2$ is
combinatorially aspherical if
$\pp_1$ is combinatorially
aspherical and the images in $T$ of $\{x_i, i\in I\}$ and $\{y_i\mid i\in I\}$ are sets of free
generators.
\end{lemma}

\begin{lemma}[Amalgamated products]\label{cch2}
Let $\pp=\la M_1, M_2\mid L_1,L_2, x_iy_i\iv, i\in I\ra$ be the {\em standard} presentation of an
amalgamated product where $x_i\in M_1, y_i\in M_2$. Then $\pp$ is combinatorially aspherical
if $\pp_1=\la M_1\mid L_1\ra$ and $\pp_2=\la M_2\mid L_2\ra$ are combinatorially aspherical, and both $\{x_i, i\in I\}$ and $\{y_i,i\in I\}$ freely generate free subgroups modulo $\pp_1$ and $\pp_2$ respectively.
\end{lemma}

\proof Let $\Delta$ be a spherical diagram over $\pp$. We can assume that $\Delta$ does not have cut-points, so that every edge of $\Delta$ belongs to the boundary of a cell (otherwise represent $\Delta$ as a connected sum of smaller spherical diagrams). We need to show that $\Delta$ is combinatorially homotopic to a trivial diagram. If $\Delta$ only has cells corresponding to the relations of $\pp_1$ (of $\pp_2$) then we can use asphericity of $\pp_1$ (of $\pp_2$). So assume $\Delta$ has cells corresponding to the relations of both $\pp_1$ and $\pp_2$. Consider a maximal holey subdiagram $\Delta_1$ of $\Delta$ filled with $\pp_1$-cells. Since $\Delta_1\ne \Delta$, it must have a non-trivial boundary component $\p$ (possibly more than one) that bounds a subdiagram $\Delta_2$ (the one that does not contain $\Delta_1$) containing cells. Let $u$ be the label of that boundary component. Note that each letter in $u$ belongs to both a relation of $\pp_1$ and a relation not in $\pp_1$, which then must be one of the relations $x_iy_i\iv$, $i\in I$ since $M_1$, $M_2$ are disjoint. Hence $u$ is a word in $\{x_i, i\in I\}$. Since $\{x_i,i\in I\}$ freely generates a free subgroup modulo $\pp_1$, it freely generates a free subgroup modulo $\pp$ (the group given by $\pp_1$ embeds into the group given by $\pp$). Hence $u$ is freely trivial. Using the diamond moves one can make the reduced boundary of $\Delta_2$ trivial. Since $\Delta_2$ has fewer cells than $\Delta$, we can use induction and conclude that $\Delta_2$ is combinatorially homotopic to a trivial diagram. Hence $\Delta$ is combinatorially homotopic to a spherical diagram with fewer cells and we can again use induction.\endproof

\section{An auxiliary group}
Let $A$ be an alphabet, $B,X,Y$ be two-letter disjoint alphabets which are disjoint from $A$.
Let
$\beta,\beta' \colon A\times X\to B^*$, $\chi, \chi'\colon A\times
X\to Y^*$, $\gamma\colon B\times Y\to Y^*, \psi\colon B\times X\to
Y^*$, $\phi\colon A\times Y\to B^*$ be maps such that the set
$\chi(A\times X)\cup \chi'(A\times X)\cup\gamma(B\times Y)\cup
\psi(B\times X)$ satisfies the small cancelation condition
$C'(\frac1{12})$, the set $\beta(A\times X)\cup \beta'(A\times
X)\cup \phi(A\times Y)$ also satisfies $C'(\frac1{12})$.\footnote{For every alphabet $T$, the notation $T^*$ stands for the set
of all words in $T$. A set of
words $Q$ satisfies the condition $C'(\lambda)$ for some $\lambda>0$
if for every two different cyclic shifts $u,v$ of two words in $Q^{\pm 1}$,
the length of their maximal common prefix is smaller than
$\lambda\min\{|u|, |v|\}$. }

Let $H(A)$ be the group given by the presentation:
$$\begin{array}{l} \la A\cup B\cup X\cup Y\mid
x a = \beta(a,x) a \beta'(a,x) \chi(a,x) x \chi'(a,x),
y a = a\phi(a,y) y,\\
x b = b \psi(b,x)x,
y b = b \gamma(b,y), \forall a\in A, x\in X, b\in B, y\in Y \ra.\end{array}
$$

Figure \ref{f0} shows the cells corresponding to the defining
relations of $H(A)$. The \vk diagrams over the presentation of
$H(A)$ are tesselated by these cells.

\begin{center}
\begin{figure}[ht!]
\unitlength .5mm 
\linethickness{0.4pt}
\ifx\plotpoint\undefined\newsavebox{\plotpoint}\fi 
\begin{picture}(280,71.5)(0,0)
\put(96.75,24.625){\vector(0,1){.141}}\put(96.75,5.25){\line(0,1){38.75}}
\put(163,24.625){\vector(0,1){.141}}\put(163,5.25){\line(0,1){38.75}}
\put(227.25,25.375){\vector(0,1){.141}}\put(227.25,6){\line(0,1){38.75}}
\put(119.375,44){\vector(1,0){.141}}\put(96.75,44){\line(1,0){45.25}}
\put(185.625,44){\vector(1,0){.141}}\put(163,44){\line(1,0){45.25}}
\put(249.875,44.75){\vector(1,0){.141}}\put(227.25,44.75){\line(1,0){45.25}}
\put(249.625,6){\vector(1,0){.141}}\put(227,6){\line(1,0){45.25}}
\put(118.75,48.5){\makebox(0,0)[cc]{$a$}}
\put(185,48.5){\makebox(0,0)[cc]{$b$}}
\put(249.25,49.25){\makebox(0,0)[cc]{$b$}}
\put(92,24.25){\makebox(0,0)[cc]{$y$}}
\put(158.25,24.25){\makebox(0,0)[cc]{$x$}}
\put(222.5,25){\makebox(0,0)[cc]{$y$}}
\put(249,2.75){\makebox(0,0)[cc]{$b$}}
\put(144.75,23.75){\makebox(0,0)[cc]{$y$}}
\put(211,23.75){\makebox(0,0)[cc]{$x$}}
\put(18.313,48.063){\vector(1,3){.141}}\multiput(14.375,37.875)(.067307692,.174145299){117}{\line(0,1){.174145299}}
\put(34.063,63.687){\vector(2,1){.141}}\multiput(22.25,58.25)(.145833333,.06712963){162}{\line(1,0){.145833333}}
\put(34.563,10.188){\vector(2,-1){.141}}\multiput(22.5,15.625)(.148919753,-.06712963){162}{\line(1,0){.148919753}}
\put(73.938,49.063){\vector(-1,2){.141}}\multiput(79,37.75)(-.06705298,.149834437){151}{\line(0,1){.149834437}}
\put(108.063,5.188){\vector(1,0){.141}}\multiput(96.75,5.25)(11.3125,-.0625){2}{\line(1,0){11.3125}}
\put(174.313,5.188){\vector(1,0){.141}}\multiput(163,5.25)(11.3125,-.0625){2}{\line(1,0){11.3125}}
\put(131.625,8.875){\vector(2,1){.141}}\qbezier(119.125,5.125)(133.063,6.563)(141.25,17.25)
\put(197.875,8.875){\vector(2,1){.141}}\qbezier(185.375,5.125)(199.313,6.563)(207.5,17.25)
\put(141.5,30.625){\vector(0,1){.141}}\multiput(141.25,17.125)(.0625,3.375){8}{\line(0,1){3.375}}
\put(207.75,30.625){\vector(0,1){.141}}\multiput(207.5,17.125)(.0625,3.375){8}{\line(0,1){3.375}}
\put(107.25,2.75){\makebox(0,0)[cc]{$a$}}
\put(173.5,2.75){\makebox(0,0)[cc]{$b$}}
\put(134.375,10){\makebox(0,0)[lc]{$\phi(a,y)$}}
\put(200.625,10){\makebox(0,0)[lc]{$\psi(b,x)$}}
\put(16.125,48.875){\makebox(0,0)[cc]{$x$}}
\put(32.125,66.125){\makebox(0,0)[cc]{$a$}}
\put(4.33,23){\makebox(0,0)[cc]{$\beta(a,x)$}}
\put(31.625,8.75){\makebox(0,0)[cc]{$a$}}
\put(75.875,2.625){\makebox(0,0)[cc]{$\beta'(a,x)$}}
\put(74.05,23.875){\makebox(0,0)[rc]{$\chi(a,x)$}}
\put(75.875,50.25){\makebox(0,0)[cc]{$x$}}
\put(61.75,74.05){\makebox(0,0)[cc]{$\chi'(a,x)$}}
\qbezier(272.5,45)(278.625,37.375)(271.25,25.25)
\put(267.75,14){\vector(0,1){.141}}\qbezier(271.25,25.25)(263.75,12.5)(272.25,5.75)
\put(280,20.25){\makebox(0,0)[cc]{$\gamma(b,y)$}}
\put(16,22.25){\vector(1,-3){.141}}\qbezier(14.75,38)(13.125,17.625)(23,15.75)
\put(75.5,23.5){\vector(1,2){.141}}\qbezier(69.25,14.25)(76.875,20.875)(79,38)
\put(58,68){\vector(-2,1){.141}}\qbezier(69,60)(58.625,71.5)(45.75,69)
\put(63,5.75){\vector(3,1){.141}}\qbezier(46.5,4.75)(68.125,1.875)(69.25,14.5)
\end{picture}
\caption{$(A,X)$-cell, $(A,Y)$-cell, $(B,X)$-cell and $(B,Y)$-cell}
\label{f0}
\end{figure}
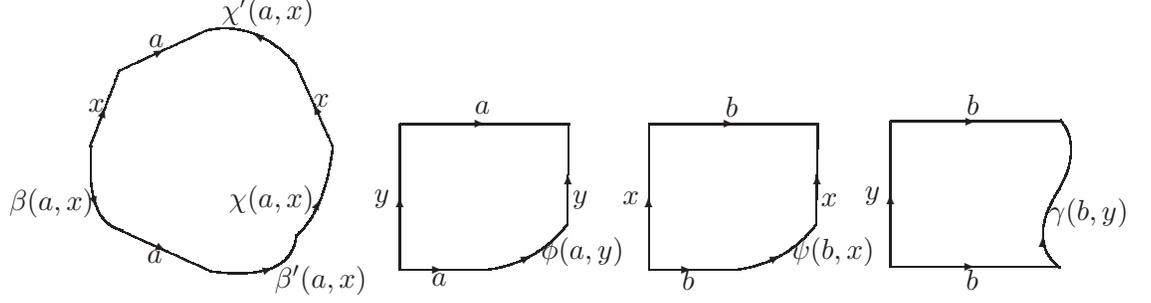
\end{center}

The words $\beta(a,x), \beta'(a,x), \chi(a,x), \chi'(a,x),
\phi(a,y), \psi(a,y), \gamma(b,y)$ are called the {\em large
sections} of the defining relators and the corresponding sections of
boundaries of cells in the \vk diagrams will be called {\em large
sections} of the cells. The following lemma immediately follows from
the fact that the set of all large sections of defining relators
satisfies the property $C'(\frac1{12})$.

\begin{lemma}\label{lsec} The presentation of $H(A)$ satisfies $C'(\frac1{12})$. If a large section of a cell $\pi$ in a \vk diagram $\Delta$ over the presentation of $H(A)$ shares a subpath of at least $\frac1{12}$ of its length with the boundary of another cell $\pi'$ of $\Delta$, then $\pi, \pi'$ cancel.
\end{lemma}

For every \vk diagram $\Delta$ over the presentation of $H(A)$,
and every $p\in A\cup B\cup X\cup Y$ we can consider $p$-bands in $\Delta$. We shall also call it a $A$-band ($B$-band, $X$-band) if $p\in A$ (resp. $B,X$). The
$A$-bands can include $(A,X)$- and $(A,Y)$-cells,
$X$-bands can include $(A,X)$- and $(B,X)$-cells. Maximal
$A$- or $X$-bands can start and end on the boundary of $\Delta$.  A
$B$-band may contain $(B,X)$-, $(B,Y)$-cells. It can
start and end either on $\partial\Delta$ or on the boundary of an
$(A,X)$-cell or on the boundary of an $(A,Y)$-cell.

\begin{lemma}\label{noannul} In a reduced\footnote{We call a diagram {\em
reduced} if it does not contain two mirror image cells sharing an
edge.} \vk diagram over $H(A)$, there are no
\begin{enumerate}
\item $A$-annuli,
\item $X$-annuli,
\item $B$-annuli,
\item $(A,X)$-annuli,
\item $(B,X)$-annuli
\end{enumerate}

\end{lemma}

\proof We prove all five statement by simultaneous induction on the number of cells in the diagram. Suppose that a reduced \vk diagram over the presentation of
$H(A)$ has an annulus of one of the types mentioned in the lemma.
Let $\Delta''$ be the inside subdiagram of $\aaa$, and $\Delta'$ be the diagram bounded by $\aaa$.  We can assume that  there is no annulus of one of the types (1)-(5) of the lemma with a
smaller inside subdiagram.

1. Suppose that $\aaa$ is an $A$-annulus, $a\in A$. The diagram $\Delta''$
does not contain $(A,X)$- or   $(A,Y)$-cells because its boundary
does not contain $A$-edges, and because, by the minimality assumption,
$\Delta''$ does not contain $A$-annuli. Hence all cells in
$\Delta''$ are $(B,X)$- and $(B,Y)$-cells.

Suppose that $\Delta''$ contains cells.

Since $\Delta''$ does not
contain $B$-annuli, the boundary of $\Delta''$ contains $B$-edges.
If the diagram $\Delta'$ contains an $(A,X)$- or a $(B,X)$-cell, then
the $X$-band containing this cell must intersect $\aaa$ twice
creating an $(A,X)$-annulus with smaller inside subdiagram than
$\Delta''$. Hence $\Delta''$ does not contain $(A,X)$- or
$(B,X)$-cells. Hence all cells in $\Delta''$ are $(B,Y)$-cells, and
all cells in $\aaa$ are $(A,Y)$-cells. Since $\Delta''$ does not
contain $B$-annuli, $\partial\Delta''$ must contain $B$-edges, and
so the label of $\partial\Delta''$ is a product of words of the from the set
$\phi(A\times Y)Y$ and their inverses. By the
Greendlinger lemma \cite{LS} one of the large sections $\gamma(B\times Y)$ shares a subword of at
least $1/2$ of its length with the label of $\partial\Delta''$ which
contradicts the assumption that the set of all large sections satisfies
$C'(\frac1{12})$.

If $\Delta''$ does not contain cells, then its boundary label is freely trivial and by Lemma \ref{lsec} two neighbor cells in $\aaa$ cancel.

\begin{figure}[htp!]
\unitlength .5mm 
\linethickness{0.4pt}
\ifx\plotpoint\undefined\newsavebox{\plotpoint}\fi 
\begin{picture}(238.17,130.375)(0,0)
\put(25.313,49.313){\vector(1,3){.141}}\multiput(21.375,39.125)(.067307692,.174145299){117}{\line(0,1){.174145299}}
\put(226.188,48.563){\vector(-1,3){.141}}\multiput(230.125,38.375)(-.067307692,.174145299){117}{\line(0,1){.174145299}}
\put(41.063,64.938){\vector(2,1){.141}}\multiput(29.25,59.5)(.145833333,.06712963){162}{\line(1,0){.145833333}}
\put(210.438,64.188){\vector(-2,1){.141}}\multiput(222.25,58.75)(-.145833333,.06712963){162}{\line(-1,0){.145833333}}
\put(41.563,11.438){\vector(2,-1){.141}}\multiput(29.5,16.875)(.148919753,-.06712963){162}{\line(1,0){.148919753}}
\put(209.938,10.688){\vector(-2,-1){.141}}\multiput(222,16.125)(-.148919753,-.06712963){162}{\line(-1,0){.148919753}}
\put(80.938,50.313){\vector(-1,2){.141}}\multiput(86,39)(-.06705298,.149834437){151}{\line(0,1){.149834437}}
\put(170.563,49.563){\vector(1,2){.141}}\multiput(165.5,38.25)(.06705298,.149834437){151}{\line(0,1){.149834437}}
\put(23.125,50.125){\makebox(0,0)[cc]{$x$}}
\put(228.375,49.375){\makebox(0,0)[]{$x$}}
\put(39.125,67.375){\makebox(0,0)[cc]{$a$}}
\put(212.375,66.625){\makebox(0,0)[]{$a$}}
\put(11.33,24.25){\makebox(0,0)[cc]{$\beta(a,x)$}}
\put(241.17,23.5){\makebox(0,0)[]{$\beta(a,x)$}}
\put(38.625,10){\makebox(0,0)[cc]{$a$}}
\put(212.875,9.25){\makebox(0,0)[]{$a$}}
\put(84.875,5.875){\makebox(0,0)[cc]{$\beta'(a,x)$}}
\put(168.625,5.125){\makebox(0,0)[]{$\beta'(a,x)$}}
\put(81.05,25.125){\makebox(0,0)[rc]{$\chi(a,x)$}}
\put(170.45,24.375){\makebox(0,0)[l]{$\chi(a,x)$}}
\put(82.875,51.5){\makebox(0,0)[cc]{$x$}}
\put(168.625,50.75){\makebox(0,0)[]{$x$}}
\put(79.75,71.3){\makebox(0,0)[cc]{$\chi'(a,x)$}}
\put(172.75,70.55){\makebox(0,0)[]{$\chi'(a,x)$}}
\put(23,23.5){\vector(1,-3){.141}}\qbezier(21.75,39.25)(20.125,18.875)(30,17)
\put(228.5,22.75){\vector(-1,-3){.141}}\qbezier(229.75,38.5)(231.375,18.125)(221.5,16.25)
\put(82.5,24.75){\vector(1,2){.141}}\qbezier(76.25,15.5)(83.875,22.125)(86,39.25)
\put(169,24){\vector(-1,2){.141}}\qbezier(175.25,14.75)(167.625,21.375)(165.5,38.5)
\put(65,69.25){\vector(-2,1){.141}}\qbezier(76,61.25)(65.625,72.75)(52.75,70.25)
\put(186.5,68.5){\vector(2,1){.141}}\qbezier(175.5,60.5)(185.875,72)(198.75,69.5)
\put(70,7){\vector(3,1){.141}}\qbezier(53.5,6)(75.125,3.125)(76.25,15.75)
\put(181.5,6.25){\vector(-3,1){.141}}\qbezier(198,5.25)(176.375,2.375)(175.25,15)
\qbezier(86.5,38.25)(121.375,55.125)(165.75,38.5)
\qbezier(76,61.5)(117.25,79.5)(175.5,60.5)
\qbezier(53,70.5)(115.125,129)(198.75,69.5)
\qbezier(29.5,59.5)(110.5,193.375)(222.5,58.75)
\put(120,109.25){\makebox(0,0)[cc]{$\aaa_1$}}
\put(120.5,55){\makebox(0,0)[cc]{$\aaa_2$}}
\put(120.75,83.5){\makebox(0,0)[cc]{$\Delta''$}}
\end{picture}

\caption{}\label{p6}
\end{figure}

\begin{figure}[htp!]

\unitlength .5mm 
\linethickness{0.4pt}
\ifx\plotpoint\undefined\newsavebox{\plotpoint}\fi 
\begin{picture}(241.17,114.5)(0,0)
\put(28.313,84.313){\vector(1,3){.141}}\multiput(24.375,74.125)(.067307692,.174145299){117}{\line(0,1){.174145299}}
\put(229.188,83.563){\vector(-1,3){.141}}\multiput(233.125,73.375)(-.067307692,.174145299){117}{\line(0,1){.174145299}}
\put(44.063,99.938){\vector(2,1){.141}}\multiput(32.25,94.5)(.145833333,.06712963){162}{\line(1,0){.145833333}}
\put(213.438,99.188){\vector(-2,1){.141}}\multiput(225.25,93.75)(-.145833333,.06712963){162}{\line(-1,0){.145833333}}
\put(44.563,46.438){\vector(2,-1){.141}}\multiput(32.5,51.875)(.148919753,-.06712963){162}{\line(1,0){.148919753}}
\put(212.938,45.688){\vector(-2,-1){.141}}\multiput(225,51.125)(-.148919753,-.06712963){162}{\line(-1,0){.148919753}}
\put(83.938,85.313){\vector(-1,2){.141}}\multiput(89,74)(-.06705298,.149834437){151}{\line(0,1){.149834437}}
\put(173.563,84.563){\vector(1,2){.141}}\multiput(168.5,73.25)(.06705298,.149834437){151}{\line(0,1){.149834437}}
\put(26.125,85.125){\makebox(0,0)[cc]{$x$}}
\put(231.375,84.375){\makebox(0,0)[]{$x$}}
\put(42.125,102.375){\makebox(0,0)[cc]{$a$}}
\put(215.375,101.625){\makebox(0,0)[]{$a$}}
\put(13.33,59.25){\makebox(0,0)[cc]{$\beta(a,x)$}}
\put(246.17,58.5){\makebox(0,0)[]{$\beta(a,x)$}}
\put(41.625,45){\makebox(0,0)[cc]{$a$}}
\put(215.875,44.25){\makebox(0,0)[]{$a$}}
\put(87.875,40.875){\makebox(0,0)[cc]{$\beta'(a,x)$}}
\put(170.625,40.125){\makebox(0,0)[]{$\beta'(a,x)$}}
\put(84.05,60.125){\makebox(0,0)[rc]{$\chi(a,x)$}}
\put(173.45,59.375){\makebox(0,0)[l]{$\chi(a,x)$}}
\put(85.875,86.5){\makebox(0,0)[cc]{$x$}}
\put(171.625,85.75){\makebox(0,0)[]{$x$}}
\put(82.75,106.3){\makebox(0,0)[cc]{$\chi'(a,x)$}}
\put(176.75,105.55){\makebox(0,0)[]{$\chi'(a,x)$}}
\put(26,58.5){\vector(1,-3){.141}}\qbezier(24.75,74.25)(23.125,53.875)(33,52)
\put(231.5,57.75){\vector(-1,-3){.141}}\qbezier(232.75,73.5)(234.375,53.125)(224.5,51.25)
\put(85.5,59.75){\vector(1,2){.141}}\qbezier(79.25,50.5)(86.875,57.125)(89,74.25)
\put(172,59){\vector(-1,2){.141}}\qbezier(178.25,49.75)(170.625,56.375)(168.5,73.5)
\put(68,104.25){\vector(-2,1){.141}}\qbezier(79,96.25)(68.625,107.75)(55.75,105.25)
\put(189.5,103.5){\vector(2,1){.141}}\qbezier(178.5,95.5)(188.875,107)(201.75,104.5)
\put(73,42){\vector(3,1){.141}}\qbezier(56.5,41)(78.125,38.125)(79.25,50.75)
\put(184.5,41.25){\vector(-3,1){.141}}\qbezier(201,40.25)(179.375,37.375)(178.25,50)
\qbezier(89.5,73.25)(124.375,90.125)(168.75,73.5)
\qbezier(79,96.5)(120.25,114.5)(178.5,95.5)
\put(123.5,90){\makebox(0,0)[cc]{$\aaa_2$}}
\qbezier(57,41)(141.625,-3.625)(200.75,40.25)
\qbezier(29.75,53.25)(130.5,-49.25)(225.25,51.25)
\put(119.75,9.75){\makebox(0,0)[cc]{$\aaa_1$}}
\put(120.5,45.25){\makebox(0,0)[cc]{$\Delta''$}}
\end{picture}

\caption{}\label{p7}
\end{figure}

\begin{figure}[htp!]

\unitlength .5mm 
\linethickness{0.4pt}
\ifx\plotpoint\undefined\newsavebox{\plotpoint}\fi 
\begin{picture}(236.5,100.625)(0,44)
\put(84.188,84.313){\vector(-1,3){.141}}\multiput(88.125,74.125)(-.067307692,.174145299){117}{\line(0,1){.174145299}}
\put(175.813,83.563){\vector(1,3){.141}}\multiput(171.875,73.375)(.067307692,.174145299){117}{\line(0,1){.174145299}}
\put(68.438,99.938){\vector(-2,1){.141}}\multiput(80.25,94.5)(-.145833333,.06712963){162}{\line(-1,0){.145833333}}
\put(191.563,99.188){\vector(2,1){.141}}\multiput(179.75,93.75)(.145833333,.06712963){162}{\line(1,0){.145833333}}
\put(67.938,46.438){\vector(-2,-1){.141}}\multiput(80,51.875)(-.148919753,-.06712963){162}{\line(-1,0){.148919753}}
\put(192.063,45.688){\vector(2,-1){.141}}\multiput(180,51.125)(.148919753,-.06712963){162}{\line(1,0){.148919753}}
\put(28.563,85.313){\vector(1,2){.141}}\multiput(23.5,74)(.06705298,.149834437){151}{\line(0,1){.149834437}}
\put(231.438,84.563){\vector(-1,2){.141}}\multiput(236.5,73.25)(-.06705298,.149834437){151}{\line(0,1){.149834437}}
\put(86.375,85.125){\makebox(0,0)[]{$x$}}
\put(173.625,84.375){\makebox(0,0)[]{$x$}}
\put(70.375,102.375){\makebox(0,0)[]{$a$}}
\put(189.625,101.625){\makebox(0,0)[]{$a$}}
\put(99.17,59.25){\makebox(0,0)[]{$\beta(a,x)$}}
\put(159.83,58.5){\makebox(0,0)[]{$\beta(a,x)$}}
\put(70.875,45){\makebox(0,0)[]{$a$}}
\put(189.125,44.25){\makebox(0,0)[]{$a$}}
\put(25.625,40.875){\makebox(0,0)[]{$\beta'(a,x)$}}
\put(235.375,40.125){\makebox(0,0)[]{$\beta'(a,x)$}}
\put(28.45,60.125){\makebox(0,0)[l]{$\chi(a,x)$}}
\put(231.55,59.375){\makebox(0,0)[r]{$\chi(a,x)$}}
\put(26.625,86.5){\makebox(0,0)[]{$x$}}
\put(233.375,85.75){\makebox(0,0)[]{$x$}}
\put(30.75,106.3){\makebox(0,0)[]{$\chi'(a,x)$}}
\put(230.25,105.55){\makebox(0,0)[]{$\chi'(a,x)$}}
\put(86.5,58.5){\vector(-1,-3){.141}}\qbezier(87.75,74.25)(89.375,53.875)(79.5,52)
\put(173.5,57.75){\vector(1,-3){.141}}\qbezier(172.25,73.5)(170.625,53.125)(180.5,51.25)
\put(27,59.75){\vector(-1,2){.141}}\qbezier(33.25,50.5)(25.625,57.125)(23.5,74.25)
\put(233,59){\vector(1,2){.141}}\qbezier(226.75,49.75)(234.375,56.375)(236.5,73.5)
\put(44.5,104.25){\vector(2,1){.141}}\qbezier(33.5,96.25)(43.875,107.75)(56.75,105.25)
\put(215.5,103.5){\vector(-2,1){.141}}\qbezier(226.5,95.5)(216.125,107)(203.25,104.5)
\put(39.5,42){\vector(-3,1){.141}}\qbezier(56,41)(34.375,38.125)(33.25,50.75)
\put(220.5,41.25){\vector(3,1){.141}}\qbezier(204,40.25)(225.625,37.375)(226.75,50)
\qbezier(88,74.5)(129.875,88)(172.25,73.5)
\qbezier(81,94)(116.5,102.125)(180,93.75)
\qbezier(80.5,94)(124.5,129.625)(179.5,93.75)
\qbezier(57.75,105.25)(121.25,150.625)(203.75,104.5)
\put(128.5,119){\makebox(0,0)[cc]{$\aaa_1$}}
\put(126,87.25){\makebox(0,0)[cc]{$\aaa_2$}}
\put(126.75,104){\makebox(0,0)[cc]{$\Delta''$}}
\end{picture}

\caption{}\label{p8}
\end{figure}

\begin{figure}[htp!]
\unitlength .5mm 
\linethickness{0.4pt}
\ifx\plotpoint\undefined\newsavebox{\plotpoint}\fi 
\begin{picture}(236.5,110.75)(0,23)
\put(84.188,84.313){\vector(-1,3){.141}}\multiput(88.125,74.125)(-.067307692,.174145299){117}{\line(0,1){.174145299}}
\put(175.813,83.563){\vector(1,3){.141}}\multiput(171.875,73.375)(.067307692,.174145299){117}{\line(0,1){.174145299}}
\put(68.438,99.938){\vector(-2,1){.141}}\multiput(80.25,94.5)(-.145833333,.06712963){162}{\line(-1,0){.145833333}}
\put(191.563,99.188){\vector(2,1){.141}}\multiput(179.75,93.75)(.145833333,.06712963){162}{\line(1,0){.145833333}}
\put(67.938,46.438){\vector(-2,-1){.141}}\multiput(80,51.875)(-.148919753,-.06712963){162}{\line(-1,0){.148919753}}
\put(192.063,45.688){\vector(2,-1){.141}}\multiput(180,51.125)(.148919753,-.06712963){162}{\line(1,0){.148919753}}
\put(28.563,85.313){\vector(1,2){.141}}\multiput(23.5,74)(.06705298,.149834437){151}{\line(0,1){.149834437}}
\put(231.438,84.563){\vector(-1,2){.141}}\multiput(236.5,73.25)(-.06705298,.149834437){151}{\line(0,1){.149834437}}
\put(86.375,85.125){\makebox(0,0)[]{$x$}}
\put(173.625,84.375){\makebox(0,0)[]{$x$}}
\put(70.375,102.375){\makebox(0,0)[]{$a$}}
\put(189.625,101.625){\makebox(0,0)[]{$a$}}
\put(99.17,59.25){\makebox(0,0)[]{$\beta(a,x)$}}
\put(159.83,58.5){\makebox(0,0)[]{$\beta(a,x)$}}
\put(70.875,45){\makebox(0,0)[]{$a$}}
\put(189.125,44.25){\makebox(0,0)[]{$a$}}
\put(25.625,40.875){\makebox(0,0)[]{$\beta'(a,x)$}}
\put(235.375,40.125){\makebox(0,0)[]{$\beta'(a,x)$}}
\put(28.45,60.125){\makebox(0,0)[l]{$\chi(a,x)$}}
\put(231.55,59.375){\makebox(0,0)[r]{$\chi(a,x)$}}
\put(26.625,86.5){\makebox(0,0)[]{$x$}}
\put(233.375,85.75){\makebox(0,0)[]{$x$}}
\put(30.75,106.3){\makebox(0,0)[]{$\chi'(a,x)$}}
\put(230.25,105.55){\makebox(0,0)[]{$\chi'(a,x)$}}
\put(86.5,58.5){\vector(-1,-3){.141}}\qbezier(87.75,74.25)(89.375,53.875)(79.5,52)
\put(173.5,57.75){\vector(1,-3){.141}}\qbezier(172.25,73.5)(170.625,53.125)(180.5,51.25)
\put(27,59.75){\vector(-1,2){.141}}\qbezier(33.25,50.5)(25.625,57.125)(23.5,74.25)
\put(233,59){\vector(1,2){.141}}\qbezier(226.75,49.75)(234.375,56.375)(236.5,73.5)
\put(44.5,104.25){\vector(2,1){.141}}\qbezier(33.5,96.25)(43.875,107.75)(56.75,105.25)
\put(215.5,103.5){\vector(-2,1){.141}}\qbezier(226.5,95.5)(216.125,107)(203.25,104.5)
\put(39.5,42){\vector(-3,1){.141}}\qbezier(56,41)(34.375,38.125)(33.25,50.75)
\put(220.5,41.25){\vector(3,1){.141}}\qbezier(204,40.25)(225.625,37.375)(226.75,50)
\qbezier(88,74.5)(129.875,88)(172.25,73.5)
\qbezier(81,94)(116.5,102.125)(180,93.75)
\put(126,87.25){\makebox(0,0)[cc]{$\aaa_2$}}
\qbezier(81.25,52.5)(130.375,29.375)(177,52.75)
\qbezier(55.75,41)(132.5,-5.125)(204.25,40.25)
\put(125.5,62.25){\makebox(0,0)[cc]{$\Delta''$}}
\put(126.75,27.25){\makebox(0,0)[cc]{$\aaa_1$}}
\end{picture}

\caption{}\label{p9}
\end{figure}

2. The case when $\aaa$ is an $X$-annulus is ruled out similarly.

3. Suppose that $\aaa$ is a $B$-annulus. Since $\Delta$ does not
have $(B,X)$-annuli with smaller inside subdiagrams than $\Delta''$,
$\Delta'$ does not contain $(B,X)$-cells. Hence $\Delta'$ consists
of $(B,Y)$-cells. But then the boundary of $\Delta''$ consists of
$Y$-edges, hence $\Delta''$ does not contain $B$-bands. Thus
$\Delta''$ does not contain cells. Hence the label of
$\partial\Delta''$ is freely trivial. The fact that $\gamma(B\times
Y)$ satisfies the condition $C'(\frac1{12})$ implies that this set
freely generates the free group, hence the label of $\Delta'$ is
freely trivial as well, and so $\aaa$ contains a pair of cells that
cancel by Lemma \ref{lsec}.

4. Suppose that $\aaa$ is an $(A,X)$-annulus, composed of an
$A$-band $\aaa_1$ and a $X$-band $\aaa_2$. By the minimality of
$\Delta''$, we can assume that $\aaa$ contains only two cells having
$X$- or $A$-edges, the two common cells $\pi_1, \pi_2$ of $\aaa_1$
and $\aaa_2$. Hence all cells in $\Delta''$ are $(B,Y)$-cells. Now
the condition $C'(\frac1{12})$ satisfied by the set of large
sections immediately implies that $\Delta''$ does not contain cells.
Thus the label of the boundary of $\Delta''$ is freely trivial.
There are four possibilities for the diagram $\Delta'$ depending on
which edges on $\partial\pi_1$ and $\partial\pi_2$ are connected by
$\aaa_1$ and $\aaa_2$. These possibilities are depicted on Figures
\ref{p6}-\ref{p9}. In each of these Figures let $\p_1$ be the
intersection of $\partial\Delta''$ with $\partial\aaa_1$ and $\p_2$
be the intersection of $\partial\Delta''$ with $\aaa_2$. It is easy
to see that in each of four cases $\Lab(\p_1)$ should be conjugate to
$\Lab(\p_2)$ in the free group. Unless both
$\Lab(\p_1)$ and $\Lab(\p_2)$ are empty, this leads to a contradiction as follows.

In the case of Figure \ref{p6}, $\Lab(\p_1)$ is a product of
words from the set $\phi(\{a\}\times Y)Y$ and their inverses, $\Lab(\p_2)$ is a
word in $B$. In the case of Figure \ref{p7},  $\Lab(\p_1)$ is a
product of words from the set $\phi(\{a\}\times Y)Y$ and their inverses ($a$
here is the label of the $a$-edges in $\pi_1$ and $\pi_2$),
$\Lab(\p_2)$ is a product of word of the form $B\psi(B\times \{x\})$ and their
inverses ($x$ here is the label of $x$-edges in $\pi_1, \pi_2$). In
the case of Figure \ref{p8}, $\Lab(\p_1)$ is a word in $Y$, while
$\Lab(\p_2)$ is a word in $B$. In the case of Figure \ref{p9},
$\Lab(\p_1)$ is a word in $Y$, $\Lab(\p_2)$ is a product of words from $B\psi(B\times\{x\})$ and their inverses. In each of these cases the words $\Lab(\p_1)$ and $\Lab(\p_2)$ are obviously non-conjugate in the free group.

5. Suppose that $\aaa$ is a $(B,X)$-annulus composed of a $B$-band
$\aaa_1$ and an $X$-band $\aaa_2$. Then $\aaa_2$ cannot contain
$(A,X)$-cells (otherwise the $A$-band containing that cell will form an $(A,X)$-annulus with $\aaa_2$ whose
inside subdiagram would be smaller than $\Delta''$). Hence all cells in
$\aaa_1$ are $(B,Y)$-cells and all cells in $\aaa_2$ are
$(B,X)$-cells. If $\aaa_2$ is not empty, then the $B$-band starting
on the boundary of a cell in $\aaa_2$ would form a $(B,X)$-annulus
with $\aaa_2$ whose inside subdiagram would be smaller than
$\Delta''$. Thus $\aaa_2$ is empty and $\pi_1, \pi_2$ cancel, a
contradiction.
\endproof

\begin{lemma}\label{l:int}
The sets $A\cup B$, $X\cup Y$ freely generate free subgroups  of $H(A)$, and $\la
A\cup B\ra\cap \la X\cup Y\ra=\{1\}$.
\end{lemma}

\proof Suppose that there exists a reduced non-empty diagram
$\Delta$ over $H(A)$ with $\partial(\Delta)=\p\q$, $\Lab(\p)$ is a reduced group word in $A\cup B$, $\Lab(\q)$ is a reduced group word in $X\cup Y$, one
but not both of these words may be empty. If $\Delta$ contains an
$(A,X)$-cell $\pi$, then consider the maximal $A$-band $\aaa$
containing that cell. That band must start and end on $\p$ because
$q$ does not contain $A$-edges. Then the maximal $X$-band $\xxx$
containing $\pi$ intersects $\aaa$. By Lemma \ref{noannul}, $\xxx$
must end on $\p$. This is a contradiction since $\p$ does not contain
$X$-edges. Hence $\Delta$ does not contain $(A,X)$-cells.

If $\Delta$ contains an $(A,Y)$-cell $\pi$, then let $\aaa$ be the maximal $A$-band
containing $\pi$. Let $\Delta'$ be the subdiagram of $\Delta$
bounded by the median of $\aaa$ and $\p$. The boundary of $\Delta'$ is of the form $\p'\q'$ where $\p'$ is a subpath of $\p$, $\q'$ is a subpath of a side of $\aaa$. Note that $\q'$ must contain $y$-edges otherwise two neighbor cells in $\aaa$ cancel (by Lemma \ref{lsec}). Since $\p$ does not contain $y$-edges, the diagram $\Delta'$ must contain cells.

We can assume that $\Delta'$ does not contain $(A,Y)$-cells (otherwise we can consider a smaller subdiagram). Hence it consists of $(B,Y)$-cells. The condition $C'(\frac1{12})$ implies that the large section of one of the $(B,Y)$-cells in $\Delta$ shares a subpath of at least $\frac16$ of its length with $\q'$ which is impossible because the set of all large sections satisfies $C'(\frac1{12})$. This contradiction
proves both statements of
the lemma. \endproof

\begin{df} \label{df:34} If $W$ is a word in an alphabet $Z$, $Z'\subset Z$, then the {\em projection} of $W$ onto the alphabet $Z'$ is the word obtained by deleting all letters of $Z\setminus Z'$ (and their inverses) from $W$.
\end{df}

\begin{lemma}\label{l0}
Every positive word $W$ in $A\cup B\cup X\cup Y$ is equal in $H(A)$ to a
unique word of the form $UV$ where $U\in (A\cup B)^*$, $V\in (X\cup Y)^*$. The projections of $U$ and $W$ onto the alphabet $A$ are equal and the projections of $V$ and $W$ on the alphabet $X$ are equal.
\end{lemma}

\proof For the proof of existence and the statement about projections, it is enough to prove the following

{\bf Claim.} For every positive word $W$ in the alphabet $X\cup Y$ and every $p\in A\cup B$ there
exist positive words $U\in (A\cup B)^*$, $V\in (X\cup Y)^*$ such that  $$Wp= UV$$
in $H(A)$; the projections onto $X$ of  $V$ and $W$ are the same;
if $p\in B$, then $U\equiv p$; if $p\in A$, then $U$ contains only
one $A$-letter, $p$.

If $W$ is empty, then there is nothing to prove. Let $W\equiv W'z$. We consider several cases depending on which set $z$ or $p$ belongs to. In each case, the proof is by induction on the pair $(m,n)$
where $m$ is the number of $X$-letters and $n$ is the number of
$Y$-letters in $W$.

{\bf Step 1.} We prove that the Claim is true if $p\in B$. If $m=0$, then the result follows from
the $(B,Y)$-relations. Suppose we have proved the statement for $m$
and suppose that $W$ contains $m+1$ occurrences of letters from $X$.
Then $W=W'xW''$ where $x\in X$, $W''\in Y^*$. Therefore
$Wp=W'xW''p=W'xpV_1$ for some word $V_1\in Y^*$ since $W''$ does not
contain $X$-letters. By the $(p,x)$-relation, we have
$Wp=W'p\psi(p,x)xV_1$. Since $W'$ has $m$ $X$-letters, the result
follows by induction. Thus our Claim is true if $p\in B$.

{\bf Step 2.} We prove that the Claim is true if $p\in A, z\in Y$.
We have: $Wp=W'zp=W'p\phi(p,z)z$. Since $W'$ has fewer letters
than $W$, we deduce that $Wp=U_1V_1\phi(p,z)z$ where $U_1\in (A\cup B)^*,
V_1\in (X\cup Y)^*$, the only $A$-letter in $U_1$ is $p$, and the number of
$X$-letters in $V_1$ is the same as in $W'$. Since $\phi(p,z)\in
B^*$, we can finish the proof by applying the statement of Step 1 several times.

{\bf Step 3.} Finally we prove the Claim if $p\in A, z\in X$. In this case we have

\begin{equation}\label{eq45}
Wp=W'zp=W'\beta(p,z)p\beta'(p,z)\chi(p,z)z\chi'(p,z).
\end{equation}
Since $\beta(p,z)$ is a positive word in the alphabet $B$, we can apply the statement of Step 1 several times, and find $V_1\in (X\cup Y)^*$ such that $W'\beta(p,z)=\beta(p,z)V_1$, and the number of $X$-letters in $V_1$ is the same as in $W'$. Then from (\ref{eq45}), we obtain
\begin{equation}\label{eq46}
Wp=\beta(p,z)V_1 p\beta'(p,z)\chi(p,z)z\chi'(p,z).
\end{equation}
Since $V_1$ contains fewer $X$-letters than $W$, we can apply first the statements of Steps 2 and  1, and then the inductive hypothesis and obtain that
$$V_1p=U_1V_2$$
where $U_1\in (A\cup B)^*$ contains exactly one $A$-letter, $p$, $V_2\in (X\cup Y)^*$ contains the same number of $X$-letters as $V_1$. Now from (\ref{eq46}) we obtain
$$
Wp=\beta(p,z)U_1V_2\beta'(p,z)\chi(p,z)z\chi'(p,z)
$$
and it remains to apply the statement of Step 1 again several times.

This completes the proof of existence and the statement about projections.

To show uniqueness of $U, V$ we need to prove that there are no reduced group words $U,U'$ over $A\cup B$, $V,V'$ over $X\cup Y$ such that $UV=U'V'$ in $H(A)$
and either $U\not= U'$ or $V\not= V'$, but that follows
immediately from Lemma \ref{l:int}.
\endproof

\begin{rk}\label{rk:0}  Since the presentation of the group $H(A)$ satisfies $C'(\frac1{12})$, it is combinatorially aspherical (in fact even diagrammatically aspherical) \cite{LS}. Hence Lemma \ref{l0} implies that for every positive $W$ in $A\cup B\cup X\cup Y$ there exists unique (up to combinatorial homotopy) \vk diagram $\Upsilon(W)$ with boundary label $W\iv UV$ where $U=U(W)$ is a word in $A\cup B$, $V=V(W)$ is a word in $X\cup Y$. Words $U, V$ are uniquely determined by $W$.
\end{rk}

\begin{rk} \label{rk:1}
It is easy to see that the lengths of words $U, V$ are in general at
least exponential in terms of the length of $W$ (in fact the maximal
sum $|U|+|V|$ for $|W|=m$ is at least $\exp O(m^2)$). But the proof
of Lemma \ref{l0} gives an algorithm of computing $U, V$ given $W$.
\end{rk}

\begin{notation}\label{not1} Let $u_a, a\in A$, be positive words in $A$ satisfying $C'(\frac1{12})$, $v_a, a\in A$, be positive words in $X$ satisfying $C'(\frac1{12})$. The words $u_a, a\in A$, and $v_a, a\in A$, will be called $A$- and $X$-blocks respectively. Let $E$ be the subgroup of $H(A)$ generated by $u_av_a$. We shall denote $u_av_a$ by $\mu(a)$.
Let $\N$ be a normal subgroup of $E$ generated by some set $\mu(R)$ where $R$ is a set of words in $A$, and $\ncn$ be the normal closure of $\N$ in $H(A)$.
\end{notation}

\begin{lemma}\label{l:mu} If two words $v_a$, $v_b$ share a subword of length $\frac1{12}\min(|v_a|, |v_b|)$, then $a=b$ (the same statement for $u$ instead of $v$ is also true).
\end{lemma}
\proof Indeed, by $C'(\frac1{12})$ in this case $v_a\equiv v_b$. Hence $a=b$.\endproof

\begin{lemma}\label{l:free} $E$ is a free subgroup of $H(A)$ freely generated by  $\{\mu(a), a\in A\}$.
\end{lemma}

\proof Indeed, let $\Delta$ is reduced \vk diagram over the presentation of $H(A)$ whose reduced boundary label is a product of words $\mu(a)=u_av_a$. If $\Delta$ has cells, then by the Greendlinger lemma \cite{LS}, one of the cells shares a non-trivial subpath of its large section with $\partial\Delta$, but large sections consist of $B$- and $Y$-edges, a contradiction. Hence $\Delta$ has no cells and the boundary label of $\Delta$ is freely trivial. The small cancelation condition $C'(\frac1{12})$ then implies that $\partial\Delta$ is empty. Hence the words $\mu(a), a\in A$, freely generate $E$.
\endproof

\begin{lemma}\label{l:det} Let $W=\mu(r)$ for some positive word $r$ in the alphabet $A$. Let $U=U(W), V=V(W)$ be the words determined by Lemma \ref{l0}. Then each of the three words $U,V,W$ uniquely determines the two other words and (up to combinatorial homotopy) the diagram over $H(A)$ for the equality $W=UV$.
\end{lemma}

\proof Indeed, by Lemma \ref{l0}, the projections of $U$ and $W$ on $A$ coincide. If $r=a_1\ldots a_l$, $a_i\in A$, then both projections are equal to $u_{a_1}\ldots u_{a_l}$. Since the words $u_{a_i}$ satisfy $C'(\frac1{12})$, $U_A$ determines $r$. Hence $U$ determines $W$. Since $W$ determines $V$ by Lemma \ref{l0}, $U$ determines $W$ and $V$.  The fact that $V$ determines $U,W$ is proved similarly (one needs to consider the projections onto $X$). The fact that $W$ determines $U,V$ is in Lemma \ref{l0}.

The fact that $W,U,V$ uniquely (up to combinatorial homotopy) determine the diagram $\Upsilon(W)$ for the equality $W=UV$ is in Remark \ref{rk:0}.
\endproof

\begin{notation} \label{not6} Consider a \vk diagram $\Delta$ over the presentation that
consists of the defining relations of $H(A)$ and all reduced words  from
$\mu(R)$ (Notation \ref{not1}). The cells corresponding to relations from $\mu(R)$ will be called
$\N$-{\em cells.}

A holey subdiagram of $\Delta$ where the label of each boundary component is freely equal to a product of words $\mu(a)$ and their inverses will be called a {\em holey $E$-subdiagram}. In particular, an $\N$-cell is an $E$-subdiagram. Using diamond moves, we can make all boundary components of all maximal holey $E$-subdiagrams reduced.
The small cancelation condition $C'(\frac1{12})$
implies that an $A$- or $X$-block in a product of words of the form $\mu(a)=u_av_a$ and their inverses shares
less than 1/6 of its length with the neighbor blocks. Thus the  reduced label of the boundary component of a holey $E$-subdiagram is a concatenation (without cancelations) of words $u_a'v_a'$ and their inverses where $u_a'$ is a subword of an $A$-block of length at least $5/6$ of the length of the block and $v_a'$ is a subword of an $X$-block  of length at least $5/6$ of the length of the block. The words $u_a', v_a'$ are called $A$- and $X$-subblocks. Thus each boundary label of $\Delta'$ is a reduced word that is a concatenation of
$A$-{subblocks} and $X$-{subblocks}.

So we will assume that this property holds and call a maximal holey $E$-subdiagrams $\N$-{\em subdiagrams} of $\Delta$. Note that by the congruence extension property that we shall prove later, the label of every boundary component of every $\N$-subdiagrams is in $\N$ (Lemma \ref{l:cep}).
Consider the graph $\mho=\mho(\Delta)$ whose
vertices are the holey $\N$-subdiagrams of $\Delta$.  Two vertices are
connected if there exists an $X$-band connecting the corresponding
subdiagrams. Note that $\mho$ is a planar graphs since $X$-bands do not intersect.

We say that two paths $\p_1, \p_2$ in a diagram $\Delta$  {\em share large $X$-portion} if $p_1=\q_1\p\q_1', \p_2=\q_2'\p\q_2$, $\Lab(\p)$ is a subword of a word $v_a$ (resp. $u_a$), $a\in A$, of length at least $\frac1{12}|v_a|$ (resp. $\frac1{12}|u_a|$).
\end{notation}

We shall need the following lemma.

\begin{lemma}\label{lp1.1} Suppose that $\Delta$ is reduced.

(a) Let $\Pi$ be a holey $E$-subdiagram of $\Delta$. Then no $X$-band or $A$-band can start and end on $\partial\Pi$.

(a') Suppose that $\partial\Delta=\p\q$ where $\Lab(\p)$ is a subword of a product of words $\mu(a)$. Then no $X$-band or $A$-band can start and end on $\p$.

(b) Let $\Pi, \Pi'$ be holey $E$-subdiagrams in $\Delta$. Let $e_1, e_2$ be
consecutive $X$-edges (resp. $A$-edges) of $\partial\Pi$ and the $X$-bands
$\xxx,\xxx'$ starting on $e_1$ and $e_2$ end on $\partial\Pi'$. Then
$\xxx,\xxx'$ are empty bands.

(b') Suppose that $\partial\Delta=\p\q$ where $\Lab(\p)$ is a subword of a product of words $\mu(a)$. Let $\Pi$ be a holey $E$-subdiagram of $\Delta$. Let $e_1,e_2$ be consecutive $X$-edges (resp. $A$-edges) of $\partial\Pi$ and the $X$-bands (resp. $A$-bands) starting on $e_1$, $e_2$ end on $\p$. Then $\xxx, \xxx'$ are empty bands.


(c) Suppose that $\partial\Delta=\p\q$ where $\Lab(\p)$ is a subword of a product of words $\mu(a)$, $a\in A$, and $\Lab(\q)$ does not have letters from $X$ (resp. letters from $A$). Suppose $\Delta$ has $\N$-cells. Then there exists a holey $\N$-subdiagram that shares a big $X$-portion (resp. big $A$-portion) of its boundary with $\p$.
\end{lemma}

\proof We shall prove the five statements by a simultaneous
induction on the number of cells in $\Delta$. For diagrams without
cells all five statements are obviously true. Suppose that $\Delta$
is a minimal counterexample (with respect to the number of cells).

Suppose that $\Delta$ {\em does not satisfy (a)}, and an $X$-band
$\xxx$ starts and ends on the same block subdiagram  $\Pi$ in $\Delta$.  Let
$e_1$ be the start edge and $e_2$ be the end edge of $\xxx$, and (without loss of generality) assume that $e_1$ precedes $e_2$ in $\partial\Pi$ (traced counterclockwise).
Let $\p$ be the subpath of $\partial\Pi$ between $e_1$ and $e_2$. Let $\Delta'$ be the subdiagram of $\Delta$ bounded by $\p$ and the median of $\xxx$. Suppose that $\Delta'$ has $\N$-cells. Then by (c) there exists an $\N$-cell $\pi$ that shares a big $X$-portion of its boundary with $\p$. By Lemma \ref{l:mu} the union of $\Pi$ and $\pi$ is an $E$-subdiagram which contradicts maximality of $\Pi$. Hence $\Delta'$ does not contain $\N$-cells, so $\Delta'$ is a diagram over the presentation of
$H(A)$. Therefore, since $X$-bands do not intersect, every $X$-band of $\Delta'$ starting on $\p$ must end on $\p$. Thus we can
assume that there are no $X$-edges on $\p$.

Hence all edges on $\p$ (if any) are
$A$-edges.  Since the boundary of $\Delta'$ contains no $X$-edges,
$\Delta'$ does not contain $(A,X)$- and $(B,X)$-cells by Lemma
\ref{noannul}. Since the presentation of $H(A)$ satisfies
$C'(\frac1{12})$ (Lemma \ref{lsec}), by the Greendlinger lemma \cite{LS} the
boundary $\partial\pi_1$ of one of the cells in $\Delta'$ shares a subpath of
length at least $\frac9{12}$ of the length of $\partial\pi_1$ with the
boundary of $\Delta'$. Therefore there exists a large section from
$\phi(A\times Y)$ or $\gamma(B\times Y)$ which shares a subpath of at least
$\frac16$ of its length with a product of at most two large sections
from the set $\chi'(A\times X)\cup \beta(A,X)\cup \beta'(A\times
X)\cup \chi(A\times X)$. That is impossible because the set of all
large sections of defining relators satisfies the small cancelation
condition $C'(\frac1{12})$.

The case when $\Delta$ {\em does not satisfy (a')} is completely analogous to the previous  case. Only instead of joining $\Pi$ and $\pi$ together, we would need to cut off $\pi$ from $\Delta$ reducing the number of cells.

Suppose now that $\Delta$ {\em does not satisfy (b)}. Let $e_1,e_2$ (resp. $e_1', e_2'$) be the start (resp. end) edges of $\xxx_1, \xxx_2$. We
suppose without loss of generality that the labels of $e_1$, $e_2$
are positive. We also assume that $e_1$ precedes $e_2$ on $\partial\Pi$. Then $e_2'$ precedes $e_1'$ on $\partial\Pi'$. Let $\p$ (resp. $\p'$) be the subpath between $e_1$, $e_2$ (resp.  $e_2'$ and $e_1'$) on $\partial\Pi$ (resp. $\partial\Pi'$). By (c), the subdiagram bounded by the medians of $\xxx, \xxx'$, $\p,\p'$  does not contain $\N$-cells. By (a) $e_2', e_1'$
are consecutive $X$-edges of $\partial\Pi'$ (that is there are no $X$-edges of $\partial\Pi$ between them).

Since $\partial\Delta'$ does not contain $X$-edges, $\Delta'$ does
not contain  as $(A,X)$- and $(B,X)$-cells by Lemma \ref{noannul}.
Thus $\Delta'$ is a reduced diagram over the presentation involving
only $(A,Y)$- and $(B,Y)$-cells. This and the small cancelation
property $C'(\frac1{12})$ imply that $\Delta'$ does not have cells.
Hence if, say, $\xxx$ is not empty, a cell in $\xxx$ shares more
than $1/12$ of one of its large sections with a large section of
another cell of $\xxx$ or with a large section of a cell in $\xxx'$.
By Lemma \ref{lsec}, these two cells cancel, a contradiction with
the assumption that $\Delta$ is reduced. Hence $\xxx,\xxx'$ are
empty, a contradiction.

Suppose that $\Delta$ {\em does not satisfy (b')}. The only difference with the previous case is that the end edges of $\xxx$ and $\xxx'$ are consecutive $x$-edges on $\partial\Delta$ in this case, and we need to use (a') instead of (a). Since the label of $\partial\Delta$ is a product of words of the form $u_av_a$ in this case, the proof proceeds the same way as in the previous case.


Finally suppose that $\Delta$ {\em does not satisfy (c)}. Assume $\Lab(\q)$ does not have letters from $X$ (the other case is similar).

By (a), the graph $\mho$ from Notation \ref{not6} does not have a vertex connected to itself. By the classical result of Heawood  \cite{He} $\mho$
has a vertex of degree at most 5.

Suppose that $\Delta$ contains two $\N$-subdiagrams $\Pi$ and $\Pi'$
such that at least
one fifth of all $X$-bands starting on $\partial\Pi$ end
on $\partial\Pi'$. By (b) then $\Pi$ and
$\Pi'$ share large $X$-portions of their boundaries. Therefore the boundaries of the union of
$\Pi$ and $\Pi'$ have labels from $E$ (by Lemma \ref{l:mu}), which contradicts the
maximality of $\Pi'$ and $\Pi'$ (as vertices of $\mho$).

Therefore there is an $\N$-cell in $\Delta$ and at least one fifth of all $X$-bands starting on $\partial\Pi$ end on $\partial\Delta$. Therefore there exists a subpath $\mathfrak{v}$ on the boundary whose label is an $X$-subblock $\mathfrak{w}$ and at least $|\mathfrak{w}|/5$ of consecutive $X$-bands starting on $\mathfrak{w}$ end on $\partial\Delta$. By (a'), the end edges of these bands are consecutive edges of $\partial\Delta$. Therefore the end edges of these bands form a path $\p$ whose label is $\Lab(\mathfrak{w})$. By (b') all these $X$-bands are empty, so $\Pi$ shares a big portion of its boundary with $\p$, a contradiction.
\endproof

\begin{lemma}\label{l:cep} The subgroup $E$ of $H(A)$ satisfies the {\em congruence extension property},  that is the intersection of $\ncn$ with $E$ is $\N$ (i.e. $E/\N$
naturally embeds into $H(A)/\ll\!\! \N\!\!\gg$) for every normal subgroup $\N$ of $E$.
\end{lemma}

\proof Suppose that $W\in E$ belongs to $\ncn$. Then we may assume that $W$ is a product of words $u_av_a$ and their inverses, and there exists a \vk diagram $\Delta$ over the presentation that
consists of the defining relations of $H(A)$ and all words  from
$\N$ and boundary label $W$. We need to show that $W\in \N$. By contradiction assume that $W$ is a counterexample, and $\Delta$ is minimal (with respect to the number of cells) diagram for $W=1$ for all counterexamples $W$. If $\Delta$ does not contain $\N$-cells, then $\Delta$ is a diagram over the presentation of $H(A)$ and so it does not contain cells by Lemma \ref{l:free}.

Suppose that $\Delta$ contains $\N$-cells. Consider the graph $\mho$ as in Notation \ref{not6}. By \cite{He} there exists a vertex $\Pi$ in $\mho$ of degree at most 5. Therefore at least $\frac15$ of consecutive $X$-bands starting on $\partial\Pi$ end either on the boundary of another $\N$-subdiagram $\Pi'$ or on $\partial\Delta$. The first possibility would mean, by Lemma \ref{lp1.1}, that $\Pi$ and $\Pi'$ share large $X$-portions of their boundaries.
Hence by Lemma \ref{l:mu} the union of $\Pi$ and $\Pi'$ is a holey $E$-subdiagram of $\Delta$ which contradicts the maximality of $\Pi, \Pi'$.
If the second possibility occurs then
by Lemma \ref{lp1.1} (b'), $\Pi$ shares a large $X$-portion of its boundary with the boundary of $\Delta$. Then (again by Lemma \ref{l:mu}) we can cut $\pi$ off $\Delta$ and produce a smaller diagram $\Delta'$ with boundary label in $E$. By the minimality assumption for $\Delta$, $\Lab(\partial\Delta')$ must belong to $\N$. But then $W\in \N$ as well, a contradiction.
\endproof

\begin{lemma} \label{p1} For every normal subgroup $\N$ of $E$,
$\la A\cup B\ra\cap \ncn=\{1\}$, $\la X\cup Y\ra \cap \ncn =\{1\}$,
\end{lemma}
\proof
Suppose that the boundary label of $\Delta$ is a reduced group word in
$A\cup B$. By Lemma \ref{lp1.1} (c), $\Delta$ has no $\N$-cells. If
$\Delta$ has $(A,X)$-cells or $(B,X)$-cells, then it has $X$-bands,
which must start and end on the boundary of $\Delta$ by Lemma
\ref{noannul}, a contradiction.  Hence $\Delta$ consists of $(B,Y)$-cells. These relations satisfy the small cancelation
condition $C'(\frac1{12})$. By the Greendlinger lemma \cite{LS} there
exists a cell $\pi$ in $\Delta$ such that at least $\frac9{12}$ of
$\partial\pi$ is contained in $\partial\Delta$. But $\frac9{12}$ of
the boundary of a $(B,Y)$-cell contains $Y$-edges, a
contradiction.

Suppose that the boundary label of $\Delta$ is a reduced word in
$X\cup Y$. As in the previous paragraph, $\Delta$ does not contain
$\N$-cells, $(A,X)$- and $(A,Y)$-cells. Therefore $\Delta$ contains
only $(B,X)$- and $(B,Y)$-cells. Then the $B$-bands must start and
end on the boundary of $\Delta$ by Lemma \ref{noannul}, a contradiction.
\endproof

\begin{lemma}\label{p2} For every word $U$ in $\la A\cup B\ra$ there exists at most one pair of reduced words $W\in E$ and $V\in \la X\cup Y\ra$ such that $W=UV$ in $H(A)$.
\end{lemma}

\proof Indeed, if $W=UV$ and $W'=UV'$, then $W\iv W'=V\iv V'\in \la X\cup Y\ra\cap E=\{1\}$ by Lemma \ref{p1} (take $\N=E$ there). Hence $W=W', V=V'$ in the free group. \endproof

\begin{notation}\label{n:1} Consider any group $\Gamma=\la A\mid R\ra$. We always assume that the presentation of $\Gamma$ is {\em positive}, i.e. consists of positive words in the alphabet $A$. For this, we assume that $A$ is divided into two parts of equal sizes $A^+, A^-$ with a bijection $\bar{\quad}  \colon A^+\leftrightarrow A^-$, and $R$ contains all relations of the form $a\bar a=1, a\in A$. Clearly every group with combinatorially aspherical presentation has a positive combinatorially aspherical presentation.

Let $A'$ be a copy of $A$, $\e$ be a bijection $A\to A'$, $R'=\e(R)$ be the set of words $R$ rewritten in $A'$,  $\Gamma'=\la A'\mid R'\ra$.
\end{notation}

Let $H'(\Gamma)$ be the group $\la \Gamma', H(A), q\mid q\e(a)q\iv =
\mu(a), a\in A'\ra$, that is $H'(\Gamma)$ is the factor-group of the
free product $\Gamma'*H(A)* {\mathbb Z}$ (where the copy of $\Z$ is generated by $q$) by the {\em conjugacy
$q$-relations} $q\e(a)q\iv = \mu(a)$.

\begin{lemma}\label{l64.1} Let $\N$ be the normal subgroup of $H(A)$ generated (as a normal subgroup) by $\mu(R)$, $H'=H(A)/\N$. The group $H'(\Gamma)$
is isomorphic to the HNN-extension of the free product $H'*\Gamma$ with free letter $q$ and associated subgroups $\Gamma'$ and $\la \mu(A)\ra \N/\N$.
\end{lemma}

\proof Indeed it is enough to establish that the groups $\Gamma$ and $\la \mu(A)\ra \N/\N$ are isomorphic. But that follows from the congruence extension property satisfied by the subgroup $E=\la\mu(A)\ra$ in $H(A)$ by Lemma \ref{l:cep}.
\endproof

Relations of $H'(\Gamma)$ and cells in \vk diagramm corresponding to the relations from $R$ of $\Gamma$ will be called $\Gamma$-relations and $\Gamma$-cells.

\begin{lemma} \label{lm:67} Suppose that every disc subdiagram $\Delta'$ of a diagram $\Delta$ whose boundary consists of $A'$-edges consists of $\Gamma$-cells. Suppose also that $\partial\Delta$ does not consist of $A'$-edges and does not contain $q$-edges. Then every $\Gamma$-cell in $\Delta$ is in the inside diagram of a $q$-annulus of $\Delta$.
\end{lemma}
\proof Indeed, let $\Delta_1$ be a maximal holey subdiagram of $\Delta'$ that consists of $\Gamma$-cells. Every boundary component of $\Delta_1$ consists of $A'$-edges. By our assumption then $\Delta_1$ does not have holes and is a disc subdiagram.  At least one edge on $\partial\Delta_1$ must belong to the boundary of a $q$-cell (since $\partial\Delta'$ does not consist of $A'$-edges). The $q$-band containing that cell is an annulus since $\partial\Delta$ does not contain $q$-edges. That annulus cannot have outside boundary consisting of $A'$-edges by our assumption (see Figure \ref{p92} a) ). So its inside diagram must contain $\Delta_1$ (see Figure \ref{p92} b) ).

\begin{figure}[ht]
\unitlength .45mm 
\linethickness{0.4pt}
\ifx\plotpoint\undefined\newsavebox{\plotpoint}\fi 
\begin{picture}(300.75,85.625)(0,5)
\qbezier(6,46.25)(-4,78.125)(32,65.5)
\qbezier(211.5,44.75)(201.5,76.625)(237.5,64)
\qbezier(32,65.5)(44.125,61.5)(59.75,63.5)
\qbezier(237.5,64)(249.625,60)(265.25,62)
\qbezier(59.75,63.5)(73.125,66.25)(68,46)
\qbezier(265.25,62)(278.625,64.75)(273.5,44.5)
\qbezier(68,46)(61.625,25.75)(31.75,33.5)
\qbezier(273.5,44.5)(267.125,24.25)(237.25,32)
\qbezier(31.75,33.5)(8.625,39.125)(6,46.25)
\qbezier(237.25,32)(214.125,37.625)(211.5,44.75)
\put(68.25,48){\vector(1,4){.176}}\multiput(67.5,45.25)(.0833333,.3055556){18}{\line(0,1){.3055556}}
\put(273.75,46.5){\vector(1,4){.176}}\multiput(273,43.75)(.0833333,.3055556){18}{\line(0,1){.3055556}}
\put(87.25,77.5){\vector(2,1){.176}}\qbezier(68.75,51.25)(77.125,90.25)(126,78.25)
\put(146,67){\vector(4,-3){.176}}\qbezier(126,78.25)(145.5,73.625)(167,42.5)
\put(153.5,31.25){\vector(-3,-1){.176}}\qbezier(167,42.5)(173.75,33)(99.5,16.5)
\put(71.25,19.5){\vector(-4,3){.176}}\qbezier(99.5,16.5)(58.875,8.25)(67.75,45)
\put(119.25,74.25){\vector(1,0){.176}}\qbezier(74.75,55)(124,99.625)(154.25,42.75)
\put(143.25,37.25){\vector(-3,-1){.176}}\qbezier(154.25,42.75)(156.25,39.75)(106.25,26.75)
\qbezier(75.25,55.25)(63.625,16.375)(106.5,27)
\put(70.25,44.25){\vector(-4,1){.176}}\multiput(72.75,43.75)(-.416667,.083333){12}{\line(-1,0){.416667}}
\put(72.125,53){\vector(-1,0){.176}}\put(74.75,53){\line(-1,0){5.25}}
\put(75.875,60.25){\vector(-4,1){.176}}\multiput(79.75,59.75)(-.645833,.083333){12}{\line(-1,0){.645833}}
\put(82,67.625){\vector(-1,1){.176}}\multiput(86,64)(-.09302326,.08430233){86}{\line(-1,0){.09302326}}
\put(90.75,73.125){\vector(-1,2){.176}}\multiput(93,68.5)(-.08333333,.1712963){54}{\line(0,1){.1712963}}
\put(70.5,30.875){\vector(-1,-1){.176}}\multiput(74,34)(-.09333333,-.08333333){75}{\line(-1,0){.09333333}}
\put(75.875,22.875){\vector(-1,-1){.176}}\multiput(79.5,27)(-.08430233,-.09593023){86}{\line(0,-1){.09593023}}
\put(94.25,75){\makebox(0,0)[cc]{$q$}}
\put(82,20.75){\makebox(0,0)[cc]{$q$}}
\put(66.25,50){\makebox(0,0)[cc]{$a$}}
\put(84.75,48.5){\makebox(0,0)[cc]{$\mu(a)$}}
\put(271.75,80.25){\vector(-1,2){.176}}\qbezier(274.5,49.25)(280.5,91.125)(251.5,89.5)
\put(218.75,85.75){\vector(-3,-1){.176}}\qbezier(251.5,89.5)(213.375,93.5)(196.75,66.5)
\put(187.75,32.25){\vector(1,-4){.176}}\qbezier(196.75,66.5)(170.125,21.25)(214,20)
\put(239,18.5){\vector(1,0){.176}}\qbezier(214,20)(241.375,18.375)(259.25,17.25)
\put(267.25,23.75){\vector(2,3){.176}}\qbezier(259.25,17.25)(268.375,16.75)(273,44.25)
\put(278.25,83){\vector(-1,3){.176}}\qbezier(280,46)(285.375,95.25)(262.25,95.5)
\put(219.5,91){\vector(-4,-1){.176}}\qbezier(262.25,95.5)(209.75,96.875)(196.25,74.75)
\put(181.5,46.75){\vector(-1,-2){.176}}\qbezier(196.25,74.75)(175.5,40.25)(178.75,31.75)
\put(193,19.25){\vector(2,-1){.176}}\qbezier(178.75,31.75)(185.125,15.625)(223,14)
\put(261,13.25){\vector(4,1){.176}}\qbezier(223,14)(272.5,7.75)(276,23.5)
\put(275.5,39.5){\vector(-2,1){.176}}\multiput(278.75,38)(-.18055556,.08333333){36}{\line(-1,0){.18055556}}
\put(277.125,48.5){\vector(-4,-1){.176}}\multiput(280,49)(-.479167,-.083333){12}{\line(-1,0){.479167}}
\put(278.5,60.25){\vector(-4,-1){.176}}\multiput(281,60.75)(-.416667,-.083333){12}{\line(-1,0){.416667}}
\put(277.5,71.125){\vector(-4,-1){.176}}\multiput(280.25,71.75)(-.3666667,-.0833333){15}{\line(-1,0){.3666667}}
\put(272.375,86.25){\vector(-2,-1){.176}}\multiput(276,88.25)(-.15104167,-.08333333){48}{\line(-1,0){.15104167}}
\put(272.625,27.75){\vector(-3,1){.176}}\multiput(276.5,26.5)(-.2583333,.0833333){30}{\line(-1,0){.2583333}}
\put(266.875,17.625){\vector(-3,2){.176}}\multiput(270,15.75)(-.13888889,.08333333){45}{\line(-1,0){.13888889}}
\put(270.5,44.75){\makebox(0,0)[cc]{$a$}}
\put(291.75,43.75){\makebox(0,0)[cc]{$\mu(a)$}}
\put(270.5,89){\makebox(0,0)[cc]{$q$}}
\put(269.2,20.1){\makebox(0,0)[cc]{$q$}}
\put(33.75,49.75){\makebox(0,0)[cc]{$\Delta_1$}}
\put(240.5,48.5){\makebox(0,0)[cc]{$\Delta_1$}}
\put(78,6){\makebox(0,0)[cc]{a)}}
\put(231.5,6.5){\makebox(0,0)[cc]{b)}}
\multiput(276,23.75)(.08333333,.48958333){48}{\line(0,1){.48958333}}
\put(279.375,44.75){\vector(1,4){.176}}\multiput(278.5,40.75)(.0833333,.3809524){21}{\line(0,1){.3809524}}
\end{picture}

\caption{}\label{p92}
\end{figure}
\endproof

\begin{lemma} \label{l64} Every diagram over the presentation of $H'(\Gamma)$ whose boundary consists of $A'$-edges is combinatorially homotopic to a diagram consisting of $\Gamma$-cells.
The subset $X\cup Y$   of $H'(\Gamma)$ freely generates a free subgroup.
\end{lemma}
\proof  Let $H''$ be the group given by the presentation of $H'(\Gamma)$ without the $\Gamma$-relations. Then $H''$ is the HNN-extension of the free product $\la A'\ra*H(A)$ with free letter $q$ and associated subgroups $\la \e(A)^{(1)}\ra$ and $\la \mu(A)^{(1)}\ra$ which are both free and freely generated by their respective generating sets by the definition of $\e$ and $\mu$. Hence by Lemmas \ref{cch1}, \ref{cch2}, the presentation of $H''$ is aspherical. The standard properties of HNN-extensions and free products show that the set $A'$ freely generates a free subgroup in $H''$.

With every diagram $\Delta$ over the presentation of $H'(\Gamma)$ we associate its  {\em 2-weight}, a pair of numbers $(m,n)$ where $n$ is the number of $\Gamma$-cells in $\Delta$ and $m$ is the number of other cells in $\Delta$. We order all pairs $(m,n)$ lexicographically.

Now consider any reduced \vk diagram $\Delta$ over the presentation of $H'(\Gamma)$ whose boundary consists of $A'$-edges. Suppose that $\Delta$ is not combinatorially homotopic to a diagram containing only $\Gamma$-cells
and is a smallest 2-weight diagram with this property and boundary label in $A'$.

If $\Delta$ does not contain $\Gamma$-cells, then the diagram is combinatorially homotopic to the trivial diagram because it is a diagram over the presentation of $H''$.

Suppose that $\Delta$ contains $\Gamma$-cells.

If $\Delta$ does not contain conjugacy $q$-cells, then the boundary of every maximal holey subdiagrams of $\Delta$ consisting of $\Gamma$-cells must coincide with the boundary of $\Delta$, therefore $\Delta$ consists of $\Gamma$-cells.
Thus we can assume that $\Delta$ contains conjugacy $q$-cells. Hence $\Delta$ has a $q$-annulus $\aaa$ consisting of conjugacy $q$-cells. Then the label of the outer boundary of $\aaa$ is either a word in $A'$ or is from $E$. Suppose that the first option holds. If $\partial\Delta$ is not the outer boundary of $\aaa$, we can use the minimality of $\Delta$ and conclude that the diagram bounded by the outer boundary of $\aaa$ is combinatorially homotopic to a diagram consisting of $\Gamma$-cells. That would reduce the number of non-$\N$-cells in $\Delta$ and the 2-weight of $\Delta$ which is impossible. Thus in the first case $\partial\Delta$ is the outer boundary of $\aaa$. The boundary label of the inside subdiagram $\Delta'$ bounded by the median of $\aaa$ is from $E$. If the second option is true, we can take $\Delta'$ to be the subdiagram bounded by the outer boundary of $\aaa$. Note that in both cases $\Delta'$ contains cells.

If $\Delta'$ does not contain conjugacy $q$-cells, then it cannot contain $\Gamma$-cells either (otherwise consider a maximal holey subdiagram consisting of $\Gamma$-cells), and so it is a diagram over the presentation of $H(A)$. Since the group $E=\la\mu(A)\ra$ in $H(A)$ is free (by Lemma \ref{l:free}), we would have that the label of $\Delta'$ is freely trivial, so $\aaa$ contains two neighbor cells that cancel. Hence we can assume that $\Delta'$ contains conjugacy $q$-cells and $q$-annuli. By the minimality of $\Delta$, the external boundary of each of these annuli must have label from $E$ and the internal boundary must consist of $A'$-edges. By Lemma \ref{l64.1}, the label of the external boundary of any $q$-annulus in $\Delta$ must then be a word from $\N$.

By Lemma \ref{lm:67}, every $\Gamma$-cell in $\Delta'$ is in the inside subdiagram of a $q$-annulus.
Maximal subdiagrams of $\Delta'$ bounded by $q$-annuli will be called $q$-subdiagrams, the $q$-annulus bounding an $q$-subdiagram is called the {\em main} $q$-annulus of the subdiagram. Using diamond moves, we can make the boundary of every $q$-subdiagram reduced. Note that all $q$-subdiagrams are $E$-subdiagrams of $\Delta$ as in Notation \ref{not6}.

The diagram $\Delta'$ is tesselated by $q$-subdiagrams and cells corresponding to the relations of $H(A)$ so the situation is the same as the one considered in Notation \ref{not6} since we can view $\Delta'$ as a reduced diagram over the presentation consisting of $\N$-relations and relations from $H(A)$.

As in Notation \ref{not6} consider the graph $\mho$ whose vertices are all $q$-subdiagrams of $\Delta'$ bounded by $q$-annuli, and two vertices are connected if there is an $X$-band connecting them. Since $\mho$ is a planar graph without multiple edges and vertices adjacent to themselves (by Lemma \ref{lp1.1}, (a) ), there must be (by \cite{He}) a vertex of degree at most 5. Hence either at least $\frac3{12}$ of consecutive $X$-bands starting on one of the $X$-blocks of the boundary of a $q$-subdiagram $\Pi_1$  end on an $X$-block of the boundary of another $q$-subdiagram $\Pi_2$ or $\frac3{12}$ of the $X$-bands starting on one of the $X$-blocks of the boundary of a $q$-subdiagram $\Pi_1$ end on one of the $X$-blocks of the boundary of $\Delta'$. By Lemma \ref{lp1.1} (b) the $X$-bands connecting these blocks are empty. Therefore $\Pi_1$ shares a large $X$-portion of its boundary either with another $q$-subdiagram $\Pi_2$ or with the boundary of $\Delta'$. Applying Lemma \ref{l:mu}, we conclude that in the first case a cell of the main $q$-band of $\Pi_1$ cancels with a cell of the main $q$-band of $\Pi_2$,  and, in the second case, a cell of the main $q$-band of $\Pi_1$ cancels with a cell in $\aaa$. This contradiction completes the proof of the first statement of the lemma.

The second statement immediately follows from Lemmas \ref{l64.1} and \ref{p1}.
\endproof

\section{The main construction}

Let $\gi=\la A\mid R\ra$ be a finitely generated recursively
presented group, so that $R$ is a recursive set of defining
relations.  In this section, we describe an embedding of $\gi$ into
a finitely presented group $\da$.

\subsection{Preliminaries on $S$-machines}

\subsubsection{A definition of $S$-machines} Following \cite{SBR, BORS} we shall give two (equivalent) definitions of $S$-machines (a slightly different definition can be found in \cite{OSnlogn}, \cite{OScol}, \cite{SaICM} and other papers but since we are going to use results of \cite{SBR}, we give definitions that are closer to \cite{SBR}).
Let $n$ be a natural number. A {\em hardware} of an $S$-machine is
a pair $(Z,Q)$ where $Z$ is an $m-1$-vector of
(not necessary disjoint) sets $Z_i$ of {\em tape letters}, $Q$ is a $m$-vector of disjoint sets
$Q_i$ of {\em state letters}. The sets $\bigcup Q_i$ and $\bigcup Z_i$ are also disjoint.

The {\em language of
admissible words} is $L(\sss)=Q_1F(Z_1)Q_2...F(Z_{n-1})Q_m$ where
$F(Z_j)$ is the language of all reduced group words in the
alphabet $Z_j\cup Z_j\iv$.

If $1\le i\le j\le m$ and $W= q_1u_1q_2...u_{m-1}q_m$ is an
admissible word, $q_i\in Q_i, u_i\in (Z_i\cup Z_i\iv)^*$ then the
subword $q_iu_i...q_j$ of $W$ is called the $(Q_i,Q_j)$-{\em subword} of
$W$.

An $S$-machine with hardware $\sss$ is a rewriting system. The
objects of this rewriting system are all admissible words.

The rewriting rules, or {\em $S$-rules}, have the following form:
$$[U_1\to V_1,...,U_n\to V_n]$$
where the following conditions hold:
\begin{itemize}
\item Each $U_i$ is a subword of an admissible word starting with
a $Q_\ell$-letter and ending with a $Q_r$-letter $\ell=\ell(i), r=r(i)$.
\item If $i<j$ then $r(i)<\ell(j)$.
\item Each $V_i$ is also a subword of an admissible word whose $Q$-letters
belong to $Q_{\ell(i)}\cup...\cup Q_{r(i)}$ and which contains a $Q_{\ell(i)}$-letter
and a $Q_{r(i)}$-letter.
\item $V_1$ must start with a $Q_1$-letter and
$V_n$ must end with a $Q_{m}$-letter.
\end{itemize}

To apply an $S$-rule to a word $W$ means to replace simultaneously subwords
$U_i$ by subwords $V_i$, $i=1,...,n$. In particular, this means that our rule is
not applicable if one of the $U_i$'s is not a subword of $W$. The following convention
is important:

{\bf After every application of a rewriting rule,
the word is automatically reduced.}

With every $S$-rule $\tau$ we associate the inverse $S$-rule $\tau\iv$
in the following way: if $$\tau=[U_1\to x_1V_i'y_1,\ U_2\to x_2V_2'y_2,...,
U_n\to x_nV_n'y_n]$$ where $V_i'$ starts with a $Q_{\ell(i)}$-letter and
ends with a $Q_{r(i)}$-letter, then $$\tau\iv=[V_1'\to x_1\iv U_1y_1\iv,\
V_2'\to x_2\iv U_2y_2\iv,...,V_n'\to x_n\iv U_ny_n\iv].$$

It is clear that $\tau\iv$ is an $S$-rule,  $(\tau\iv)\iv=\tau$, and that
rules $\tau$ and $\tau\iv$ cancel each other (meaning
that if we apply $\tau$ and then $\tau\iv$, we return to the
original word).

The following convention is also important:

{\bf We always assume that an $S$-machine is symmetric, that is if an
$S$-machine contains a rewriting rule $\tau$, it also contains the rule
$\tau\iv$.} Among any pair of mutually inverse rules we pick one which we call {\em positive rule}, the other rule is called {\em negative}.

We define the history of
a computation of an $S$-machine as the sequence (word)
of rules used in this computation.
A computation is called {\em reduced} if the history of this computation is
reduced, that is if two mutually inverse rules are never
applied next to each other.

For some $S$-machines we distinguish {\em input} and {\em stop} admissible words.

\begin{rk} \label{rk:45}
We always assume that:
\begin{itemize}
\item There is only one stop word, it does not contain tape letters;
\item If we remove tape letters from any two input words, we obtain the same word $q_1\ldots q_m$ which we call the {\em input base};
\item A negative rule cannot apply to any input word;
\item There exists only one positive rule that applies to an input configuration, it has the form $[q_1\ldots q_{m}\to \tilde q_1\ldots \tilde q_{n}]$ where $\tilde q_i\ne q_i$ for some $i$. This rule will be called the {\em transition rule} of the $S$-machine.
\end{itemize}
(The last two properties can be obtained by adding the {\em transition rule} to the $S$-machine. Note that it does not affect any other properties of $S$-machines used in this paper.)
\end{rk}

If an $S$-machine $\sss$ has a stop word $W_0$, then we say that an admissible word $W$ is {\em accepted} by $\sss$ if there exists a computation of $\sss$ starting with $W$ and ending with $W_0$. That computation is called an {\em accepting computation} for $W$.

\subsubsection{Recursively enumerable sets and $S$-machines}

One of the main results of \cite{SBR} implies that for every recursively enumerable set of words $L$ in an alphabet $A$ there exists an $S$-machine $\sss$ recognizing $L$ in the following sense.\footnote{More precisely, the $S$-machine $\sss$ we use is the $S$-machine from \cite{SBR} without $\alpha$- and $\omega$-sectors. These sectors are needed in \cite{SBR} only to control the Dehn function of the resulting group. The facts from \cite{SBR} that we are using here remain true. We could use literally the same $S$-machines as in \cite{SBR} but it would make our proof unnecessarily more cumbersome because the input words would contain powers of $\alpha$ and $\omega$.}

\begin{prop} \label{mach1} The $S$-machine $\sss$ has a stop word $W_0$. For every positive word $u$ in the alphabet $A$, there exists an input word $I(u)$ and
\begin{enumerate}

\item\label{11} The input word $I(u)$ has the form $q_1uq_2\ldots q_m$, $q_i\in Q_i, i=1,\ldots,m$;
\item\label{12} A word $u$ belongs to $L$ if and only if $I(u)$ is accepted by $\sss$;
\item\label{13} If $I(u)$ is accepted by $\sss$, then there exists only one reduced computation accepting $I(u)$.
\item\label{14} There is only one reduced computation connecting $W_0$ with itself, the empty one.
\end{enumerate}
\end{prop}
\proof Take a deterministic Turing machine $T$ recognizing $L$. Convert it into a symmetric Turing machine $T'$ using \cite[Lemma 3.1]{SBR}. Then use \cite[Proposition 4.1]{SBR} to convert $T'$ into an $S$-machine.\footnote{The conversion in \cite{SBR} was very complicated because we needed to control the speed of the $S$-machine. Since we do not care about the speed of the $S$-machines in this paper, we could use simpler but (exponentially) slower $S$-machines constructed in \cite{OSnlogn}. Still the $S$-machine from \cite{SBR} is useful for us because we can use some facts about it proved in \cite{SBR}.} The fact that this $S$-machine satisfies the conditions (1) and (2) of the proposition immediately follows from \cite[Lemma 3.1 and Proposition 4.1]{SBR}.\endproof

\subsubsection{$S$-machines as HNN-extensions of free groups}
\label{ss:hnn}

Another, probably even easier, way to look at $S$-machines is to consider them as multiple HNN-extensions of free groups (see \cite{SaICM, Sasur}). Let $\sss$ be an $S$-machine with the set of tape letters $Z=\cup_{i=1}^{m-1} Z_i$, set of state letters\footnote{$\sqcup$ denotes the disjoint union.} $Q=\sqcup_{i=1}^m Q_i$ and the set of rules $\Theta$. The set of all positive rules is denoted by $\Theta^+$.

The generating set of the group is $Q\cup Z\cup \Theta^+$. The relations are:

$$U_i\theta=\theta V_i,\,\,\,\, i=1,\ldots, n$$
(these relations will be called $(Q,\Theta)$-{\em relations}),

$$ \theta a=a\theta$$
for all $a\in Z$, $\theta\in\Theta^+$ (these relations will be called $(Z,\Theta)$-{\em relations}).

For simplicity and following \cite{SaICM, Sasur}, we shall call this group an $S$-machine too, and denote it by the same letter $\sss$.


\begin{lemma}[Lemma 7.6 \cite{SBR}]\label{noannul3} In any reduced diagram over the presentation of an $S$-machine $\sss$, there are no $\Theta$-annuli, $Q_i$-annuli and $Z$-annuli.
\end{lemma}


\begin{figure}[ht]
\unitlength .8mm 
\linethickness{0.4pt}
\ifx\plotpoint\undefined\newsavebox{\plotpoint}\fi 
\begin{picture}(76.25,87.5)(0,10)
\put(26.875,49.625){\vector(-1,4){.07}}\multiput(32.25,7.5)(-.0336990596,.2641065831){319}{\line(0,1){.2641065831}}
\put(68.625,49.625){\vector(1,4){.07}}\multiput(63.25,7.5)(.0336990596,.2641065831){319}{\line(0,1){.2641065831}}
\put(47.875,7.375){\vector(1,0){.07}}\multiput(32.25,7.25)(3.90625,.03125){8}{\line(1,0){3.90625}}
\put(47.625,91.75){\vector(1,0){.07}}\put(21.25,91.75){\line(1,0){52.75}}
\put(47.875,15.875){\vector(1,0){.07}}\multiput(31.25,15.75)(4.15625,.03125){8}{\line(1,0){4.15625}}
\put(47.75,24.25){\vector(1,0){.07}}\multiput(30.25,24)(2.3333333,.0333333){15}{\line(1,0){2.3333333}}
\put(47.875,82.5){\vector(1,0){.07}}\put(22.75,82.5){\line(1,0){50.25}}
\put(46.25,9.25){\makebox(0,0)[cc]{$W_1$}}
\put(46.75,18.75){\makebox(0,0)[cc]{$W_2$}}
\put(45.75,27.25){\makebox(0,0)[cc]{$W_3$}}
\put(45.75,85.75){\makebox(0,0)[cc]{$W_{g-1}$}}
\put(46.25,94.5){\makebox(0,0)[cc]{$W_g$}}
\put(29.25,11.5){\makebox(0,0)[cc]{$\theta_1$}}
\put(26.75,19.75){\makebox(0,0)[cc]{$\theta_2$}}
\put(16,86.75){\makebox(0,0)[cc]{$\theta_g$}}
\put(67.25,11.75){\makebox(0,0)[cc]{$\theta_1$}}
\put(67.75,21){\makebox(0,0)[cc]{$\theta_2$}}
\put(76.25,87.5){\makebox(0,0)[cc]{$\theta_g$}}
\put(44.5,54.25){\makebox(0,0)[cc]{$\ldots$}}
\end{picture}
\caption{}\label{p11}
\end{figure}

Consider now an arbitrary computation $C = (W_1, W_2, \ldots, W_g$) of an $S$-machine $\sss$
with a history word $h$.
With every $i=1,\ldots,g-1$
we associate the $\Theta$-band $\ttt_i$ with the boundary label $\theta_i\iv W_i\theta_iW_{i+1}\iv$ where $\theta_i$ is the $i$-letter in $h$.
We can ``concatenate" all these bands to obtain a \vk diagram with boundary of the form $\p_1\iv \p_2\p_3\p_4\iv$ where $\p_1, \p_4$ are labeled by $h$, $\Lab(\p_2)\equiv W_1$, $\Lab(\p_4)\equiv W_g$ (see Figure \ref{p11}). This diagram is called a {\em computational sector} corresponding to the computation $C$.

In general by a {\em sector} we mean a reduced diagram $\Delta$ over the presentation of $\sss$
with boundary divided into four parts,
$\partial(\Delta)=\p_1\iv \p_2\p_3 \p_4\iv$,
such that the following properties hold:
\begin{itemize}
\item $\Lab(\p_1), \Lab(\p_3)$ are reduced group words in $\Theta^+$;
\item $\Lab(\p_2), \Lab(\p_4)$ are admissible words.
\end{itemize}

The following lemma is essentially \cite[Proposition 9.1]{SBR}.

\begin{lemma}\label{p9.1} Every sector is combinatorially homotopic (even without insertions of cancelable cells) to a computational sector corresponding to a reduced computation connecting $\Lab(\p_2)$ with $\Lab(\p_4)$. The history of that computation is $\Lab(\p_1)\equiv \Lab(\p_3)$.
\end{lemma}

\subsection{Four $S$-machines}

\subsubsection{The $S$-machines $\sss_1$, $\sss_2$}
By Lemma \ref{l0}, the set of words $U(R)=\{U(\mu(r))\mid r\in R\}$ is
recursive (see the notation in Remark \ref{rk:0}). Let $\sss_1$ be an $S$-machine recognizing
that set of words and satisfying the conditions of
Lemma \ref{mach1}.

Let $\sss_2$ be an $S$-machine recognizing the set
$R$ written in the copy $\hat A$ of $A$ and satisfying conditions of Lemma \ref{mach1}. By splitting the state
letters (i.e. replacing in admissible words and in all rules a state
letter $q$ by a product of new state letters $q'q''$, we can assume
that input bases (see Remark \ref{rk:45}) of $\sss_1$ and $\sss_2$ are the same words, and the set of state letters of $\sss_i$ has $m$ parts, $i=1,2$.

The tape alphabet of $\sss_1$ is $Z(\sss_1)=Z_1\cup\ldots\cup Z_{m-1}$ (we shall assume, without loss of generality, that $Z_1=\emptyset$), the
state alphabet is $Q(\sss_1)=Q_1\sqcup\ldots\sqcup Q_m$. The tape alphabet of
$\sss_2$ is $Z(\sss_2)=\hat Z_1\cup\ldots \cup\hat Z_{m-1}$ (we shall assume that $\hat Z_2=\emptyset$), the state
alphabet is $Q(\sss_2)=\hat Q_1\sqcup\ldots \sqcup\hat Q_m$.

We shall assume that $A\cup B$ is contained in $Z_2$, there exists an injective map $\e$ from $A$ to $\hat Z_1$.
The input configuration of $\sss_1$ corresponding to a word $w$ in
the alphabet $A\cup B$ has the form $$I_1(w)\equiv q_1q_2 w
q_3\ldots q_m$$ (the word $w$ is between $q_2$ and $q_3$, there are
no more  tape letters in that word). The input configuration of
$\sss_2$ corresponding to a word $w$ in the alphabet $A$ is
$$I_2(w)\equiv q_1\varepsilon(w)q_2q_3\ldots q_m.$$ Thus the input bases of these $S$-machines are the same:  $q_1\ldots q_m$. We also shall assume that
$Q_i\cap \hat Q_i=\{q_i\}$: $Q_i$ and $\hat Q_i$ do not share
letters except the state letters $q_i$ of the input bases.
The stop words of $\sss_1$ and $\sss_2$, are denoted by $W_0(\sss_1)$ and $W_0(\sss_2)$
respectively.

Pick a number $N\ge 12$ and consider two new $S$-machines $\bar
\sss_1$ and $\bar \sss_2$.

\subsubsection{The $S$-machine $\bar\sss_1$}\label{ss:1}

For each $i=1,\ldots, 2N$ let $Z^{(i)}$ be a disjoint copy of $Z$,
$Q^{(i)}$ be a disjoint copy of $Q$. Then the tape alphabet of
$\bar\sss_1$ is

$$
\begin{array}{l}
Z(\bar\sss_1)=(Z_1^{(1)} \cup \ldots Z_{m-1}^{(1)})\cup
\emptyset \cup \emptyset \cup (Z_{m-1}^{(2)}\cup\ldots\cup Z_1^{(2)})\cup \emptyset \cup\\
\ldots\\
\cup (Z_1^{(2N-1)}\cup\ldots Z_{m-1}^{(2N-1)}) \cup
\emptyset \cup \emptyset \cup (Z_{m-1}^{(2N)}\cup\ldots\cup Z_1^{(2N)}) \cup \emptyset\end{array},$$
the state alphabet is
$$ \begin{array}{l} Q(\bar\sss_1)=(Q_1^{(1)}\sqcup\ldots\sqcup
Q_m^{(1)})\sqcup\{k_1,\bar k_1\} \sqcup (Q_{m}^{(2)}\sqcup\ldots Q_1^{(2)})\sqcup\{t_1\}\sqcup\\
\ldots\\
\sqcup (Q_1^{(2N-1)}\sqcup\ldots\sqcup
Q_m^{(2N-1)})\sqcup\{k_N,\bar k_N\} \sqcup (Q_{m}^{(2N)}
\sqcup\ldots\sqcup Q_1^{(2N)}) \sqcup\{t_N\}\end{array}
$$
The admissible words of the $S$-machine $\bar\sss_1$ are described as follows.
For every word
$W$ in $Z(\sss_1)\cup Q(\sss_1)$ let $W^{(i)}$ be the corresponding copy of that
word in the alphabet $Z(\sss_1)^{(i)}\cup Q(\sss_1)^{(i)}$. Also for every word $W$
let $\overleftarrow{W}$ be word $W$ read from right to left. If
$W$ is an admissible word of $\sss_1$ which is not an input word, then the corresponding
admissible word of $\bar\sss_1$ has the form
$$W(\bar\sss_1)\equiv W^{(1)}\bar k_1\overleftarrow{W}^{(2)}t_1W^{(3)}\bar k_2\overleftarrow{W}^{(4)}\ldots
t_{N-1}W^{(2N-1)}\bar k_N\overleftarrow{W}^{(2N)}t_N.$$ If $I(u)$ is an input word, then the corresponding admissible word of $\bar\sss_1$ is $$\bar I(u)\equiv I(u)^{(1)}k_1\overleftarrow{I(u)}^{(2)}t_1I(u)^{(3)}k_2\overleftarrow{I(u)}^{(4)}\ldots
t_{N-1}I(u)^{(2N-1)}k_N\overleftarrow{I(u)}^{(2N)}t_N$$
(in this case the $k$-letters are $k_i$ instead of $\bar k_i$.)

The stop word of
$\bar\sss_1$, which will be called {\em the first hub} is
$W_0(\bar\sss_1)$, the admissible word of $\bar\sss_1$ corresponding to
the stop word $W_0$ of $\sss_1$. The rules of $\bar\sss_1$ naturally
correspond to rules of $\sss_1$: if $\theta=[U_1\to V_1, \ldots,
U_n\to V_n]$ is a not the transition rule of $\sss_1$ which we shall denote by $\tau_1$ (see Remark \ref{rk:45}), then the corresponding rule
$\bar\theta$ of $\bar\sss_1$ is

\begin{equation}\label{e1} \bar\theta=\left[\begin{array}{l}U_1^{(1)}\to V_1^{(1)},\ldots, U_n^{(1)}\to V_n^{(1)}, k_1\to k_1, \\ \overleftarrow{U_n}^{(2)}\to \overleftarrow{V_n}^{(2)}, \ldots,  \overleftarrow{U_1}^{(2)}\to \overleftarrow{V_1}^{(2)}, t_1\to t_1,\\ \ldots \\
U_1^{(2N-1)}\to V_1^{(2N-1)},\ldots, U_n^{(2N-1)}\to V_n^{(2N-1)}, k_N\to k_N,\\
\overleftarrow{U_n}^{(2N)}\to  \overleftarrow{V_n}^{(2N)}, \ldots,
 \overleftarrow{U_1}^{(2N)}\to \overleftarrow{V_1}^{(2N)}, t_N\to
t_N\end{array}\right].\end{equation} Thus the rule $\bar\theta$
simultaneously executes copies of the rule $\theta$ on all
$(Q_1^{(i)},Q_m^{(i)})$-subwords of an admissible word, $i=1,\ldots,2N$. Essentially $\bar\sss_1$ runs simultaneously $N$ copies of $\sss_1$, which we denote by $\sss_1^{(2i-1)}$, $i=1,\ldots,N$ and $N$ copies of the mirror images $\overleftarrow\sss_1^{(2i)}$ of $\sss_1$, $i=1,\ldots,N$.

If $\theta$ is the transition rule $\tau_1=[q_1\ldots q_m\to\tilde q_1\ldots \tilde q_m]$, then the corresponding {\em transition} rule $\bar\sss_1$ is

\begin{equation}\label{e1.1}\bar\tau_1=\left[\begin{array}{l}q_1^{(1)}\ldots q_m^{(1)}\to \tilde q_1^{(1)}\ldots \tilde q_m^{(1)}, k_1\to \bar k_1, \\ q_m^{(2)} \ldots  q_1^{(2)}\to \tilde q_m^{(2)}\ldots \tilde q_1^{(2)}, t_1\to t_1,\\ \ldots \\
q_1^{(2N-1)}\ldots q_m^{(2N-1)}\to \tilde q_1^{(2N-1)}\ldots \tilde q_m^{(2N-1)}, k_N\to \bar k_N,\\
q_m^{(2N)} \ldots,
 q_1^{(2N)}\to \tilde q_m^{(2N)}\ldots\tilde q_1^{(2N)}, t_N\to
t_N\end{array}\right]\end{equation}
(that is the rule changes all $q_j^{(i)}$ to $\tilde q_j^{(i)}$ and all $k_i$ to $\bar k_i$).

\subsubsection{The $S$-machine $\bar\sss_2$}\label{ss:2}

This machine is constructed using $\sss_2$ in a similar way as
$\bar\sss_1$ is constructed from $\sss_1$, only the set of rules is
constructed somewhat differently.

For each $i=1,\ldots, 2N$ let $\hat Z^{(i)}$ be a disjoint copy of
$\hat Z$, $\hat Q^{(i)}$ be a disjoint copy of $\hat Q$. We identify $\hat Z^{(1)}$ with $\hat Z$. Then the
tape alphabet of $\bar\sss_2$ is

$$\begin{array}{l}Z(\bar\sss_2)=\underbrace{\emptyset\cup\ldots \cup \emptyset}_{m-1}\cup \emptyset\cup\emptyset\cup  (\hat Z_{m-1}^{(2)}\cup\ldots \cup\hat Z_1^{(2)})\cup \emptyset \\
\cup(\hat Z_1^{(3)}\cup\ldots \cup\hat Z_{m-1}^{(3)})\cup \emptyset\cup\emptyset\cup (\hat Z_{m-1}^{(4)}\cup\ldots\cup \hat Z_1^{(4)})\cup \emptyset\\
\ldots\\
\sqcup (\hat Z_1^{(2N-1)}\cup\ldots\cup \hat Z_{m-1}^{(2N-1)})\cup \emptyset
\sqcup\emptyset\cup (\hat Z_{m-1}^{(2N)}\cup\ldots\cup \hat Z_1^{(2N)})\cup \emptyset,\end{array}$$
the state alphabet is
$$ \begin{array}{l} Q(\bar\sss_2)=(\hat Q_1^{(1)}\sqcup\ldots\sqcup \hat Q_m^{(1)})\sqcup\{k_1,\hat k_1\} \sqcup (\hat Q_{m}^{(2)}\sqcup\ldots \hat Q_1^{(2)})\sqcup\{t_1\}\sqcup\\
\ldots\\
\sqcup (\hat Q_1^{(2N-1)}\sqcup\ldots\sqcup\hat  Q_m^{(2N-1)})\sqcup\{k_N,\hat k_N\} \sqcup (\hat Q_{m}^{(2N)}\sqcup\ldots\sqcup \hat Q_1^{(2N)}) \sqcup\{t_N\}.\end{array}
$$
The  description of the admissible words of $\bar\sss_2$ is the
following. For every word $W$ in the alphabet $\hat Z\cup \hat Q$ let $W^{(i)}$
be the corresponding copy of that word in the alphabet $\hat
Z^{(i)}\cup \hat Q^{(i)}$. For every word $W$ in the alphabet $\hat
Z\cup \hat Q$ let $q(W)\equiv W_{\hat Q}$ be the word $W$ with all letters from $\hat
Z$ deleted (the projection of $W$ onto $\hat Q$).
If $W$ is an admissible word of $\sss_2$ but not an input word, then the
corresponding admissible word of $\bar\sss_2$ has the form
$$W(\bar\sss_2)\equiv q(W^{(1)})\hat k_1\overleftarrow{W}^{(2)}t_1W^{(3)}\hat k_2\overleftarrow{W}^{(4)}\ldots
t_{N-1}W^{(2N-1)}\hat k_N\overleftarrow{W}^{(2N)}t_N.$$ The admissible word $\bar I_2(u)$ corresponding to the input word $I_2(u)$ of $\sss_2$ is $$q(I_2(u))k_1\overleftarrow{I_2(u)}^{(2)}t_1I_2(u)^{(3)} k_2\overleftarrow{I_2(u)}^{(4)}\ldots
t_{N-1}I_2(u)^{(2N-1)}k_N\overleftarrow{I_2(u)}^{(2N)}t_N.$$
The stop word of
$\bar\sss_2$ - is the word $W_0(\sss_2)$ which will be called {\em the second hub}.

The rules of $\bar\sss_2$ correspond to
rules of $\sss_2$.  If $\theta=[U_1\to V_1, \ldots, U_n\to V_n]$
is a rule of $\sss_2$ but not the transition rule, which we shall denote $\tau_2$, then the corresponding rule $\bar\theta$ of
$\bar\sss_2$ is (compare with (\ref{e1})):

\begin{equation}\label{e2} \bar\theta=\left[\begin{array}{l}q(U_1^{(1)})\tool q(V_1^{(1)}),\ldots, q(U_n^{(1)})\tool q(V_n^{(1)}), \hat k_1\to \hat k_1, \\ \overleftarrow{U_n}^{(2)}\to \overleftarrow{V_n}^{(2)}, \ldots,  \overleftarrow{U_1}^{(2)}\to \overleftarrow{V_1}^{(2)}, t_1\to t_1,\\ \ldots \\
U_1^{(2N-1)}\to V_1^{(2N-1)},\ldots, U_n^{(2N-1)}\to V_n^{(2N-1)}, \hat k_N\to \hat k_N,\\
\overleftarrow{U_m}^{(2N)}\to  \overleftarrow{V_m}^{(2N)}, \ldots,
 \overleftarrow{U_1}^{(2N)}\to \overleftarrow{V_1}^{(2N)}, t_N\to t_N\end{array}\right].\end{equation}
Thus the $S$-machine $\bar\sss_2$ does not insert of delete any tape letters in the $(Q_1^{(1)},  Q_m^{(2)})$-subwords of admissible words.
If $\theta$ is the transition rule $\tau_2=[q_1\ldots q_m\to \tilde q_1\ldots \tilde q_m]$, then the corresponding {\em transition} rule of $\bar\sss_2$ is
$$\bar\tau_2=\left[\begin{array}{l}q_1^{(1)}\ldots q_m^{(1)}\to \tilde q_1^{(1)}\ldots \tilde q_m^{(1)}, k_1\to \hat k_1, \\ q_m^{(2)} \ldots   q_1^{(2)}\to \tilde q_m^{(2)}\ldots \tilde q_1^{(2)}, t_1\to t_1,\\ \ldots \\
q_1^{(2N-1)}\ldots q_m^{(2N-1)}\to \tilde q_1^{(2N-1)}\ldots \tilde q_m^{(2N-1)}, k_N\to \hat k_N,\\
q_m^{(2N)} \ldots
 q_1^{(2N)}\to \tilde q_m^{(2N)}\ldots \tilde q_1^{(2N)}, t_N\to
t_N\end{array}\right].$$

\subsection{The composition of the $S$-machines and the auxiliary group}\label{ss:3}

Let $A$ be the generating set of our group $\gi$. Let $A\cup B\cup X\cup Y$ be the generating set of $H(A)$. For every $i=1,\ldots, 2N$ we define copies $A^{(i)}$, $B^{(i)}, X^{(i)}, Y^{(i)}$ of the sets $A, B, X, Y$, and the generating set of the group $H^{(i)}(A)$. The relations of $H^{(i)}(A)$ are the relations of $H(A^{(i)})$ if $i$ is odd. If $i$ is even, then the relations of $H^{(i)}(A)$ are obtained from the relations of $H(A^{(i)})$ by reading them from right to left.
Clearly, there exists an isomorphism from $H^{(i)}(A)$ onto $H(A^{(i)})$ for even $i$ given by the map $z\mapsto z\iv$.

We add the generators $X^{(i)}, Y^{(i)}$, relators of all groups $H^{(i)}(A)$ to the already introduced generators and relators. In addition we shall need all relations of the following form

\begin{equation}\label{e3}
z^{(2i-1)}(q_3^{(2i-1)}\ldots q_m^{(2i-1)}k_iq_m^{(2i)}\ldots
q_3^{(2i)})=(q_3^{(2i-1)}\ldots q_m^{(2i-1)}k_iq_m^{(2i)}\ldots
q_3^{(2i)})z^{(2i)}
\end{equation}
for all $z\in X\cup Y$, and all relations of the form

\begin{equation}\label{e4}
\varepsilon(a)^{(2i-1)}q_2^{(2i-1)}=q_2^{(2i-1)}u_a^{(2i-1)}v_a^{(2i-1)}
\end{equation}
for all $a\in A$, $i=1,\ldots,N$, where $u_av_a=\mu(a)$ is the
generator of the subgroup $E$ corresponding to $a\in A$

\begin{equation}\label{e5}
q_2^{(2i)}\varepsilon(a)^{(2i)}=\overleftarrow{v_a}^{(2i)}\overleftarrow{u_a}^{(2i)}q_2^{(2i)}
\end{equation}
for all $a\in A$, $i=1,\ldots,N$.

\subsection{The group $\da$}\label{ss:da}

\begin{df}
The group $\da$ is the group given by the defining relations from
Sections \ref{ss:1}, \ref{ss:2}, \ref{ss:3} and the two hub
relations $W_0(\bar\sss_1)=1$ and $W_0(\bar\sss_2)=1$. Thus the finite presentation of $\da$ consists of:
\begin{itemize}
\item defining relations of $\bar\sss_i$, $i=1,2$ (called the $\bar\sss_i$-{\em relations}),
\item the defining relations of $H^{(i)}(A)$ (called the $H^{(i)}(A)$-{\em relations}),
\item the relations (\ref{e3}) (called the {\em gluing $k_i$-relations}),
\item the relations (\ref{e4}) and (\ref{e5}) (called the {\em conjugacy $q_2^{(i)}$-relations}),
\item the two hub relations.
\end{itemize}
\end{df}

\section{The main result}

We are going to prove the following theorem. Recall that $\varepsilon$ is an injective  map from $A$ to $\hat Z_1$. We shall identify  $\hat Z$ with $\hat Z^{(1)}\subset \da$.

\begin{theorem} \label{th:1} (1)
The map $\varepsilon\colon A\to \da$ extends to a homomorphism
$\gi\to \da$ which will be denoted by $\varepsilon$ as well.

(2) The homomorphism $\varepsilon$ is injective;


(3) If the presentation  $\la A\mid R\ra$ of $\gi$ is combinatorially aspherical, then the (finite) presentation of $\da$ is combinatorially
aspherical.
\end{theorem}

\subsection{The map $\e$ is a homomorphism}

Part (1) of Theorem \ref{th:1} is given by the following lemma.

\begin{lemma}\label{l:phi} For every $r\in R$, $\varepsilon(r)=1$ in $\da$.
\end{lemma}
\proof Let $r\in R$. By the definition of $\sss_1, \sss_2$, $I_1(r)$ is accepted by $\sss_1$, $I_2(r)$ is accepted by $\sss_2$. Therefore $\bar I_1(r)$ is accepted by $\bar\sss_1$ and $\bar I_2(r)$ is accepted by $\bar\sss_2$. Then $\bar I_1(r)=\bar I_2(r)=1$ in $\da$. In fact these equalities are true modulo relations of the $S$-machines $\bar\sss_1, \bar\sss_2$ and the two hub relations; the corresponding \vk diagrams are obtained  from the computational sectors (Figure \ref{p11}) by identifying the left and right sides, and gluing the cells corresponding to the hub relator in the hall of the resulting annulus (see \cite[Section 9]{SBR}).

By the definition of $U(r)$ and $I_1(r)$ there exist words $V^{(i)}(r)$ in the alphabet $X^{(i)}\cup Y^{(i)}$ such that $\mu(r)$  is equal to $U(r)V(r)$ in $H(A)$. Therefore modulo the $H^{(i)}(A)$-relations, $i=1,\ldots,2N$, the
gluing relations and the conjugacy relations the word $\bar I_1(r)$ is equal to the word

\begin{equation}\label{e99} \begin{array}{l} q_1^{(1)}\varepsilon(r)^{(1)}q_2^{(1)}\ldots q_m^{(1)}k_1q_m^{(2)}\ldots q_2^{(2)}\overleftarrow{\varepsilon(r)}^{(2)}t_1\\ \ldots \\ \\q_1^{(2N-1)}\varepsilon(r)^{(2N-1)}q_2^{(2N-1)}\ldots q_m^{(2N-1)}k_Nq_m^{(2N)}\ldots q_2^{(2N)}\overleftarrow{\varepsilon(r)}^{(2N)}t_N.\end{array}\end{equation}
The word (\ref{e99}) can be obtained from the word $\bar I_2(r)$ by inserting the word $\varepsilon(r)^{(1)}\equiv \varepsilon(r)$ between $q_1^{(1)}$ and $q_2^{(2)}$. Since both $\bar I_2(r)$ and the word (\ref{e99}) are equal to 1 in $\da$, the word $\varepsilon(r)$ is equal to 1 in $\da$ as well.
\endproof

\begin{notation}\label{notphi} The proof of Lemma \ref{l:phi} gives a \vk diagram for the relation $\e(r)=1$ for every $r\in R$. That \vk diagram will be denoted by $\Psi(r)$ and will be called {\em standard}. It contains two hubs $\pi_1$ and $\pi_2$ connected by $t_1$-,\ldots, $t_N$-, and $k_1$-,\ldots,$k_N$-bands. It is convenient to make the notation independent on how we enumerate the hubs in $\Psi(r)$, so we assume that $\pi_1$ is a $j$-hub, and $\pi_2$ is a $3-j$-th hub (for some $j=1,2$).  If $i<N$, then the subdiagram bounded by  the $t_i$-band and the $k_{i+1}$-band and $\partial(\ttt_1), \partial(\ttt_2)$ is denoted by $\Psi_i(r)$. The subdiagram $\Psi_i(r)$ contains two maximal transition $\theta$-bands $\aaa_1, \aaa_2$. The medians of these bands divide $\Psi_i(r)$ into three parts $\Psi_i'(r), \Psi_i''(r), \Psi_i'''(r)$ counting from $\pi_1$ to $\pi_2$. The first and the third subdiagrams are computational sectors corresponding to the computation of $\sss_j^{(2i+1)}$ accepting  $I_j^{(i+1)}(U(r))$ and the computation of $\sss_{3-j}^{(2i+1)}$ accepting $I_{3-j}^{(i+1)}(r)$ respectively where $j=1$ or $2$ (recall that $i<N$). The union of $\Psi_i'(r)\cup \aaa_1\cup\Psi''(r)\cup \aaa_2$ is denoted by $\Phi_i(r)$. The complement of that subdiagram in $\Psi(r)$ will be denoted by $\bar\Phi_i(r)$.
\end{notation}

\subsection{Expanded presentation of $\da$}

Let us denote the finite presentation of $\da$ given in \ref{ss:da} by $\pp$, and the presentation obtained from $\pp$ by adding all $\Gamma$-{\em relations} $\e(r), r\in R$ by $\pp'$.

\begin{notation}\label{n:bands}
We shall study diagrams over $\pp$ and $\pp'$. We say that two diagrams $\Delta_1, \Delta_2$ over $\pp'$ are {\em  combinatorially $\pp$-homotopic} if one can transform $\Delta_1$ to $\Delta_2$ using diamond moves, insertion and deletion of cancelable cells corresponding to relations from $\pp$. If $\Delta$ is a diagram over $\pp$, we can consider the following bands in $\Delta$.

\begin{itemize}
\item $\theta$-bands, $\theta\in\bar\Theta_1^+\cup\bar\Theta_2^+$ consisting of $(\theta,a)$, $(\theta,q)$-, and $(\theta,k)$-cells. These bands can start and end on the boundary of $\Delta$.
\item $Q_i^{(j)}$-bands, $i=1, 3,4,\ldots,m$, $j=1,\ldots,2N$ consisting of $(Q,\Theta)$-cells, and gluing $k_j$-cells ($i=3,\ldots,m$). These bands can start (end) on the hub cells and on the boundary of $\Delta$.

\item $Q_2^{(j)}$-bands, $j=1,\ldots,2N$ consisting of $(Q,\Theta)$-cells and conjugacy $q_2^{(i)}$-cells.    These bands can start (end) on the hub cells and on the boundary of $\Delta$.

\item $k_i$-bands consisting of $(k_i,\Theta)$-cells and gluing $k_i$-cells. These bands can start (end) on the hub cells and on the boundary of $\Delta$.

\item $t_i$-bands consisting of $(t_i,\Theta)$-cells. These bands can start (end) on the hub cells and on the boundary of $\Delta$

\item $Z$-bands consisting of $(Z,\Theta)$-cells, $(A^{(i)},X^{(i)})$- and $(A^{(i)},Y^{(i)})$-cells. These bands can start (end) on the boundary of $\Delta$ or on a $(Q,\theta)$-cell or on a conjugacy $q_2^{(i)}$-cell.

\item $X^{(i)}$-bands consisting of $(A^{(i)},X^{(i)})$- and $(B^{(i)},X^{(i)})$-cells. These bands can start (end) on the boundary of $\Delta$ or on the boundary of a $k_i$-gluing cell, or on the boundary of a $q_{2}^{(i)}$-conjugacy cell.
\end{itemize}
\end{notation}

\begin{notation}\label{not90}
Let $\Delta$ be a reduced diagram having two hubs $\pi_1, \pi_2$ connected by a $t_i$-band $\xxx_1$ and a $k_{i+1}$-band $\xxx_2$ where $1 \le i < N$. Suppose that $\Psi$ is bounded by a side of $\xxx_1$, a side of $\xxx_2$ and parts of $\partial\pi_1$, $\partial\pi_2$, and does not contain hubs except for $\pi_1, \pi_2$. We shall call such a diagram a {\em $2$-hub $i$-diagram}. The subdiagram obtained from $\Delta$ by removing $\pi_1, \pi_2, \xxx_1,\xxx_2$ is denoted by $\Delta'$ and is called the {\em inside subdiagram} of $\Delta$. The band $\xxx_2$ consists of $(k_{i+1},\theta)$-cells, transition cells and gluing $k_{i+1}$-cells. The number of transition cells in $\kkk$ will be called the {\em complexity} of $\Delta$.
\end{notation}

Note that the diagram $\Psi_i(r)$ from  Notation \ref{notphi} is a 2-hub $i$-diagram of complexity 2.

\begin{notation}\label{n:kkk} Let us define several groups given by subpresentations of $\pp'$.
\begin{itemize}

\item[{\bf 1.}] For every $j=2,\ldots,N$, let $K_1^{(j)}$ be the group given by all the relations of $H^{(2j-1)}(A)$ and the conjugacy $q_2^{(2j-1)}$-relations.
\item[{\bf 2.}] For every $j=1,\ldots,N$ let $K_2^{(j)}$ be the group given by all the relations of $\overleftarrow{H}^{(2i)}(A)$ and the conjugacy $q_2^{(2j)}$-relations.
\item[{\bf 3.}] For every $j=2,\ldots, N$ let $K^{(j)}$ be the group given by all the relations of $K_i^{(j)}$, $i=1,2$, and the gluing $k_j$-relations.
\item[{\bf 4.}] Let $K_1^{(1)}$ be the group given by all the relations of $H^{(1)}(A)$, defining relations of $\e\Gamma$ and the conjugacy $q_2^{(1)}$-relations.
\item[{\bf 5.}] Let $K^{(1)}$ be the group given by all the relations of $K_1^{(1)}$, relations of $K_2^{(1)}$, and the gluing $k_1$-relation.
\item[{\bf 6.}] Let $\kkk$ be the free product of all $K^{(i)}$, $i=1,\ldots,N$.
\end{itemize}
Clearly, the presentation of $\kkk$ is a subpresentation of $\pp'$.
\end{notation}

\begin{lemma}\label{l61} For every $j=2,\ldots,N$, the set $X^{(2j-1)}\cup Y^{(2j-1)}$  freely generates a free subgroup in the group $K_1^{(j)}$.
\end{lemma}

\proof Indeed, that group is an HNN-extension of the free product of $H^{(2j-1)}*\la \e(A)^{(2j-1)}\ra$ (the second factor is a free group generated by $\e(A)^{(2j-1)}$)
with free letter $q_2^{(2j-1)}$ and associated subgroups generated by sets $\e(A)^{(2j-1)}$ and $\mu(A)^{(2j-1)}$ respectively. Both sets freely generate free subgroups in the free product by the definitions of $\e$ and $\mu$.
\endproof

The following lemma is proved in the same way as Lemma \ref{l61}

\begin{lemma}\label{l62} For every $j=1,\ldots,N$, the set $X^{(2j)}\cup Y^{(2j)}$ freely generates a free subgroup in the group $K_2^{(j)}$.
\end{lemma}

\begin{lemma}\label{l63} For every $j=2,\ldots,N$, $K^{(j)}$ is constructed as follows.
First take the HNN-extension $\hat K^{(j)}$ of the free product of $K_1^{(j)}*K_2^{(1)}$ with free letter $\hat k_j$ and free associated subgroups generated by $X^{(2j-1)}\cup Y^{(2j-1)}$ and
$X^{(2j)}\cup Y^{(2j)}$ respectively. Then $K^{(j)}$ is an amalgamated product of $\hat K^{(j)}$ and the free group freely generated by $q_3^{(2j-1)},\ldots,q_m^{(2j-1)},k_m, q_3^{(2j)},\ldots,$ $q_m^{(2j)}$ with cyclic associated subgroups generated by $\hat k_j$ and $q_3^{(2j-1)}\ldots q_m^{(2j-1)}k_jq_m^{(2j)}\ldots q_3^{(2j)}$ respectively.
\end{lemma}
\proof This immediately follows from the definition of $K^{(j)}$ and Lemmas \ref{l61}, \ref{l62}. \endproof

\begin{lemma}\label{l65} The group $K^{(1)}$ is constructed as follows. First take the HNN-extension $\hat K^{(1)}$ of the free product of $H'(\Gamma)*K_2^{(1)}$ with free letter $\hat k_1$ and free associated subgroups freely generated by $X^{(1)}\cup Y^{(1)}$ and
$X^{(2)}\cup Y^{(2)}$ respectively. Then $K^{(1)}$ is an amalgamated product of $\hat K^{(1)}$ and the free group freely generated by $q_3^{(1)},\ldots, q_m^{(1)},k_1,q_3^{(2)},\ldots,q_m^{(2)}$ with cyclic associated subgroups generated by $\hat k_1$ and $q_3^{(1)}\ldots q_m^{(1)}k_1q_m^{(2)}\ldots q_3^{(2)}$ respectively.
\end{lemma}
\proof Indeed, the group $K_1^{(1)}$ is clearly isomorphic to $H'(\Gamma)$. Then by Lemma \ref{l64.1}, the group $H(A)/\N$ naturally embeds into $K_1^{(1)}$ where $\N$ is the normal subgroup of $H(A)$ generated by $\mu(R)$ (see Lemma \ref{l64.1}). This, part (1) of Lemma \ref{p1}, and Lemma \ref{l64} imply that the sets $X^{(1)}\cup Y^{(1)}$ and
$X^{(2)}\cup Y^{(2)}$ freely generate free subgroups of $K_1^{(1)}$ and $K_2^{(2)}$ respectively. Then it is clear that the presentation of $\hat K^{(1)}$ is the standard presentation of the HNN-extension and the presentation of $K^{(1)}$ is the standard presentation of an amalgamated product.
\endproof

\begin{lemma}\label{l:69} The set $\bigcup (A^{(j)}\cup B^{(j)})$ generates a free subgroup in $\kkk$.
\end{lemma}

\proof This immediately follows from  Lemma \ref{l:int},  \ref{p1} and \ref{l64.1}.\endproof


\begin{lemma}\label{l:theta} Let $\Delta$ be a reduced diagram over $\pp'$ without hubs. Then
 it does not have $\Theta$-annuli.
\end{lemma}

\proof Suppose that $\Delta$ is a counterexample with the smallest number of cells and $\ttt$ is a $\Theta$-annulus in $\Delta$. Let $\Delta'$ be the inside subdiagram of $\ttt$.

Since $\partial\Delta'$ does not contain $\Theta$-edges, $\Delta'$ does not have any cells corresponding to the defining relations of $\bar\sss_i$, $i=1,2$ (because $\Delta'$ has fewer cells than $\Delta$). Therefore $\Delta'$ is a diagram over the presentation of the group $\kkk$ (since $\Delta$ does not contain hubs, and all other relations have been ruled out). Note that different cells of $\ttt$ cannot have non-$\Theta$-edges in common, otherwise they cancel. Hence every cell in $\ttt$ has a common edge with a cell corresponding to a relation of $\kkk$. Therefore $\ttt$ corresponds to a transition rule of $\sss_1$ or of $\sss_2$.

In the first case
the boundary label of $\Delta'$ is a word in $$\bigcup_{i=1}^{2N}(A\cup B)^i\cup \bigcup_{i=1}^{2N}\{q_1^{(i)},\ldots,q_m^{(i)}\}\cup \{k_1,\ldots,k_N\}.$$ In the second case it is a word in $$\bigcup_{i=2}^{2N}(\e(A)^{(i)} \cup \bigcup_{i=1}^{2N}\{q_1^{(i)},\ldots,q_m^{(i)}\}\cup \{k_1,\ldots,k_N\}.$$ We {\bf claim} that both sets freely generate free subgroups in $\kkk$. Indeed, since $\kkk$ is a free product, it is enough to show that each of the sets

\begin{equation}\label{eq1} A^{(2i-1)}\cup B^{(2i-1)}\cup\{q_1^{(2i-1)},\ldots, q_m^{(2i-1)}\}\cup \{k_i\}\cup A^{(2i)}\cup B^{(2i)}\cup\{q_1^{(2i)},\ldots,q_m^{(2i)}\}\end{equation} the set
\begin{equation}\label{eq2} \e(A)^{(2i-1)}\cup \{q_1^{(2i-1)},\ldots,q_m^{(2i-1)}\}\cup \{k_i\}\cup \{q_1^{(2i)},\ldots,q_m^{(2i)}\}\cup\e(A)^{(2i)}, i\ge 2,\end{equation} and the set
\begin{equation}\label{eq3} \{q_1^{(1)},\ldots,q_m^{(1)}\}\cup \{k_1\}\cup \{q_1^{(2)},\ldots,q_m^{(2)}\}\cup\e(A)^{(2)}\end{equation}
freely generates a free subgroup in the corresponding $K^{(i)}$, $i=1,\ldots,N$. For the set (\ref{eq1}) it follows the representation of $K^{(i)}$ as an amalgamated product (Lemmas \ref{l63}, \ref{l65}) and Lemma \ref{l:int} if $i\ge 2$ and Lemmas \ref{l:int} and \ref{p1} if $i=1$. For the sets (\ref{eq2}), (\ref{eq3}) it follows from the representation of $K^{(i)}$ as an amalgamated product, the fact that the subgroup $E$ of $H(A)$ is freely generated by $\mu(a), a\in A$, (by construction) and the fact that the set $\e(A)^{(j)}$ is conjugated to the set $\mu(A)^{(j)}$ in $K_i^{(j)}$ where $i=1, j\ge 2$ or $i=2, j\ge 1$ by $q_2^{(j)}$.

The claim shows that the inner boundary of $\ttt$ has freely trivial label, hence $\ttt$ has two cells that cancel, a contradiction.
\endproof

\begin{lemma}\label{l:k1}  Let $\Delta$ be a reduced diagram over $\pp'$ without hubs and $\theta$-edges. Then
 $\Delta$ does not have $k_1$-annuli.
\end{lemma}

\proof Suppose that $\Delta$ contains a $k_1$-annulus $\ttt$. Let $\Delta'$ be the inside diagram of that annulus. Note that $\partial\Delta'$ consists of $X^{(1)}\cup Y^{(1)}$-edges. Since $\Delta$ does not contain hubs and $\Theta$-edges, it is a diagram over the presentation of $\kkk$. Using the fact that $\kkk$ is a free product, we conclude that it is a diagram over $K^{(1)}$. But by Lemma \ref{l65}, $X^{(1)}\cup Y^{(1)}$  freely generates a free subgroup in $K^{(1)}$. Hence the label of $\partial\Delta$ is freely trivial, and so $\ttt$ is not reduced, a contradiction.\endproof

The following Lemma is proved in \cite{SBR}.

\begin{lemma}\label{l:hub} Suppose that $\Delta$ is a reduced diagram over $\pp'$ without $t_j$-edges on the boundary, $j=1,\ldots,N$. Suppose that $\Delta$ contains a hub. Then there exists a hub $\pi_1$ in $\Delta$, and an $i<N$, such that the $t_i$-band $\ttt_1$ and the $t_{i+1}$-band $\ttt_2$ starting on $\partial\pi_1$ end on the boundary of another  hub $\pi_2$ of $\Delta$ and there are no hubs in the subdiagram bounded by the medians of $\ttt_1, \ttt_2$, and the parts of $\partial\pi_1$, $\partial\pi_2$ connecting the start edges of $\ttt_1$, $\ttt_2$ and not containing $t$-edges.
\end{lemma}

\subsection{Sectors} Here we shall generalize Lemma \ref{p9.1}.

\begin{lemma}\label{l:sector} Suppose that a reduced diagram $\Delta$ over $\pp'$  without hubs has boundary $\p\q\rf\iv \s\iv$ where $\p$ and  $\rf$ are sides of a $t_i$- and a $k_{i+1}$-bands $\ttt_1$, $\ttt_2$ respectively, and $\Lab(\p)$ and $\Lab(\rf)$ do not contain transition rules $\bar\tau_1^{\pm 1}$, $\Lab(\q)\equiv t_iW_1\bar k_{i+1}$ and $\Lab(\s)\equiv t_iW_2\bar k_{i+1}$ where $W_1, W_2$ are admissible words of $\sss_1^{(2i+1)}$ for some $i=1,\ldots, N-1$. Then the subdiagram of $\Delta$ bounded by the median of $\ttt_1$, the median of $\ttt_2$, $\q, \s$ is combinatorially $\pp$-homotopic to a computational sector of $\sss_1^{(2i+1)}$.
\end{lemma}

\proof
By Lemma \ref{l:theta} every $\Theta$-band in $\Delta$ starts on $\p$ or on $\rf$. Suppose that  a $\Theta$-band starts and ends on $\p$. Since $\Theta$-bands do not intersect, there exists a $\Theta$-band $\aaa$ whose start and end edges belong to the neighbor cells of $\ttt_1$. Then these cells cancel, a contradiction.

Similarly a $\Theta$-band cannot start and end on $\rf$. Therefore $\Theta$-bands that start on $\p$ (on $\rf$) end on $\rf$ (on $\p$). Let $\aaa_1$, $\aaa_2$ be two consecutive $\Theta$-bands in $\Delta$. These bands start on two consecutive $\Theta$-edges of $\p$ and end on two consecutive $\Theta$-edges of $\rf$. The boundary of the subdiagram $\Delta'$ bounded by the medians of $\aaa_1,\aaa_2$ and $\p, \rf$ consists of edges with labels in $Z\cup Q$. The subdiagram $\Delta'$ does not contain hubs and $\Theta$-edges, so it is a diagram over the presentation of $\kkk$. Therefore the boundary label of $\Delta'$ is freely trivial, and $\Delta'$ does not contain cells.

Hence all cells in $\Delta$ correspond to relations of $\bar\sss_1$. It remains to use Lemma \ref{p9.1}.
\endproof

The following lemma is proved in the same way as Lemma \ref{l:sector}.

\begin{lemma}\label{l:sector1} Suppose that a reduced diagram $\Delta$ over $\pp'$  without hubs has boundary $\p\q\rf\iv \s\iv$ where $\p$ and  $\rf$ are sides of a $t_i$- and a $k_{i+1}$-bands $\ttt_1$, $\ttt_2$ respectively, and $\Lab(\p)$ and $\Lab(\rf)$ do not contain transition rules $\bar\tau_2^{\pm 1}$, $\Lab(\q)\equiv t_iW_1\hat k_{i+1}$ and $\Lab(\mathfrak{s})\equiv t_iW_2\hat k_{i+1}$ where $W_1, W_2$ are admissible words of $\sss_2^{(2i+1)}$ for some $i=1,\ldots, N-1$. Then the subdiagram of $\Delta$ bounded by the median of $\ttt_1$, the median of $\ttt_2$, $\q, \mathfrak{s}$ is combinatorially $\pp$-homotopic to a computational sector of $\sss_2^{(2i+1)}$.
\end{lemma}

\subsection{Diagrams with boundary label over $\e(A)$}

\begin{lemma}\label{l:trans} Let $\Delta$ be a diagram over $\pp'$ without hubs. Then it does not contain $(\Theta, k_{i})$-annuli, $N>i\ge 1$.
\end{lemma}

\proof Suppose that there exists a $(\Theta,k)$-annulus $\ttt$ in $\Delta$ composed of a $\Theta$-band $\aaa$ and a $k$-band $\bbb$. Since $\Theta$-bands do not intersect, we can assume that the inside diagram $\Delta'$ of $\ttt$ does not have $\Theta$-edges (by Lemma \ref{l:theta}). Therefore $\Delta'$ is a diagram over the presentation of the group $\kkk$.

Suppose that  $\Delta'$ contains $k_j$-edges. Then it contains a $k_j$-annulus for some $j$. By Lemma \ref{l:k1}, $j>1$. Consider an innermost such annulus. The reduced boundary label of its inside diagram is a word in $X^{(s)}\cup Y^{(s)}$ for some $s>1$. But $X^{(s)}\cup Y^{(s)}$ freely generate a free subgroup of $\kkk$ by Lemma \ref{l:int}, so that annulus contains two cells that cancel. Hence $\Delta'$ does not contain $k_j$-edges.

Therefore $\Delta'$ is a diagram over the presentation of one of the $H(A)^{(s)}$, $s>1$. The label of the boundary of $\Delta'$ is then a product $UV$ where $U$ is a word in $A^{(i)}\cup B^{(i)}$ and $V$ is a word in $X^{(i)}\cup Y^{(i)}$. By Lemma \ref{l:int} then $U=V=1$ in the free group. Hence either $\aaa$ or $\bbb$ contains two cells that cancel, a contradiction.
\endproof

\begin{lemma}\label{l:e} Let $\Delta$ be a reduced diagram over the presentation $\pp'$ with $\partial\Delta$ consisting of $\e(A)$-edges (that includes the case when $\Delta$ is spherical). Suppose that $\Delta$ does not contain hubs. Then $\Delta$ is combinatorially $\pp$-homotopic to a diagram consisting of $\Gamma$-cells.
\end{lemma}

\proof Suppose that $\Delta$ contains $\Theta$-edges. Then by Lemma \ref{l:theta},  $\partial\Delta$ contains $\Theta$-edges, a contradiction. Hence all cells in $\Delta$ correspond to the relations of $\kkk$. Since $\kkk$ is a free product of $K^{(i)}$, $i=1,\ldots,N$, and the boundary label is in $K^{(1)}$, all cells in $\Delta$ correspond to relations of $K^{(1)}$ (consider a maximal subdiagram consisting of cells corresponding to relations of other free factors; all the boundary components of that subdiagram must be empty). By Lemma \ref{l:k1}, $\Delta$ does not contain $k_1$-edges. Hence $\Delta$ is a diagram over the presentation of the free product of $K_1^{(1)}$ and $K_2^{(1)}$. Since the boundary label in $\Delta$ is from $K_1^{(1)}$, all cells in $\Delta$ correspond to relations of $K_1^{(1)}$ which is a copy of $H'(\Gamma)$. It remains to use the first statement of Lemma \ref{l64}.
\endproof

A not necessary reduced diagram over $\pp'$ is called {\em normal} if every hub is contained in a standard subdiagram (i.e. subdiagram of the form $\Psi(r)$, $r\in R$, see Notation \ref{not90}), standard subdiagrams  do not share cells and every $2$-hub $i$-subdiagram is reduced, $i<N$.

\begin{lemma} \label{l:e4} Let $1\le i<N$, then every reduced 2-hub $i$-diagram $\Delta$ (see Notation \ref{not90}) over $\pp'$ is combinatorially $\pp$-homotopic to a normal diagram.
\end{lemma}

\proof Suppose first that the complexity of $\Delta$ is 0. Then both hubs in $\Delta$ are first hubs or both of them are second hubs. Assume both are first hubs (the other case is similar). By Lemma \ref{l:sector}, then the subdiagram of $\Delta$ bounded by the medians of $\xxx_1, \xxx_2$ and parts of the boundaries of the hubs (and not containing the hubs) is combinatorially homotopic to a computational sector of $\sss_1^{(2i+1)}$. By Proposition \ref{mach1}, (4), that computation is empty, so $\ttt_1$ and $\ttt_2$ are empty, and the two hubs in $\Delta$ cancel, a contradiction. The case when both hubs are second hubs is similar only one needs to use Lemma \ref{l:sector1} instead of \ref{l:sector}.

\begin{figure}[ht]
\unitlength .55mm 
\linethickness{0.4pt}
\ifx\plotpoint\undefined\newsavebox{\plotpoint}\fi 
\begin{picture}(256.911,111.25)(0,0)
\put(40.911,31.25){\line(0,1){.5147}}
\put(40.897,31.765){\line(0,1){.5132}}
\put(40.854,32.278){\line(0,1){.5101}}
\put(40.784,32.788){\line(0,1){.5055}}
\multiput(40.686,33.293)(-.03138,.12484){4}{\line(0,1){.12484}}
\multiput(40.561,33.793)(-.030528,.098352){5}{\line(0,1){.098352}}
\multiput(40.408,34.285)(-.029885,.080446){6}{\line(0,1){.080446}}
\multiput(40.229,34.767)(-.029349,.067449){7}{\line(0,1){.067449}}
\multiput(40.023,35.239)(-.032994,.065743){7}{\line(0,1){.065743}}
\multiput(39.792,35.7)(-.031973,.05586){8}{\line(0,1){.05586}}
\multiput(39.537,36.147)(-.031093,.048025){9}{\line(0,1){.048025}}
\multiput(39.257,36.579)(-.033673,.046252){9}{\line(0,1){.046252}}
\multiput(38.954,36.995)(-.032537,.039907){10}{\line(0,1){.039907}}
\multiput(38.628,37.394)(-.031519,.034607){11}{\line(0,1){.034607}}
\multiput(38.282,37.775)(-.033365,.032831){11}{\line(-1,0){.033365}}
\multiput(37.915,38.136)(-.035111,.030957){11}{\line(-1,0){.035111}}
\multiput(37.528,38.476)(-.040427,.03189){10}{\line(-1,0){.040427}}
\multiput(37.124,38.795)(-.046789,.032923){9}{\line(-1,0){.046789}}
\multiput(36.703,39.092)(-.04852,.030315){9}{\line(-1,0){.04852}}
\multiput(36.266,39.364)(-.056368,.031068){8}{\line(-1,0){.056368}}
\multiput(35.815,39.613)(-.066267,.03193){7}{\line(-1,0){.066267}}
\multiput(35.351,39.837)(-.079233,.032967){6}{\line(-1,0){.079233}}
\multiput(34.876,40.034)(-.080917,.028584){6}{\line(-1,0){.080917}}
\multiput(34.391,40.206)(-.098831,.028939){5}{\line(-1,0){.098831}}
\put(33.896,40.351){\line(-1,0){.5013}}
\put(33.395,40.468){\line(-1,0){.507}}
\put(32.888,40.558){\line(-1,0){.5112}}
\put(32.377,40.62){\line(-1,0){.5138}}
\put(31.863,40.654){\line(-1,0){.5149}}
\put(31.348,40.66){\line(-1,0){.5144}}
\put(30.834,40.637){\line(-1,0){.5124}}
\put(30.321,40.587){\line(-1,0){.5089}}
\put(29.813,40.508){\line(-1,0){.5038}}
\multiput(29.309,40.402)(-.12432,-.03339){4}{\line(-1,0){.12432}}
\multiput(28.811,40.269)(-.097847,-.03211){5}{\line(-1,0){.097847}}
\multiput(28.322,40.108)(-.079954,-.031178){6}{\line(-1,0){.079954}}
\multiput(27.842,39.921)(-.066968,-.030432){7}{\line(-1,0){.066968}}
\multiput(27.374,39.708)(-.057053,-.029793){8}{\line(-1,0){.057053}}
\multiput(26.917,39.47)(-.055338,-.032869){8}{\line(-1,0){.055338}}
\multiput(26.475,39.207)(-.047517,-.031863){9}{\line(-1,0){.047517}}
\multiput(26.047,38.92)(-.041133,-.030973){10}{\line(-1,0){.041133}}
\multiput(25.636,38.61)(-.039378,-.033177){10}{\line(-1,0){.039378}}
\multiput(25.242,38.278)(-.034095,-.032073){11}{\line(-1,0){.034095}}
\multiput(24.867,37.926)(-.032289,-.03389){11}{\line(0,-1){.03389}}
\multiput(24.512,37.553)(-.033426,-.039166){10}{\line(0,-1){.039166}}
\multiput(24.177,37.161)(-.031234,-.040935){10}{\line(0,-1){.040935}}
\multiput(23.865,36.752)(-.032165,-.047314){9}{\line(0,-1){.047314}}
\multiput(23.575,36.326)(-.03322,-.055128){8}{\line(0,-1){.055128}}
\multiput(23.31,35.885)(-.030155,-.056862){8}{\line(0,-1){.056862}}
\multiput(23.068,35.43)(-.030857,-.066773){7}{\line(0,-1){.066773}}
\multiput(22.852,34.963)(-.031685,-.079754){6}{\line(0,-1){.079754}}
\multiput(22.662,34.484)(-.032731,-.097641){5}{\line(0,-1){.097641}}
\multiput(22.499,33.996)(-.027342,-.099285){5}{\line(0,-1){.099285}}
\put(22.362,33.499){\line(0,-1){.5032}}
\put(22.253,32.996){\line(0,-1){.5084}}
\put(22.171,32.488){\line(0,-1){.5121}}
\put(22.117,31.976){\line(0,-1){2.0547}}
\put(22.184,29.921){\line(0,-1){.5076}}
\put(22.27,29.414){\line(0,-1){.5021}}
\multiput(22.384,28.911)(.02831,-.099013){5}{\line(0,-1){.099013}}
\multiput(22.526,28.416)(.033683,-.097317){5}{\line(0,-1){.097317}}
\multiput(22.694,27.93)(.032463,-.079441){6}{\line(0,-1){.079441}}
\multiput(22.889,27.453)(.031508,-.066468){7}{\line(0,-1){.066468}}
\multiput(23.11,26.988)(.030709,-.056565){8}{\line(0,-1){.056565}}
\multiput(23.355,26.535)(.030006,-.048712){9}{\line(0,-1){.048712}}
\multiput(23.625,26.097)(.032625,-.046998){9}{\line(0,-1){.046998}}
\multiput(23.919,25.674)(.031632,-.040628){10}{\line(0,-1){.040628}}
\multiput(24.235,25.268)(.030734,-.035307){11}{\line(0,-1){.035307}}
\multiput(24.573,24.879)(.032619,-.033573){11}{\line(0,-1){.033573}}
\multiput(24.932,24.51)(.034406,-.031739){11}{\line(1,0){.034406}}
\multiput(25.311,24.161)(.0397,-.03279){10}{\line(1,0){.0397}}
\multiput(25.708,23.833)(.041434,-.03057){10}{\line(1,0){.041434}}
\multiput(26.122,23.527)(.047826,-.031398){9}{\line(1,0){.047826}}
\multiput(26.553,23.245)(.055656,-.032327){8}{\line(1,0){.055656}}
\multiput(26.998,22.986)(.065532,-.033411){7}{\line(1,0){.065532}}
\multiput(27.456,22.752)(.067262,-.029777){7}{\line(1,0){.067262}}
\multiput(27.927,22.544)(.080254,-.030395){6}{\line(1,0){.080254}}
\multiput(28.409,22.361)(.098156,-.031153){5}{\line(1,0){.098156}}
\multiput(28.9,22.206)(.12464,-.03217){4}{\line(1,0){.12464}}
\put(29.398,22.077){\line(1,0){.5049}}
\put(29.903,21.976){\line(1,0){.5096}}
\put(30.413,21.902){\line(1,0){.5129}}
\put(30.926,21.857){\line(1,0){1.0294}}
\put(31.955,21.85){\line(1,0){.5134}}
\put(32.468,21.889){\line(1,0){.5105}}
\put(32.979,21.956){\line(1,0){.5061}}
\multiput(33.485,22.051)(.12504,.03058){4}{\line(1,0){.12504}}
\multiput(33.985,22.173)(.098544,.029902){5}{\line(1,0){.098544}}
\multiput(34.478,22.323)(.080634,.029373){6}{\line(1,0){.080634}}
\multiput(34.962,22.499)(.078907,.033739){6}{\line(1,0){.078907}}
\multiput(35.435,22.701)(.065952,.032575){7}{\line(1,0){.065952}}
\multiput(35.897,22.93)(.056062,.031617){8}{\line(1,0){.056062}}
\multiput(36.345,23.182)(.048222,.030787){9}{\line(1,0){.048222}}
\multiput(36.779,23.46)(.046465,.033379){9}{\line(1,0){.046465}}
\multiput(37.198,23.76)(.040113,.032283){10}{\line(1,0){.040113}}
\multiput(37.599,24.083)(.034807,.031299){11}{\line(1,0){.034807}}
\multiput(37.982,24.427)(.033043,.033156){11}{\line(0,1){.033156}}
\multiput(38.345,24.792)(.03118,.034913){11}{\line(0,1){.034913}}
\multiput(38.688,25.176)(.032146,.040223){10}{\line(0,1){.040223}}
\multiput(39.009,25.578)(.03322,.046579){9}{\line(0,1){.046579}}
\multiput(39.308,25.997)(.030623,.048326){9}{\line(0,1){.048326}}
\multiput(39.584,26.432)(.031425,.05617){8}{\line(0,1){.05617}}
\multiput(39.835,26.882)(.03235,.066062){7}{\line(0,1){.066062}}
\multiput(40.062,27.344)(.03347,.079022){6}{\line(0,1){.079022}}
\multiput(40.263,27.818)(.029098,.080734){6}{\line(0,1){.080734}}
\multiput(40.437,28.303)(.029566,.098646){5}{\line(0,1){.098646}}
\multiput(40.585,28.796)(.03016,.12514){4}{\line(0,1){.12514}}
\put(40.706,29.296){\line(0,1){.5064}}
\put(40.799,29.803){\line(0,1){.5108}}
\put(40.864,30.314){\line(0,1){.9365}}
\put(115.661,31.5){\line(0,1){.5147}}
\put(115.647,32.015){\line(0,1){.5132}}
\put(115.604,32.528){\line(0,1){.5101}}
\put(115.534,33.038){\line(0,1){.5055}}
\multiput(115.436,33.543)(-.03138,.12484){4}{\line(0,1){.12484}}
\multiput(115.311,34.043)(-.030528,.098352){5}{\line(0,1){.098352}}
\multiput(115.158,34.535)(-.029885,.080446){6}{\line(0,1){.080446}}
\multiput(114.979,35.017)(-.029349,.067449){7}{\line(0,1){.067449}}
\multiput(114.773,35.489)(-.032994,.065743){7}{\line(0,1){.065743}}
\multiput(114.542,35.95)(-.031973,.05586){8}{\line(0,1){.05586}}
\multiput(114.287,36.397)(-.031093,.048025){9}{\line(0,1){.048025}}
\multiput(114.007,36.829)(-.033673,.046252){9}{\line(0,1){.046252}}
\multiput(113.704,37.245)(-.032537,.039907){10}{\line(0,1){.039907}}
\multiput(113.378,37.644)(-.031519,.034607){11}{\line(0,1){.034607}}
\multiput(113.032,38.025)(-.033365,.032831){11}{\line(-1,0){.033365}}
\multiput(112.665,38.386)(-.035111,.030957){11}{\line(-1,0){.035111}}
\multiput(112.278,38.726)(-.040427,.03189){10}{\line(-1,0){.040427}}
\multiput(111.874,39.045)(-.046789,.032923){9}{\line(-1,0){.046789}}
\multiput(111.453,39.342)(-.04852,.030315){9}{\line(-1,0){.04852}}
\multiput(111.016,39.614)(-.056368,.031068){8}{\line(-1,0){.056368}}
\multiput(110.565,39.863)(-.066267,.03193){7}{\line(-1,0){.066267}}
\multiput(110.101,40.087)(-.079233,.032967){6}{\line(-1,0){.079233}}
\multiput(109.626,40.284)(-.080917,.028584){6}{\line(-1,0){.080917}}
\multiput(109.141,40.456)(-.098831,.028939){5}{\line(-1,0){.098831}}
\put(108.646,40.601){\line(-1,0){.5013}}
\put(108.145,40.718){\line(-1,0){.507}}
\put(107.638,40.808){\line(-1,0){.5112}}
\put(107.127,40.87){\line(-1,0){.5138}}
\put(106.613,40.904){\line(-1,0){.5149}}
\put(106.098,40.91){\line(-1,0){.5144}}
\put(105.584,40.887){\line(-1,0){.5124}}
\put(105.071,40.837){\line(-1,0){.5089}}
\put(104.563,40.758){\line(-1,0){.5038}}
\multiput(104.059,40.652)(-.12432,-.03339){4}{\line(-1,0){.12432}}
\multiput(103.561,40.519)(-.097847,-.03211){5}{\line(-1,0){.097847}}
\multiput(103.072,40.358)(-.079954,-.031178){6}{\line(-1,0){.079954}}
\multiput(102.592,40.171)(-.066968,-.030432){7}{\line(-1,0){.066968}}
\multiput(102.124,39.958)(-.057053,-.029793){8}{\line(-1,0){.057053}}
\multiput(101.667,39.72)(-.055338,-.032869){8}{\line(-1,0){.055338}}
\multiput(101.225,39.457)(-.047517,-.031863){9}{\line(-1,0){.047517}}
\multiput(100.797,39.17)(-.041133,-.030973){10}{\line(-1,0){.041133}}
\multiput(100.386,38.86)(-.039378,-.033177){10}{\line(-1,0){.039378}}
\multiput(99.992,38.528)(-.034095,-.032073){11}{\line(-1,0){.034095}}
\multiput(99.617,38.176)(-.032289,-.03389){11}{\line(0,-1){.03389}}
\multiput(99.262,37.803)(-.033426,-.039166){10}{\line(0,-1){.039166}}
\multiput(98.927,37.411)(-.031234,-.040935){10}{\line(0,-1){.040935}}
\multiput(98.615,37.002)(-.032165,-.047314){9}{\line(0,-1){.047314}}
\multiput(98.325,36.576)(-.03322,-.055128){8}{\line(0,-1){.055128}}
\multiput(98.06,36.135)(-.030155,-.056862){8}{\line(0,-1){.056862}}
\multiput(97.818,35.68)(-.030857,-.066773){7}{\line(0,-1){.066773}}
\multiput(97.602,35.213)(-.031685,-.079754){6}{\line(0,-1){.079754}}
\multiput(97.412,34.734)(-.032731,-.097641){5}{\line(0,-1){.097641}}
\multiput(97.249,34.246)(-.027342,-.099285){5}{\line(0,-1){.099285}}
\put(97.112,33.749){\line(0,-1){.5032}}
\put(97.003,33.246){\line(0,-1){.5084}}
\put(96.921,32.738){\line(0,-1){.5121}}
\put(96.867,32.226){\line(0,-1){2.0547}}
\put(96.934,30.171){\line(0,-1){.5076}}
\put(97.02,29.664){\line(0,-1){.5021}}
\multiput(97.134,29.161)(.02831,-.099013){5}{\line(0,-1){.099013}}
\multiput(97.276,28.666)(.033683,-.097317){5}{\line(0,-1){.097317}}
\multiput(97.444,28.18)(.032463,-.079441){6}{\line(0,-1){.079441}}
\multiput(97.639,27.703)(.031508,-.066468){7}{\line(0,-1){.066468}}
\multiput(97.86,27.238)(.030709,-.056565){8}{\line(0,-1){.056565}}
\multiput(98.105,26.785)(.030006,-.048712){9}{\line(0,-1){.048712}}
\multiput(98.375,26.347)(.032625,-.046998){9}{\line(0,-1){.046998}}
\multiput(98.669,25.924)(.031632,-.040628){10}{\line(0,-1){.040628}}
\multiput(98.985,25.518)(.030734,-.035307){11}{\line(0,-1){.035307}}
\multiput(99.323,25.129)(.032619,-.033573){11}{\line(0,-1){.033573}}
\multiput(99.682,24.76)(.034406,-.031739){11}{\line(1,0){.034406}}
\multiput(100.061,24.411)(.0397,-.03279){10}{\line(1,0){.0397}}
\multiput(100.458,24.083)(.041434,-.03057){10}{\line(1,0){.041434}}
\multiput(100.872,23.777)(.047826,-.031398){9}{\line(1,0){.047826}}
\multiput(101.303,23.495)(.055656,-.032327){8}{\line(1,0){.055656}}
\multiput(101.748,23.236)(.065532,-.033411){7}{\line(1,0){.065532}}
\multiput(102.206,23.002)(.067262,-.029777){7}{\line(1,0){.067262}}
\multiput(102.677,22.794)(.080254,-.030395){6}{\line(1,0){.080254}}
\multiput(103.159,22.611)(.098156,-.031153){5}{\line(1,0){.098156}}
\multiput(103.65,22.456)(.12464,-.03217){4}{\line(1,0){.12464}}
\put(104.148,22.327){\line(1,0){.5049}}
\put(104.653,22.226){\line(1,0){.5096}}
\put(105.163,22.152){\line(1,0){.5129}}
\put(105.676,22.107){\line(1,0){1.0294}}
\put(106.705,22.1){\line(1,0){.5134}}
\put(107.218,22.139){\line(1,0){.5105}}
\put(107.729,22.206){\line(1,0){.5061}}
\multiput(108.235,22.301)(.12504,.03058){4}{\line(1,0){.12504}}
\multiput(108.735,22.423)(.098544,.029902){5}{\line(1,0){.098544}}
\multiput(109.228,22.573)(.080634,.029373){6}{\line(1,0){.080634}}
\multiput(109.712,22.749)(.078907,.033739){6}{\line(1,0){.078907}}
\multiput(110.185,22.951)(.065952,.032575){7}{\line(1,0){.065952}}
\multiput(110.647,23.18)(.056062,.031617){8}{\line(1,0){.056062}}
\multiput(111.095,23.432)(.048222,.030787){9}{\line(1,0){.048222}}
\multiput(111.529,23.71)(.046465,.033379){9}{\line(1,0){.046465}}
\multiput(111.948,24.01)(.040113,.032283){10}{\line(1,0){.040113}}
\multiput(112.349,24.333)(.034807,.031299){11}{\line(1,0){.034807}}
\multiput(112.732,24.677)(.033043,.033156){11}{\line(0,1){.033156}}
\multiput(113.095,25.042)(.03118,.034913){11}{\line(0,1){.034913}}
\multiput(113.438,25.426)(.032146,.040223){10}{\line(0,1){.040223}}
\multiput(113.759,25.828)(.03322,.046579){9}{\line(0,1){.046579}}
\multiput(114.058,26.247)(.030623,.048326){9}{\line(0,1){.048326}}
\multiput(114.334,26.682)(.031425,.05617){8}{\line(0,1){.05617}}
\multiput(114.585,27.132)(.03235,.066062){7}{\line(0,1){.066062}}
\multiput(114.812,27.594)(.03347,.079022){6}{\line(0,1){.079022}}
\multiput(115.013,28.068)(.029098,.080734){6}{\line(0,1){.080734}}
\multiput(115.187,28.553)(.029566,.098646){5}{\line(0,1){.098646}}
\multiput(115.335,29.046)(.03016,.12514){4}{\line(0,1){.12514}}
\put(115.456,29.546){\line(0,1){.5064}}
\put(115.549,30.053){\line(0,1){.5108}}
\put(115.614,30.564){\line(0,1){.9365}}
\put(151.911,31.5){\line(0,1){.5147}}
\put(151.897,32.015){\line(0,1){.5132}}
\put(151.854,32.528){\line(0,1){.5101}}
\put(151.784,33.038){\line(0,1){.5055}}
\multiput(151.686,33.543)(-.03138,.12484){4}{\line(0,1){.12484}}
\multiput(151.561,34.043)(-.030528,.098352){5}{\line(0,1){.098352}}
\multiput(151.408,34.535)(-.029885,.080446){6}{\line(0,1){.080446}}
\multiput(151.229,35.017)(-.029349,.067449){7}{\line(0,1){.067449}}
\multiput(151.023,35.489)(-.032994,.065743){7}{\line(0,1){.065743}}
\multiput(150.792,35.95)(-.031973,.05586){8}{\line(0,1){.05586}}
\multiput(150.537,36.397)(-.031093,.048025){9}{\line(0,1){.048025}}
\multiput(150.257,36.829)(-.033673,.046252){9}{\line(0,1){.046252}}
\multiput(149.954,37.245)(-.032537,.039907){10}{\line(0,1){.039907}}
\multiput(149.628,37.644)(-.031519,.034607){11}{\line(0,1){.034607}}
\multiput(149.282,38.025)(-.033365,.032831){11}{\line(-1,0){.033365}}
\multiput(148.915,38.386)(-.035111,.030957){11}{\line(-1,0){.035111}}
\multiput(148.528,38.726)(-.040427,.03189){10}{\line(-1,0){.040427}}
\multiput(148.124,39.045)(-.046789,.032923){9}{\line(-1,0){.046789}}
\multiput(147.703,39.342)(-.04852,.030315){9}{\line(-1,0){.04852}}
\multiput(147.266,39.614)(-.056368,.031068){8}{\line(-1,0){.056368}}
\multiput(146.815,39.863)(-.066267,.03193){7}{\line(-1,0){.066267}}
\multiput(146.351,40.087)(-.079233,.032967){6}{\line(-1,0){.079233}}
\multiput(145.876,40.284)(-.080917,.028584){6}{\line(-1,0){.080917}}
\multiput(145.391,40.456)(-.098831,.028939){5}{\line(-1,0){.098831}}
\put(144.896,40.601){\line(-1,0){.5013}}
\put(144.395,40.718){\line(-1,0){.507}}
\put(143.888,40.808){\line(-1,0){.5112}}
\put(143.377,40.87){\line(-1,0){.5138}}
\put(142.863,40.904){\line(-1,0){.5149}}
\put(142.348,40.91){\line(-1,0){.5144}}
\put(141.834,40.887){\line(-1,0){.5124}}
\put(141.321,40.837){\line(-1,0){.5089}}
\put(140.813,40.758){\line(-1,0){.5038}}
\multiput(140.309,40.652)(-.12432,-.03339){4}{\line(-1,0){.12432}}
\multiput(139.811,40.519)(-.097847,-.03211){5}{\line(-1,0){.097847}}
\multiput(139.322,40.358)(-.079954,-.031178){6}{\line(-1,0){.079954}}
\multiput(138.842,40.171)(-.066968,-.030432){7}{\line(-1,0){.066968}}
\multiput(138.374,39.958)(-.057053,-.029793){8}{\line(-1,0){.057053}}
\multiput(137.917,39.72)(-.055338,-.032869){8}{\line(-1,0){.055338}}
\multiput(137.475,39.457)(-.047517,-.031863){9}{\line(-1,0){.047517}}
\multiput(137.047,39.17)(-.041133,-.030973){10}{\line(-1,0){.041133}}
\multiput(136.636,38.86)(-.039378,-.033177){10}{\line(-1,0){.039378}}
\multiput(136.242,38.528)(-.034095,-.032073){11}{\line(-1,0){.034095}}
\multiput(135.867,38.176)(-.032289,-.03389){11}{\line(0,-1){.03389}}
\multiput(135.512,37.803)(-.033426,-.039166){10}{\line(0,-1){.039166}}
\multiput(135.177,37.411)(-.031234,-.040935){10}{\line(0,-1){.040935}}
\multiput(134.865,37.002)(-.032165,-.047314){9}{\line(0,-1){.047314}}
\multiput(134.575,36.576)(-.03322,-.055128){8}{\line(0,-1){.055128}}
\multiput(134.31,36.135)(-.030155,-.056862){8}{\line(0,-1){.056862}}
\multiput(134.068,35.68)(-.030857,-.066773){7}{\line(0,-1){.066773}}
\multiput(133.852,35.213)(-.031685,-.079754){6}{\line(0,-1){.079754}}
\multiput(133.662,34.734)(-.032731,-.097641){5}{\line(0,-1){.097641}}
\multiput(133.499,34.246)(-.027342,-.099285){5}{\line(0,-1){.099285}}
\put(133.362,33.749){\line(0,-1){.5032}}
\put(133.253,33.246){\line(0,-1){.5084}}
\put(133.171,32.738){\line(0,-1){.5121}}
\put(133.117,32.226){\line(0,-1){2.0547}}
\put(133.184,30.171){\line(0,-1){.5076}}
\put(133.27,29.664){\line(0,-1){.5021}}
\multiput(133.384,29.161)(.02831,-.099013){5}{\line(0,-1){.099013}}
\multiput(133.526,28.666)(.033683,-.097317){5}{\line(0,-1){.097317}}
\multiput(133.694,28.18)(.032463,-.079441){6}{\line(0,-1){.079441}}
\multiput(133.889,27.703)(.031508,-.066468){7}{\line(0,-1){.066468}}
\multiput(134.11,27.238)(.030709,-.056565){8}{\line(0,-1){.056565}}
\multiput(134.355,26.785)(.030006,-.048712){9}{\line(0,-1){.048712}}
\multiput(134.625,26.347)(.032625,-.046998){9}{\line(0,-1){.046998}}
\multiput(134.919,25.924)(.031632,-.040628){10}{\line(0,-1){.040628}}
\multiput(135.235,25.518)(.030734,-.035307){11}{\line(0,-1){.035307}}
\multiput(135.573,25.129)(.032619,-.033573){11}{\line(0,-1){.033573}}
\multiput(135.932,24.76)(.034406,-.031739){11}{\line(1,0){.034406}}
\multiput(136.311,24.411)(.0397,-.03279){10}{\line(1,0){.0397}}
\multiput(136.708,24.083)(.041434,-.03057){10}{\line(1,0){.041434}}
\multiput(137.122,23.777)(.047826,-.031398){9}{\line(1,0){.047826}}
\multiput(137.553,23.495)(.055656,-.032327){8}{\line(1,0){.055656}}
\multiput(137.998,23.236)(.065532,-.033411){7}{\line(1,0){.065532}}
\multiput(138.456,23.002)(.067262,-.029777){7}{\line(1,0){.067262}}
\multiput(138.927,22.794)(.080254,-.030395){6}{\line(1,0){.080254}}
\multiput(139.409,22.611)(.098156,-.031153){5}{\line(1,0){.098156}}
\multiput(139.9,22.456)(.12464,-.03217){4}{\line(1,0){.12464}}
\put(140.398,22.327){\line(1,0){.5049}}
\put(140.903,22.226){\line(1,0){.5096}}
\put(141.413,22.152){\line(1,0){.5129}}
\put(141.926,22.107){\line(1,0){1.0294}}
\put(142.955,22.1){\line(1,0){.5134}}
\put(143.468,22.139){\line(1,0){.5105}}
\put(143.979,22.206){\line(1,0){.5061}}
\multiput(144.485,22.301)(.12504,.03058){4}{\line(1,0){.12504}}
\multiput(144.985,22.423)(.098544,.029902){5}{\line(1,0){.098544}}
\multiput(145.478,22.573)(.080634,.029373){6}{\line(1,0){.080634}}
\multiput(145.962,22.749)(.078907,.033739){6}{\line(1,0){.078907}}
\multiput(146.435,22.951)(.065952,.032575){7}{\line(1,0){.065952}}
\multiput(146.897,23.18)(.056062,.031617){8}{\line(1,0){.056062}}
\multiput(147.345,23.432)(.048222,.030787){9}{\line(1,0){.048222}}
\multiput(147.779,23.71)(.046465,.033379){9}{\line(1,0){.046465}}
\multiput(148.198,24.01)(.040113,.032283){10}{\line(1,0){.040113}}
\multiput(148.599,24.333)(.034807,.031299){11}{\line(1,0){.034807}}
\multiput(148.982,24.677)(.033043,.033156){11}{\line(0,1){.033156}}
\multiput(149.345,25.042)(.03118,.034913){11}{\line(0,1){.034913}}
\multiput(149.688,25.426)(.032146,.040223){10}{\line(0,1){.040223}}
\multiput(150.009,25.828)(.03322,.046579){9}{\line(0,1){.046579}}
\multiput(150.308,26.247)(.030623,.048326){9}{\line(0,1){.048326}}
\multiput(150.584,26.682)(.031425,.05617){8}{\line(0,1){.05617}}
\multiput(150.835,27.132)(.03235,.066062){7}{\line(0,1){.066062}}
\multiput(151.062,27.594)(.03347,.079022){6}{\line(0,1){.079022}}
\multiput(151.263,28.068)(.029098,.080734){6}{\line(0,1){.080734}}
\multiput(151.437,28.553)(.029566,.098646){5}{\line(0,1){.098646}}
\multiput(151.585,29.046)(.03016,.12514){4}{\line(0,1){.12514}}
\put(151.706,29.546){\line(0,1){.5064}}
\put(151.799,30.053){\line(0,1){.5108}}
\put(151.864,30.564){\line(0,1){.9365}}
\put(74.661,96.5){\line(0,1){.5147}}
\put(74.647,97.015){\line(0,1){.5132}}
\put(74.604,97.528){\line(0,1){.5101}}
\put(74.534,98.038){\line(0,1){.5055}}
\multiput(74.436,98.543)(-.03138,.12484){4}{\line(0,1){.12484}}
\multiput(74.311,99.043)(-.030528,.098352){5}{\line(0,1){.098352}}
\multiput(74.158,99.535)(-.029885,.080446){6}{\line(0,1){.080446}}
\multiput(73.979,100.017)(-.029349,.067449){7}{\line(0,1){.067449}}
\multiput(73.773,100.489)(-.032994,.065743){7}{\line(0,1){.065743}}
\multiput(73.542,100.95)(-.031973,.05586){8}{\line(0,1){.05586}}
\multiput(73.287,101.397)(-.031093,.048025){9}{\line(0,1){.048025}}
\multiput(73.007,101.829)(-.033673,.046252){9}{\line(0,1){.046252}}
\multiput(72.704,102.245)(-.032537,.039907){10}{\line(0,1){.039907}}
\multiput(72.378,102.644)(-.031519,.034607){11}{\line(0,1){.034607}}
\multiput(72.032,103.025)(-.033365,.032831){11}{\line(-1,0){.033365}}
\multiput(71.665,103.386)(-.035111,.030957){11}{\line(-1,0){.035111}}
\multiput(71.278,103.726)(-.040427,.03189){10}{\line(-1,0){.040427}}
\multiput(70.874,104.045)(-.046789,.032923){9}{\line(-1,0){.046789}}
\multiput(70.453,104.342)(-.04852,.030315){9}{\line(-1,0){.04852}}
\multiput(70.016,104.614)(-.056368,.031068){8}{\line(-1,0){.056368}}
\multiput(69.565,104.863)(-.066267,.03193){7}{\line(-1,0){.066267}}
\multiput(69.101,105.087)(-.079233,.032967){6}{\line(-1,0){.079233}}
\multiput(68.626,105.284)(-.080917,.028584){6}{\line(-1,0){.080917}}
\multiput(68.141,105.456)(-.098831,.028939){5}{\line(-1,0){.098831}}
\put(67.646,105.601){\line(-1,0){.5013}}
\put(67.145,105.718){\line(-1,0){.507}}
\put(66.638,105.808){\line(-1,0){.5112}}
\put(66.127,105.87){\line(-1,0){.5138}}
\put(65.613,105.904){\line(-1,0){.5149}}
\put(65.098,105.91){\line(-1,0){.5144}}
\put(64.584,105.887){\line(-1,0){.5124}}
\put(64.071,105.837){\line(-1,0){.5089}}
\put(63.563,105.758){\line(-1,0){.5038}}
\multiput(63.059,105.652)(-.12432,-.03339){4}{\line(-1,0){.12432}}
\multiput(62.561,105.519)(-.097847,-.03211){5}{\line(-1,0){.097847}}
\multiput(62.072,105.358)(-.079954,-.031178){6}{\line(-1,0){.079954}}
\multiput(61.592,105.171)(-.066968,-.030432){7}{\line(-1,0){.066968}}
\multiput(61.124,104.958)(-.057053,-.029793){8}{\line(-1,0){.057053}}
\multiput(60.667,104.72)(-.055338,-.032869){8}{\line(-1,0){.055338}}
\multiput(60.225,104.457)(-.047517,-.031863){9}{\line(-1,0){.047517}}
\multiput(59.797,104.17)(-.041133,-.030973){10}{\line(-1,0){.041133}}
\multiput(59.386,103.86)(-.039378,-.033177){10}{\line(-1,0){.039378}}
\multiput(58.992,103.528)(-.034095,-.032073){11}{\line(-1,0){.034095}}
\multiput(58.617,103.176)(-.032289,-.03389){11}{\line(0,-1){.03389}}
\multiput(58.262,102.803)(-.033426,-.039166){10}{\line(0,-1){.039166}}
\multiput(57.927,102.411)(-.031234,-.040935){10}{\line(0,-1){.040935}}
\multiput(57.615,102.002)(-.032165,-.047314){9}{\line(0,-1){.047314}}
\multiput(57.325,101.576)(-.03322,-.055128){8}{\line(0,-1){.055128}}
\multiput(57.06,101.135)(-.030155,-.056862){8}{\line(0,-1){.056862}}
\multiput(56.818,100.68)(-.030857,-.066773){7}{\line(0,-1){.066773}}
\multiput(56.602,100.213)(-.031685,-.079754){6}{\line(0,-1){.079754}}
\multiput(56.412,99.734)(-.032731,-.097641){5}{\line(0,-1){.097641}}
\multiput(56.249,99.246)(-.027342,-.099285){5}{\line(0,-1){.099285}}
\put(56.112,98.749){\line(0,-1){.5032}}
\put(56.003,98.246){\line(0,-1){.5084}}
\put(55.921,97.738){\line(0,-1){.5121}}
\put(55.867,97.226){\line(0,-1){2.0547}}
\put(55.934,95.171){\line(0,-1){.5076}}
\put(56.02,94.664){\line(0,-1){.5021}}
\multiput(56.134,94.161)(.02831,-.099013){5}{\line(0,-1){.099013}}
\multiput(56.276,93.666)(.033683,-.097317){5}{\line(0,-1){.097317}}
\multiput(56.444,93.18)(.032463,-.079441){6}{\line(0,-1){.079441}}
\multiput(56.639,92.703)(.031508,-.066468){7}{\line(0,-1){.066468}}
\multiput(56.86,92.238)(.030709,-.056565){8}{\line(0,-1){.056565}}
\multiput(57.105,91.785)(.030006,-.048712){9}{\line(0,-1){.048712}}
\multiput(57.375,91.347)(.032625,-.046998){9}{\line(0,-1){.046998}}
\multiput(57.669,90.924)(.031632,-.040628){10}{\line(0,-1){.040628}}
\multiput(57.985,90.518)(.030734,-.035307){11}{\line(0,-1){.035307}}
\multiput(58.323,90.129)(.032619,-.033573){11}{\line(0,-1){.033573}}
\multiput(58.682,89.76)(.034406,-.031739){11}{\line(1,0){.034406}}
\multiput(59.061,89.411)(.0397,-.03279){10}{\line(1,0){.0397}}
\multiput(59.458,89.083)(.041434,-.03057){10}{\line(1,0){.041434}}
\multiput(59.872,88.777)(.047826,-.031398){9}{\line(1,0){.047826}}
\multiput(60.303,88.495)(.055656,-.032327){8}{\line(1,0){.055656}}
\multiput(60.748,88.236)(.065532,-.033411){7}{\line(1,0){.065532}}
\multiput(61.206,88.002)(.067262,-.029777){7}{\line(1,0){.067262}}
\multiput(61.677,87.794)(.080254,-.030395){6}{\line(1,0){.080254}}
\multiput(62.159,87.611)(.098156,-.031153){5}{\line(1,0){.098156}}
\multiput(62.65,87.456)(.12464,-.03217){4}{\line(1,0){.12464}}
\put(63.148,87.327){\line(1,0){.5049}}
\put(63.653,87.226){\line(1,0){.5096}}
\put(64.163,87.152){\line(1,0){.5129}}
\put(64.676,87.107){\line(1,0){1.0294}}
\put(65.705,87.1){\line(1,0){.5134}}
\put(66.218,87.139){\line(1,0){.5105}}
\put(66.729,87.206){\line(1,0){.5061}}
\multiput(67.235,87.301)(.12504,.03058){4}{\line(1,0){.12504}}
\multiput(67.735,87.423)(.098544,.029902){5}{\line(1,0){.098544}}
\multiput(68.228,87.573)(.080634,.029373){6}{\line(1,0){.080634}}
\multiput(68.712,87.749)(.078907,.033739){6}{\line(1,0){.078907}}
\multiput(69.185,87.951)(.065952,.032575){7}{\line(1,0){.065952}}
\multiput(69.647,88.18)(.056062,.031617){8}{\line(1,0){.056062}}
\multiput(70.095,88.432)(.048222,.030787){9}{\line(1,0){.048222}}
\multiput(70.529,88.71)(.046465,.033379){9}{\line(1,0){.046465}}
\multiput(70.948,89.01)(.040113,.032283){10}{\line(1,0){.040113}}
\multiput(71.349,89.333)(.034807,.031299){11}{\line(1,0){.034807}}
\multiput(71.732,89.677)(.033043,.033156){11}{\line(0,1){.033156}}
\multiput(72.095,90.042)(.03118,.034913){11}{\line(0,1){.034913}}
\multiput(72.438,90.426)(.032146,.040223){10}{\line(0,1){.040223}}
\multiput(72.759,90.828)(.03322,.046579){9}{\line(0,1){.046579}}
\multiput(73.058,91.247)(.030623,.048326){9}{\line(0,1){.048326}}
\multiput(73.334,91.682)(.031425,.05617){8}{\line(0,1){.05617}}
\multiput(73.585,92.132)(.03235,.066062){7}{\line(0,1){.066062}}
\multiput(73.812,92.594)(.03347,.079022){6}{\line(0,1){.079022}}
\multiput(74.013,93.068)(.029098,.080734){6}{\line(0,1){.080734}}
\multiput(74.187,93.553)(.029566,.098646){5}{\line(0,1){.098646}}
\multiput(74.335,94.046)(.03016,.12514){4}{\line(0,1){.12514}}
\put(74.456,94.546){\line(0,1){.5064}}
\put(74.549,95.053){\line(0,1){.5108}}
\put(74.614,95.564){\line(0,1){.9365}}
\put(256.911,31.5){\line(0,1){.5147}}
\put(256.897,32.015){\line(0,1){.5132}}
\put(256.854,32.528){\line(0,1){.5101}}
\put(256.784,33.038){\line(0,1){.5055}}
\multiput(256.686,33.543)(-.03138,.12484){4}{\line(0,1){.12484}}
\multiput(256.561,34.043)(-.030528,.098352){5}{\line(0,1){.098352}}
\multiput(256.408,34.535)(-.029885,.080446){6}{\line(0,1){.080446}}
\multiput(256.229,35.017)(-.029349,.067449){7}{\line(0,1){.067449}}
\multiput(256.023,35.489)(-.032994,.065743){7}{\line(0,1){.065743}}
\multiput(255.792,35.95)(-.031973,.05586){8}{\line(0,1){.05586}}
\multiput(255.537,36.397)(-.031093,.048025){9}{\line(0,1){.048025}}
\multiput(255.257,36.829)(-.033673,.046252){9}{\line(0,1){.046252}}
\multiput(254.954,37.245)(-.032537,.039907){10}{\line(0,1){.039907}}
\multiput(254.628,37.644)(-.031519,.034607){11}{\line(0,1){.034607}}
\multiput(254.282,38.025)(-.033365,.032831){11}{\line(-1,0){.033365}}
\multiput(253.915,38.386)(-.035111,.030957){11}{\line(-1,0){.035111}}
\multiput(253.528,38.726)(-.040427,.03189){10}{\line(-1,0){.040427}}
\multiput(253.124,39.045)(-.046789,.032923){9}{\line(-1,0){.046789}}
\multiput(252.703,39.342)(-.04852,.030315){9}{\line(-1,0){.04852}}
\multiput(252.266,39.614)(-.056368,.031068){8}{\line(-1,0){.056368}}
\multiput(251.815,39.863)(-.066267,.03193){7}{\line(-1,0){.066267}}
\multiput(251.351,40.087)(-.079233,.032967){6}{\line(-1,0){.079233}}
\multiput(250.876,40.284)(-.080917,.028584){6}{\line(-1,0){.080917}}
\multiput(250.391,40.456)(-.098831,.028939){5}{\line(-1,0){.098831}}
\put(249.896,40.601){\line(-1,0){.5013}}
\put(249.395,40.718){\line(-1,0){.507}}
\put(248.888,40.808){\line(-1,0){.5112}}
\put(248.377,40.87){\line(-1,0){.5138}}
\put(247.863,40.904){\line(-1,0){.5149}}
\put(247.348,40.91){\line(-1,0){.5144}}
\put(246.834,40.887){\line(-1,0){.5124}}
\put(246.321,40.837){\line(-1,0){.5089}}
\put(245.813,40.758){\line(-1,0){.5038}}
\multiput(245.309,40.652)(-.12432,-.03339){4}{\line(-1,0){.12432}}
\multiput(244.811,40.519)(-.097847,-.03211){5}{\line(-1,0){.097847}}
\multiput(244.322,40.358)(-.079954,-.031178){6}{\line(-1,0){.079954}}
\multiput(243.842,40.171)(-.066968,-.030432){7}{\line(-1,0){.066968}}
\multiput(243.374,39.958)(-.057053,-.029793){8}{\line(-1,0){.057053}}
\multiput(242.917,39.72)(-.055338,-.032869){8}{\line(-1,0){.055338}}
\multiput(242.475,39.457)(-.047517,-.031863){9}{\line(-1,0){.047517}}
\multiput(242.047,39.17)(-.041133,-.030973){10}{\line(-1,0){.041133}}
\multiput(241.636,38.86)(-.039378,-.033177){10}{\line(-1,0){.039378}}
\multiput(241.242,38.528)(-.034095,-.032073){11}{\line(-1,0){.034095}}
\multiput(240.867,38.176)(-.032289,-.03389){11}{\line(0,-1){.03389}}
\multiput(240.512,37.803)(-.033426,-.039166){10}{\line(0,-1){.039166}}
\multiput(240.177,37.411)(-.031234,-.040935){10}{\line(0,-1){.040935}}
\multiput(239.865,37.002)(-.032165,-.047314){9}{\line(0,-1){.047314}}
\multiput(239.575,36.576)(-.03322,-.055128){8}{\line(0,-1){.055128}}
\multiput(239.31,36.135)(-.030155,-.056862){8}{\line(0,-1){.056862}}
\multiput(239.068,35.68)(-.030857,-.066773){7}{\line(0,-1){.066773}}
\multiput(238.852,35.213)(-.031685,-.079754){6}{\line(0,-1){.079754}}
\multiput(238.662,34.734)(-.032731,-.097641){5}{\line(0,-1){.097641}}
\multiput(238.499,34.246)(-.027342,-.099285){5}{\line(0,-1){.099285}}
\put(238.362,33.749){\line(0,-1){.5032}}
\put(238.253,33.246){\line(0,-1){.5084}}
\put(238.171,32.738){\line(0,-1){.5121}}
\put(238.117,32.226){\line(0,-1){2.0547}}
\put(238.184,30.171){\line(0,-1){.5076}}
\put(238.27,29.664){\line(0,-1){.5021}}
\multiput(238.384,29.161)(.02831,-.099013){5}{\line(0,-1){.099013}}
\multiput(238.526,28.666)(.033683,-.097317){5}{\line(0,-1){.097317}}
\multiput(238.694,28.18)(.032463,-.079441){6}{\line(0,-1){.079441}}
\multiput(238.889,27.703)(.031508,-.066468){7}{\line(0,-1){.066468}}
\multiput(239.11,27.238)(.030709,-.056565){8}{\line(0,-1){.056565}}
\multiput(239.355,26.785)(.030006,-.048712){9}{\line(0,-1){.048712}}
\multiput(239.625,26.347)(.032625,-.046998){9}{\line(0,-1){.046998}}
\multiput(239.919,25.924)(.031632,-.040628){10}{\line(0,-1){.040628}}
\multiput(240.235,25.518)(.030734,-.035307){11}{\line(0,-1){.035307}}
\multiput(240.573,25.129)(.032619,-.033573){11}{\line(0,-1){.033573}}
\multiput(240.932,24.76)(.034406,-.031739){11}{\line(1,0){.034406}}
\multiput(241.311,24.411)(.0397,-.03279){10}{\line(1,0){.0397}}
\multiput(241.708,24.083)(.041434,-.03057){10}{\line(1,0){.041434}}
\multiput(242.122,23.777)(.047826,-.031398){9}{\line(1,0){.047826}}
\multiput(242.553,23.495)(.055656,-.032327){8}{\line(1,0){.055656}}
\multiput(242.998,23.236)(.065532,-.033411){7}{\line(1,0){.065532}}
\multiput(243.456,23.002)(.067262,-.029777){7}{\line(1,0){.067262}}
\multiput(243.927,22.794)(.080254,-.030395){6}{\line(1,0){.080254}}
\multiput(244.409,22.611)(.098156,-.031153){5}{\line(1,0){.098156}}
\multiput(244.9,22.456)(.12464,-.03217){4}{\line(1,0){.12464}}
\put(245.398,22.327){\line(1,0){.5049}}
\put(245.903,22.226){\line(1,0){.5096}}
\put(246.413,22.152){\line(1,0){.5129}}
\put(246.926,22.107){\line(1,0){1.0294}}
\put(247.955,22.1){\line(1,0){.5134}}
\put(248.468,22.139){\line(1,0){.5105}}
\put(248.979,22.206){\line(1,0){.5061}}
\multiput(249.485,22.301)(.12504,.03058){4}{\line(1,0){.12504}}
\multiput(249.985,22.423)(.098544,.029902){5}{\line(1,0){.098544}}
\multiput(250.478,22.573)(.080634,.029373){6}{\line(1,0){.080634}}
\multiput(250.962,22.749)(.078907,.033739){6}{\line(1,0){.078907}}
\multiput(251.435,22.951)(.065952,.032575){7}{\line(1,0){.065952}}
\multiput(251.897,23.18)(.056062,.031617){8}{\line(1,0){.056062}}
\multiput(252.345,23.432)(.048222,.030787){9}{\line(1,0){.048222}}
\multiput(252.779,23.71)(.046465,.033379){9}{\line(1,0){.046465}}
\multiput(253.198,24.01)(.040113,.032283){10}{\line(1,0){.040113}}
\multiput(253.599,24.333)(.034807,.031299){11}{\line(1,0){.034807}}
\multiput(253.982,24.677)(.033043,.033156){11}{\line(0,1){.033156}}
\multiput(254.345,25.042)(.03118,.034913){11}{\line(0,1){.034913}}
\multiput(254.688,25.426)(.032146,.040223){10}{\line(0,1){.040223}}
\multiput(255.009,25.828)(.03322,.046579){9}{\line(0,1){.046579}}
\multiput(255.308,26.247)(.030623,.048326){9}{\line(0,1){.048326}}
\multiput(255.584,26.682)(.031425,.05617){8}{\line(0,1){.05617}}
\multiput(255.835,27.132)(.03235,.066062){7}{\line(0,1){.066062}}
\multiput(256.062,27.594)(.03347,.079022){6}{\line(0,1){.079022}}
\multiput(256.263,28.068)(.029098,.080734){6}{\line(0,1){.080734}}
\multiput(256.437,28.553)(.029566,.098646){5}{\line(0,1){.098646}}
\multiput(256.585,29.046)(.03016,.12514){4}{\line(0,1){.12514}}
\put(256.706,29.546){\line(0,1){.5064}}
\put(256.799,30.053){\line(0,1){.5108}}
\put(256.864,30.564){\line(0,1){.9365}}
\put(179.911,96.5){\line(0,1){.5147}}
\put(179.897,97.015){\line(0,1){.5132}}
\put(179.854,97.528){\line(0,1){.5101}}
\put(179.784,98.038){\line(0,1){.5055}}
\multiput(179.686,98.543)(-.03138,.12484){4}{\line(0,1){.12484}}
\multiput(179.561,99.043)(-.030528,.098352){5}{\line(0,1){.098352}}
\multiput(179.408,99.535)(-.029885,.080446){6}{\line(0,1){.080446}}
\multiput(179.229,100.017)(-.029349,.067449){7}{\line(0,1){.067449}}
\multiput(179.023,100.489)(-.032994,.065743){7}{\line(0,1){.065743}}
\multiput(178.792,100.95)(-.031973,.05586){8}{\line(0,1){.05586}}
\multiput(178.537,101.397)(-.031093,.048025){9}{\line(0,1){.048025}}
\multiput(178.257,101.829)(-.033673,.046252){9}{\line(0,1){.046252}}
\multiput(177.954,102.245)(-.032537,.039907){10}{\line(0,1){.039907}}
\multiput(177.628,102.644)(-.031519,.034607){11}{\line(0,1){.034607}}
\multiput(177.282,103.025)(-.033365,.032831){11}{\line(-1,0){.033365}}
\multiput(176.915,103.386)(-.035111,.030957){11}{\line(-1,0){.035111}}
\multiput(176.528,103.726)(-.040427,.03189){10}{\line(-1,0){.040427}}
\multiput(176.124,104.045)(-.046789,.032923){9}{\line(-1,0){.046789}}
\multiput(175.703,104.342)(-.04852,.030315){9}{\line(-1,0){.04852}}
\multiput(175.266,104.614)(-.056368,.031068){8}{\line(-1,0){.056368}}
\multiput(174.815,104.863)(-.066267,.03193){7}{\line(-1,0){.066267}}
\multiput(174.351,105.087)(-.079233,.032967){6}{\line(-1,0){.079233}}
\multiput(173.876,105.284)(-.080917,.028584){6}{\line(-1,0){.080917}}
\multiput(173.391,105.456)(-.098831,.028939){5}{\line(-1,0){.098831}}
\put(172.896,105.601){\line(-1,0){.5013}}
\put(172.395,105.718){\line(-1,0){.507}}
\put(171.888,105.808){\line(-1,0){.5112}}
\put(171.377,105.87){\line(-1,0){.5138}}
\put(170.863,105.904){\line(-1,0){.5149}}
\put(170.348,105.91){\line(-1,0){.5144}}
\put(169.834,105.887){\line(-1,0){.5124}}
\put(169.321,105.837){\line(-1,0){.5089}}
\put(168.813,105.758){\line(-1,0){.5038}}
\multiput(168.309,105.652)(-.12432,-.03339){4}{\line(-1,0){.12432}}
\multiput(167.811,105.519)(-.097847,-.03211){5}{\line(-1,0){.097847}}
\multiput(167.322,105.358)(-.079954,-.031178){6}{\line(-1,0){.079954}}
\multiput(166.842,105.171)(-.066968,-.030432){7}{\line(-1,0){.066968}}
\multiput(166.374,104.958)(-.057053,-.029793){8}{\line(-1,0){.057053}}
\multiput(165.917,104.72)(-.055338,-.032869){8}{\line(-1,0){.055338}}
\multiput(165.475,104.457)(-.047517,-.031863){9}{\line(-1,0){.047517}}
\multiput(165.047,104.17)(-.041133,-.030973){10}{\line(-1,0){.041133}}
\multiput(164.636,103.86)(-.039378,-.033177){10}{\line(-1,0){.039378}}
\multiput(164.242,103.528)(-.034095,-.032073){11}{\line(-1,0){.034095}}
\multiput(163.867,103.176)(-.032289,-.03389){11}{\line(0,-1){.03389}}
\multiput(163.512,102.803)(-.033426,-.039166){10}{\line(0,-1){.039166}}
\multiput(163.177,102.411)(-.031234,-.040935){10}{\line(0,-1){.040935}}
\multiput(162.865,102.002)(-.032165,-.047314){9}{\line(0,-1){.047314}}
\multiput(162.575,101.576)(-.03322,-.055128){8}{\line(0,-1){.055128}}
\multiput(162.31,101.135)(-.030155,-.056862){8}{\line(0,-1){.056862}}
\multiput(162.068,100.68)(-.030857,-.066773){7}{\line(0,-1){.066773}}
\multiput(161.852,100.213)(-.031685,-.079754){6}{\line(0,-1){.079754}}
\multiput(161.662,99.734)(-.032731,-.097641){5}{\line(0,-1){.097641}}
\multiput(161.499,99.246)(-.027342,-.099285){5}{\line(0,-1){.099285}}
\put(161.362,98.749){\line(0,-1){.5032}}
\put(161.253,98.246){\line(0,-1){.5084}}
\put(161.171,97.738){\line(0,-1){.5121}}
\put(161.117,97.226){\line(0,-1){2.0547}}
\put(161.184,95.171){\line(0,-1){.5076}}
\put(161.27,94.664){\line(0,-1){.5021}}
\multiput(161.384,94.161)(.02831,-.099013){5}{\line(0,-1){.099013}}
\multiput(161.526,93.666)(.033683,-.097317){5}{\line(0,-1){.097317}}
\multiput(161.694,93.18)(.032463,-.079441){6}{\line(0,-1){.079441}}
\multiput(161.889,92.703)(.031508,-.066468){7}{\line(0,-1){.066468}}
\multiput(162.11,92.238)(.030709,-.056565){8}{\line(0,-1){.056565}}
\multiput(162.355,91.785)(.030006,-.048712){9}{\line(0,-1){.048712}}
\multiput(162.625,91.347)(.032625,-.046998){9}{\line(0,-1){.046998}}
\multiput(162.919,90.924)(.031632,-.040628){10}{\line(0,-1){.040628}}
\multiput(163.235,90.518)(.030734,-.035307){11}{\line(0,-1){.035307}}
\multiput(163.573,90.129)(.032619,-.033573){11}{\line(0,-1){.033573}}
\multiput(163.932,89.76)(.034406,-.031739){11}{\line(1,0){.034406}}
\multiput(164.311,89.411)(.0397,-.03279){10}{\line(1,0){.0397}}
\multiput(164.708,89.083)(.041434,-.03057){10}{\line(1,0){.041434}}
\multiput(165.122,88.777)(.047826,-.031398){9}{\line(1,0){.047826}}
\multiput(165.553,88.495)(.055656,-.032327){8}{\line(1,0){.055656}}
\multiput(165.998,88.236)(.065532,-.033411){7}{\line(1,0){.065532}}
\multiput(166.456,88.002)(.067262,-.029777){7}{\line(1,0){.067262}}
\multiput(166.927,87.794)(.080254,-.030395){6}{\line(1,0){.080254}}
\multiput(167.409,87.611)(.098156,-.031153){5}{\line(1,0){.098156}}
\multiput(167.9,87.456)(.12464,-.03217){4}{\line(1,0){.12464}}
\put(168.398,87.327){\line(1,0){.5049}}
\put(168.903,87.226){\line(1,0){.5096}}
\put(169.413,87.152){\line(1,0){.5129}}
\put(169.926,87.107){\line(1,0){1.0294}}
\put(170.955,87.1){\line(1,0){.5134}}
\put(171.468,87.139){\line(1,0){.5105}}
\put(171.979,87.206){\line(1,0){.5061}}
\multiput(172.485,87.301)(.12504,.03058){4}{\line(1,0){.12504}}
\multiput(172.985,87.423)(.098544,.029902){5}{\line(1,0){.098544}}
\multiput(173.478,87.573)(.080634,.029373){6}{\line(1,0){.080634}}
\multiput(173.962,87.749)(.078907,.033739){6}{\line(1,0){.078907}}
\multiput(174.435,87.951)(.065952,.032575){7}{\line(1,0){.065952}}
\multiput(174.897,88.18)(.056062,.031617){8}{\line(1,0){.056062}}
\multiput(175.345,88.432)(.048222,.030787){9}{\line(1,0){.048222}}
\multiput(175.779,88.71)(.046465,.033379){9}{\line(1,0){.046465}}
\multiput(176.198,89.01)(.040113,.032283){10}{\line(1,0){.040113}}
\multiput(176.599,89.333)(.034807,.031299){11}{\line(1,0){.034807}}
\multiput(176.982,89.677)(.033043,.033156){11}{\line(0,1){.033156}}
\multiput(177.345,90.042)(.03118,.034913){11}{\line(0,1){.034913}}
\multiput(177.688,90.426)(.032146,.040223){10}{\line(0,1){.040223}}
\multiput(178.009,90.828)(.03322,.046579){9}{\line(0,1){.046579}}
\multiput(178.308,91.247)(.030623,.048326){9}{\line(0,1){.048326}}
\multiput(178.584,91.682)(.031425,.05617){8}{\line(0,1){.05617}}
\multiput(178.835,92.132)(.03235,.066062){7}{\line(0,1){.066062}}
\multiput(179.062,92.594)(.03347,.079022){6}{\line(0,1){.079022}}
\multiput(179.263,93.068)(.029098,.080734){6}{\line(0,1){.080734}}
\multiput(179.437,93.553)(.029566,.098646){5}{\line(0,1){.098646}}
\multiput(179.585,94.046)(.03016,.12514){4}{\line(0,1){.12514}}
\put(179.706,94.546){\line(0,1){.5064}}
\put(179.799,95.053){\line(0,1){.5108}}
\put(179.864,95.564){\line(0,1){.9365}}
\qbezier(37.5,38.5)(52.25,46)(80,45.5)
\qbezier(71.25,103.75)(86,111.25)(113.75,110.75)
\qbezier(190.75,45.75)(206.375,45.75)(239.5,36.75)
\qbezier(113.75,110.75)(129.375,110.75)(162.5,101.75)
\qbezier(38.75,37)(53.75,44)(79.75,43)
\qbezier(72.5,102.25)(87.5,109.25)(113.5,108.25)
\qbezier(190.5,43.25)(206.5,42.625)(238.5,34.5)
\qbezier(113.5,108.25)(129.5,107.625)(161.5,99.5)
\qbezier(39,26.25)(53.75,21.875)(79.5,21)
\qbezier(72.75,91.5)(87.5,87.125)(113.25,86.25)
\qbezier(190.25,21.25)(215.625,19.5)(238.5,28.75)
\qbezier(113.25,86.25)(138.625,84.5)(161.5,93.75)
\qbezier(40.5,28.5)(52.875,24.875)(79.75,23.75)
\qbezier(74.25,93.75)(86.625,90.125)(113.5,89)
\qbezier(190.5,24)(217,22.875)(237.5,30.25)
\qbezier(113.5,89)(140,87.875)(160.5,95.25)
\put(190.5,21.25){\line(0,1){.25}}
\put(113.5,86.25){\line(0,1){.25}}
\multiput(113.5,86.5)(.03125,3.09375){8}{\line(0,1){3.09375}}
\put(77.5,21.25){\line(0,1){24.5}}
\put(111.25,86.5){\line(0,1){24.5}}
\put(61.75,44.75){\line(0,-1){22.5}}
\put(95.5,110){\line(0,-1){22.5}}
\multiput(60,44.5)(.03125,-2.8125){8}{\line(0,-1){2.8125}}
\multiput(93.75,109.75)(.03125,-2.8125){8}{\line(0,-1){2.8125}}
\qbezier(80.75,42.75)(89.5,43)(99.25,37.25)
\qbezier(80,45.75)(95.75,44.125)(100.5,39)
\qbezier(79.75,24)(90.5,23.875)(97.25,28.25)
\qbezier(79.5,21)(91,20.75)(98.5,26.5)
\put(75.5,21){\line(0,1){24.75}}
\qbezier(149,38.25)(165.75,46.375)(190.5,46)
\qbezier(150.25,36.5)(169.625,44.125)(190.5,43.25)
\qbezier(151.75,30.25)(170,25.25)(190.25,24.25)
\qbezier(151,28.5)(168.375,21.75)(190.25,21)
\qbezier(121.75,32.5)(110.375,101.875)(11.5,30.75)
\qbezier(121.75,32.75)(125.625,-1.5)(59,8.25)
\qbezier(59,8.25)(40.5,10.625)(20,14.5)
\qbezier(20,14.5)(-4.25,18.875)(11.5,30.75)
\qbezier(188.25,46.25)(150.875,78.75)(58,76.25)
\qbezier(58,76.25)(7.5,64.5)(1,21.75)
\qbezier(1,21.75)(-1.5,2.625)(78,2)
\qbezier(78,2)(183.125,2)(186.75,15)
\qbezier(186.75,15)(192,41.875)(188.25,46.25)
\put(64.75,96.5){\makebox(0,0)[cc]{$\pi_1$}}
\put(169.75,96.5){\makebox(0,0)[cc]{$\pi_2$}}
\put(31.25,31.75){\makebox(0,0)[cc]{$\pi_1$}}
\put(246.75,31.5){\makebox(0,0)[cc]{$\pi_2$}}
\put(82.25,98.75){\makebox(0,0)[cc]{$\Delta_1$}}
\put(101,97.75){\makebox(0,0)[cc]{$\Delta_2$}}
\put(136,96.75){\makebox(0,0)[cc]{$\Delta_3$}}
\put(45.25,33){\makebox(0,0)[cc]{$\Delta_1$}}
\put(66.25,33.25){\makebox(0,0)[cc]{$\Delta_2$}}
\put(205.75,32.25){\makebox(0,0)[cc]{$\Delta_3$}}
\put(87.25,54){\makebox(0,0)[cc]{$\Psi(r)$}}
\put(138.25,59.5){\makebox(0,0)[cc]
{$\bar\Phi_i(r)\iv$}}
\put(90.25,94.25){\makebox(0,0)[cc]{$\aaa_1$}}
\put(118.75,93){\makebox(0,0)[cc]{$\aaa_2$}}
\put(56.5,28.75){\makebox(0,0)[cc]{$\aaa_1$}}
\put(82.5,28){\makebox(0,0)[cc]{$\aaa_2$}}
\put(140.5,111){\makebox(0,0)[cc]{$\xxx_1$}}
\put(145.75,85){\makebox(0,0)[cc]{$\xxx_2$}}
\end{picture}
\caption{}\label{fin}
\end{figure}

Suppose now that $\Delta$ has complexity $\ge 1$. Let $\aaa_1$ be the first transition $\Theta$-band in $\Delta$ counting from $\pi_1$ to $\pi_2$.

{\bf Case 1.} Let us assume that $\pi_1$ is the first hub.

Then $\aaa_1$ must be a $\bar\tau_1$-band (since the letter in $\Lab(\pi_1)$ from $\{\tilde k_{i+1}, k_{i+1}, \hat k_{i+1}\}$ is $\tilde k_{i+1}$ and relations involving $\bar\tau_2$ do not contain that letter). The side of $\aaa_1$ that is further from $\pi_1$ must contain $k_{i+1}$. Hence the next after $\aaa_1$ $\Theta$-band $\aaa_2$ in $\Delta$ (counting from $\pi_1$ to $\pi_2$ must be either a $\bar\tau_1\iv$-band or $\bar\tau_2\iv$-band. The first option is impossible because $\aaa_1$ and $\aaa_2$ intersect $\ttt_1$ in two neighbor cells, and if the first option occurs, these cells cancel which would contradict the assumption that $\Delta$ is reduced. Thus $\aaa_2$ is a $\bar\tau_2\iv$-band.

The medians of $\aaa_1$ and $\aaa_2$ cut the diagram $\Delta$ into three parts (see the top part of Figure \ref{fin}): $\Delta_1$, $\Delta_2, \Delta_3$ where $\Delta_1$ contains $\pi_1$, $\Delta_3$ contains $\pi_2$. Let $\Delta_1'$ be the diagram $\Delta_1$ without the hub cell. Note that the boundary of $\Delta_1'$ has the form $\p\q\rf\iv\mathfrak{s}\iv$ where $\p$, $\rf$ are sides of a $t_i$- and $k_{i+1}$-bands respectively, not containing transition cells, $\Lab(\q)\equiv W_0(\sss_1^{(2i+1)})$, $$\Lab(s)\equiv t_i\tilde q_1^{(2i+1)}\tilde q_2^{(2i+1)}U\tilde q_3^{(2i+1)}\ldots \tilde q_m^{(2i+1)}\tilde k_{i+1}$$ where $U$ is a word in $Z_2^{(2i+1)}$. By Lemma \ref{l:sector}, that diagram without the $t_i$-band and the $k_{i+1}$-band is combinatorially $\pp$-homotopic to a computational sector corresponding to a computation of $\sss_1^{(2i+1)}$ accepting the word $\tilde q_1^{(2i+1)}\tilde q_2^{(2i+1)}U\tilde q_3^{(2i+1)}\ldots \tilde q_m^{(2i+1)}$. By the choice of the $S$-machine $\sss_1$ then $U\equiv U(r)$ for some $r\in R$, and by Lemma \ref{mach1}, (\ref{13}), the diagram $\Delta_1$ is uniquely determined (up to $\pp'$-homotopy) by $r$.

The $k_{i+1}$-band in $\Delta_2$ consists of gluing $k_{i+1}$-cells. Indeed, no $(\Theta,k_{i+1})$-relation except transition cells involves letter $k_{i+1}$ (they involve letters $\tilde k_{i+1}$ and $\hat k_{i+1}$) and the $k_{i+1}$-band in $\Delta_2$ does not contain $\Theta$-edges by Lemma \ref{l:trans}.

Let $\Delta_2'$ be the diagram obtained from $\Delta_2$ by removing the $k_{i+1}$-band and the maximal $q_2^{(2i+1)}$-band started on $\partial\Delta_2$. By Lemma \ref{l:trans} that $k_{i+1}$-band does not contain $\Theta$-edges. Since the intersections of $\aaa_1, \aaa_2$ with the $t_i$-band $\xxx_1$ (resp. with the maximal $Q_1^{(2i+1)}$-band starting on $\partial\Delta_2$) are neighbor cells on this band, the boundary of diagram $\Delta_2'$ has the form $\q_1\rf_1\iv\mathfrak{s}\iv$ where $\Lab(\q_1)=U(r)$, $\Lab(\rf_1)$ is a word in $X^{(2i+1)}\cup Y^{(2i+1)}$, $\Lab(s)$ is a product of words of the form $\mu(a)^{(2i+1)}$, $a\in A$. By Lemma \ref{p2}, the diagram $\Delta_2'$ is determined by the word $U(r)$ up to combinatorial $\pp$-homotopy. Therefore the diagram $\Delta_2$ is determined by $r$ up to combinatorial $\pp$-homotopy (and is homotopic to $\Upsilon(\mu(r))$).

{\bf Case 2.} Now assume that $\pi_1$ is the second hub.

Then $\aaa_1$ must be a $\bar\tau_2$-band (since the letter in $\Lab(\pi_1)$ from $\{\tilde k_{i+1}, k_{i+1}, \hat k_{i+1}\}$ is $\hat k_{i+1}$ and relations involving $\bar\tau_1$ do not contain that letter). The side of $\aaa_1$ that is further from $\pi_1$ must contain $k_{i+1}$. Similarly to the previous case, the next after $\aaa_1$ $\Theta$-band $\aaa_2$ in $\Delta$ (counting from $\pi_1$ to $\pi_2$) is $\bar\tau_1\iv$-band.

As in Case 1, the medians of $\aaa_1$ and $\aaa_2$ cut the diagram $\Delta$ into three parts $\Delta_1$, $\Delta_2, \Delta_3$ where $\Delta_1$ contains $\pi_1$, $\Delta_3$ contains $\pi_2$. Let $\Delta_1'$ be the diagram $\Delta_1$ without the hub cell. Note that the boundary of $\Delta_1'$ has the form $\p\q\rf\iv\mathfrak{s}\iv$ where $\p$, $\rf$ are sides of a $t_i$- and $k_{i+1}$-bands respectively, not containing transition cells, $\Lab(\q)\equiv W_0(\sss_1^{(2i+1)})$, $$\Lab(s)\equiv t_i\tilde q_1^{(2i+1)}\tilde q_2^{(2i+1)}U\tilde q_3^{(2i+1)}\ldots \tilde q_m^{(2i+1)}\tilde k_{i+1}$$ where $U$ is a word in $Z_2^{(2i+1)}$. By Lemma \ref{l:sector1}, that diagram without the $t_i$-band and the $k_{i+1}$-band is combinatorially $\pp$-homotopic to a computational sector corresponding to a computation of $\sss_2^{(2i+1)}$ accepting the word $\hat q_1^{(2i+1)}U \hat q_2^{(2i+1)}\hat q_3^{(2i+1)}\ldots \hat q_m^{(2i+1)}$. By the choice of the $S$-machine $\sss_2$ then $U\equiv \e(r)$ for some $r\in R$, and by Lemma \ref{mach1}, (\ref{13}), the diagram $\Delta_1$ is uniquely determined (up to $\pp'$-homotopy) by $r$.

As in Case 1, the $k_{i+1}$-band in $\Delta_2$ consists of gluing $k_{i+1}$-cells.
Let $\Delta_2'$ be the diagram obtained from $\Delta_2$ by removing the $k_{i+1}$-band and the maximal $q_2^{(2i+1)}$-band started on $\partial\Delta_2$. By Lemma \ref{l:trans} that $k_{i+1}$-band does not contain $\Theta$-edges. Since the intersections of $\aaa_1, \aaa_2$ with the $t_i$-band $\xxx_1$ (resp. with the maximal $Q_1^{(2i+1)}$-band starting on $\partial\Delta_2$) are neighbor cells on this band, the boundary of diagram $\Delta_2'$ has the form $\q_1\rf_1\iv\mathfrak{s}\iv$ where $\Lab(\q_1)=U'$, $\Lab(\rf_1)$ is a word in $X^{(2i+1)}\cup Y^{(2i+1)}$, $\Lab(s)=\mu(r)$. By Lemma \ref{p2}, the diagram $\Delta_2'$ is determined by the word $\mu(r)$ up to combinatorial $\pp$-homotopy. Therefore the diagram $\Delta_2$ is determined by $r$ up to combinatorial $\pp$-homotopy (and is homotopic to $\Upsilon(\mu(r))$).

Thus in both cases we have proved that up to combinatorial homotopy the subdiagram $\Delta''=\Delta_1\cup\aaa_1\cup\Delta_2\cup\aaa_2$ of $\Delta$ coincides with the corresponding subdiagram $\Phi_i(r)$ of the diagram $\Psi(r)$ described in Notation \ref{notphi}.

Consider the following surgery (see the bottom part of Figure \ref{fin}). We cut $\Delta$ along the boundary of $\Delta''$, insert in the resulting hole the diagram $\bar\Phi_i(r)$ (the complement of $\Phi_i(r)$ in $\Psi(r)$, see Notation \ref{notphi}) and its mirror image $\bar\Phi_1(r)\iv$ so that $\Delta''\cup \bar\Phi_i(r)$ is combinatorially $\pp$-homotopic to $\Psi(r)$, and the subdiagram $\bar\Phi_i(r)$ and its mirror image $\bar\Phi_i(r)\iv$ cancel each other. The new diagram $\check\Delta$, after a combinatorial $\pp$-homotopy becomes a union of a subdiagram that is a copy of $\Psi(r)$, a $2$-hub reduced $i$-subdiagram of smaller complexity than $\Delta$, and a number of non-hub cells corresponding to relations of $\pp$. This allows us to proceed by induction on complexity of the 2-hub $i$-subdiagram.
\endproof

\begin{lemma}\label{l:final} Let $\Delta$ be any reduced diagram over $\pp$ whose boundary label is a word in $\e(A)$. Then $\Delta$ is combinatorially $\pp$-homotopic to a normal diagram that coincides with the union of its standard subdiagrams.
\end{lemma}

\proof Every normal diagram over $\pp$ can be viewed as a diagram over $\pp'$: ignore the cells of standard subdiagrams and view these subdiagrams as $\Gamma$-cells (the standard diagram $\Psi(r)$ has boundary label $\e(r)$ by Notation \ref{notphi}). Therefore if $\Delta$ does not have hubs, we can apply Lemma \ref{l:e}. The combinatorial $\pp$-homotopy from Lemma \ref{l:e} does not touch the $\Gamma$-cells (since $\pp$ does not contain $\Gamma$-relations), so $\Delta$ is combinatorially $\pp$-homotopic to a normal diagram.

Suppose that $\Delta$ has hubs. Consider the graph $\nabla=\nabla(\Delta)$ whose vertices are the hubs of $\Delta$, two hubs are adjacent is there exists a $t_i$-band, $i<N$, connecting them. Since no vertex is adjacent to itself and the graph is planar ($t_i$-bands do not intersect), we apply Heawood's theorem \cite{He} again, and conclude that $\nabla$ has a vertex of degree $5$. Since $\partial\Delta$ does not have $t_i$-edges, there exist 2 hubs $\pi_1, \pi_2$ connected by a $t_i$-band and a $t_j$-band, $i<j<N$ (recall that $N\ge 12$). Consider the subdiagram $\Delta'$ of $\Delta$ bounded by the two hubs $\pi_1, \pi_2$ and the $t_i$- and $t_j$-band. Assume that $\Delta'$ is a smallest under inclusion such subdiagram of $\Delta$. Suppose that $\Delta'$ contains hubs. Then the graph $\nabla(\Delta')$ contains a vertex of degree 5. By the minimality of $\Delta'$, no two hubs of $\Delta'$ are connected by a $t_{i'}$- and $t_{j'}$-bands, $i'<j'<N$. Hence there exists a hub $\pi_3\in \Delta'$ which has at least 3 of its $t_j$-bands, $j<N$, each connecting that hub with $\pi_1$ or $\pi_2$ (by the pigeon-hole property, since $11>2\times 5$). Then either $\pi_1$ or $\pi_2$ is connected with $\pi_3$ by two of these bands which contradicts the minimality of $\Delta'$. Hence we can assume that $\Delta'$ does not have hubs. Then the $k_{i+1}$-band starting on $\partial\pi_1$ ends on $\partial\pi_2$. Consider the subdiagram $\Delta''$ bounded by $\pi_1, \pi_2$, the $t_i$-band and the $k_{i+1}$-band connecting them. Viewed as a diagram over $\pp'$ (we again consider standard subdiagrams as $\Gamma$-cells) $\Delta''$ is a 2-hub $i$-subdiagram of $\Delta$. By Lemma \ref{l:e4}, this subdiagram is combinatorially $\pp$-homotopic to a normal diagram. Replacing $\Delta''$ by that normal diagram we obtain a diagram with the same boundary label as $\Delta$ but with hubs outside standard subdiagrams (note that this homotopy does not touch cells outside $\Delta''$), so we can proceed by induction on the number of hubs outside standard subdiagrams.
\endproof

The following two lemmas {\em complete the proof of Theorem \ref{th:1}}.

\begin{lemma}[Part (2) of Theorem \ref{th:1}]\label{l:e9} The map $\e\colon \Gamma\to \da$ is injective.
\end{lemma}

\proof Suppose that for some word $u$ in $A$, $\e(u)=1$ in $\da$. Then there exists a \vk diagram $\Delta$ over $\pp$ with boundary label $\e(u)$. By Lemma \ref{l:final}, $\Delta$ is combinatorially homotopic to a normal diagram tesselated by standard subdiagrams. If we view each of these subdiagrams as a $\Gamma$-cell, $\Delta$ becomes a diagram $\Delta_1$ over the presentation $\la \e(A)\mid \e(R)\ra$. Applying $\e\iv$ to the labels of $\Delta_1$, we obtain a \vk diagram over the presentation of $\Gamma$ with boundary label $u$. Hence $u=1$ in $\Gamma$.
\endproof

\begin{lemma}[Part (3) of Theorem \ref{th:1}]\label{l:e10} If the presentation $\la A\mid R\ra$ of $\Gamma$ is combinatorially aspherical, then the presentation $\pp$ of $\da$ is combinatorially aspherical.
\end{lemma}

\proof Consider a spherical diagram $\Delta$ over $\pp$. It is a disc diagram with empty boundary. By Lemma \ref{l:final}, $\Delta$ is combinatorially $\pp$-homotopic to a spherical normal diagram $\Delta_1$ tesselated by standard subdiagrams. Again, view $\Delta_1$ as a spherical diagram over the presentation of $\Gamma$. Since that presentation is aspherical by our assumption, $\Delta_1$ can be combinatorially deformed to a trivial diagram by diamond moves and insertions and deletions of cancelable standard subdiagrams. Therefore $\Delta$ is combinatorially homotopic to a trivial diagram.
\endproof

\begin{rk}\label{rk510} It is quite possible to show, using the ideas from  \cite{BORS}, that the embedding constructed in this paper is quasi-isometric (i.e. the subgroup $\e(\Gamma)$ embeds without distortion), and preserves solvability of the word problem. But the Dehn function of the finitely presented group $\da$ is almost always superexponential:
indeed, by Remark \ref{rk:1}, the computational sector of $\sss_1$ for the computation accepting $U(r)$ has superexponential area in terms of $|r|$, and is a subdiagram in the \vk diagram for the relation $\e(r)=1$ in $\da$. Another property that this construction most probably preserves is the {\em finite decomposition complexity} of Guentner, Tessera and Yu \cite{GTY}. Finally D. Osin asked whether one can embed every recursively presented group into a finitely presented group as  a malnormal subgroup. Again it is quite possible that $\e(\Gamma)$ embeds malnormally, and so the answer to that question is affirmative. In order to check malnormality of $\e(\Gamma)$, one needs to consider an annular diagram $\Delta$  with boundary labels $\e(u)$, $\e(v)$ where $u,v$ are words over $A$ and prove that the label of every path $\p$ in $\Delta$ connecting the inner and the outer boundaries is equal to a word $w$ over $\e(A)$ in $\da$. For this, one needs to check that various subdiagrams considered in the proof of Lemma \ref{l:final} and the previous lemmas do not contain the hole of this diagram. If so, Lemma \ref{l:final} would apply, and we would conclude that the annular diagram $\Delta$ is combinatorially homotopic to an annular diagram tesselated by subdiagrams of the form $\Psi(r), r\in R$. A combinatorial homotopy transforms $\p$ to a path whose label is equal to $\Lab(\p)$ in $\da$. This would imply that indeed $\Lab(\p)\in \e(\Gamma)$. We have not checked the details and leave it to the reader.
\end{rk}

\end{document}